\documentclass[a4paper,fleqn]{cas-sc}

\usepackage[numbers]{natbib}

\usepackage{preamble_W}

\newcommand{\setOfReals}{\mathbb{R}}
\newcommand{\setOfNaturals}{\mathbb{N}}
\newcommand{\setOfNonnegativeIntegers}{{\mathbb{N}_0}}
\newcommand{\setOfPositiveReals}{{\setOfReals_{+}}}

\newcommand{\BigO}[1]{\mathop{}\!O{\left(#1\right)}}

\newcommand{\sys}[1]{\textsc{#1}}

\newcommand{\absolute}[1]{| #1 | }

\makeatletter
\let\oldabs\abs
\def\abs{\@ifstar{\oldabs}{\oldabs*}}

\newcommand{\cadlag}{c\`adl\`ag }

\newcommand{\defeq}{\coloneqq}

\newcommand{\indicator}[2]{\mathsf{1}_{#1}{\left(#2\right)}}
\newcommand{\indic}[1]{\mathsf{1}_{#1}}

\newcommand{\differential}[1]{\mathrm{d} #1}
\newcommand{\timeDerivative}[1]{\frac{\differential}{\differential t} #1 }

\newcommand{\eqstop}{.}
\newcommand{\eqcomma}{,}

\newcommand{\norm}[1]{\left\lVert#1\right\rVert}

\newcommand{\E}{\mathbb{E}}
\newcommand{\Eof}[1]{\E\left[#1 \right]}

\newcommand{\prob}{\mathbb{P}}
\newcommand{\probOf}[1]{\prob\left(#1\right)}

\newcommand{\history}[1]{\mathcal{F}_{#1}  }
\newcommand{\trans}{^{\mathsf T}}

\newcommand{\myExp}[1]{\exp \left( #1 \right)  }

\newcommand{\ie}{\textit{i.e.}}
\newcommand{\eg}{\textit{e.g.}}

\newcommand{\f}{\frac}
\newcommand{\Rt}{\longrightarrow}
\newcommand{\rt}{\rightarrow}
\newcommand{\RT}{\Rightarrow}
\newcommand{\LRT}{\Longrightarrow}

\newcommand{\nrt}{\stackrel{n\rt \infty}\Rt}

\newcommand{\mrm}{\mathrm}

\newcommand{\leb}{\l_{\mbox{Leb}}}

\newcommand{\np}{\noindent}

\newcommand{\hs}{\hspace}

\newcommand{\om}{\omega}
\newcommand{\Om}{\Omega}
\newcommand{\g}{\gamma}

\newcommand{\ep}{\epsilon}

\def\SC{\mathcal}

\newcommand{\mfk}{\mathfrak}

\def \triple|{|\! | \! |}
\renewcommand\leq{\ensuremath{\leqslant}}
\def\lf{\left}
\def\ri{\right}
\renewcommand\geq{\ensuremath{\geqslant}}

\DeclareMathOperator*{\argmax}{arg\,max}

\def\l{\lambda}

\def\ot{\otimes}
\def\<{\langle}
\def\>{\rangle}
\def\~{\tilde}
\newcommand{\EE}{\mathbb{E}}
\newcommand{\PP}{\mathbb{P}}

\def\N{\mathbb N}

\def\BddZ{\mathbb B}
\def\R{\mathbb R}

\newcommand{\gen}{\mathcal{A}}
\newcommand{\fgen}{\mathcal{B}}

\newcommand{\occ}{\Gamma}
\newcommand{\mart}{\SC{M}}
\newcommand{\modu}{\mfk{m}}

\newcommand{\like}{\mathscr{L}}

\newcommand{\subst}{Y}   
\newcommand{\tauden}{\phi} 
\newcommand{\taudist}{\Phi} 

\newcommand{\conscon}{J}
\newcommand{\const}{\mathscr{C}}
\newcommand{\cnst}{\mathscr{B}}

\begin{document}
\let\WriteBookmarks\relax
\def\floatpagepagefraction{1}
\def\textpagefraction{.001}

\shorttitle{Statistical inference for Michaelis--Menten enzyme kinetics}    

\shortauthors{Ganguly and KhudaBukhsh}  

\title [mode = title]{Statistical inference for a multiscale stochastic model of enzyme kinetics via propagation of chaos}  

\tnotemark[1,2] 

\tnotetext[1]{Research of A. Ganguly is supported in part by NSF DMS - 1855788, NSF DMS - 2246815 and Simons Foundation (via Travel Support for Mathematicians). W. R. KhudaBukhsh was supported in part by an International Research Collaboration Fund (IRCF) awarded by the University of Nottingham, United Kingdom,  a Scheme 4 `Research in Pairs' Grant 42360 by the London Mathematical Society (LMS), and in part by the Engineering and Physical Sciences Research Council (EPSRC) [grant number EP/Y027795/1].} 

\tnotetext[2]{W. R. KhudaBukhsh thanks the Isaac Newton Institute for Mathematical Sciences (INI) (EPSRC [grant number EP/R014604/1]) for the support and the hospitality during an INI Retreat programme when a  part of the work was undertaken. 
}

%

\author[1]{Arnab Ganguly}

\cormark[1]


\ead{aganguly@lsu.edu}

\ead[url]{https://sites.google.com/site/gangulymath/}

\credit{Methodology, Writing}

\affiliation[1]{organization={Department of Mathematics, Louisiana State University},
            addressline={303 Lockett Hall}, 
            city={Baton Rouge},
            postcode={70803-4918}, 
            state={Louisiana},
            country={United States of America}}

\author[2]{Wasiur R. KhudaBukhsh}[orcid=0000-0003-1803-0470]

\cormark[1]


\ead{wasiur.khudabukhsh@nottingham.ac.uk}

\ead[url]{https://www.wasiur.xyz/}

\credit{Methodology, Writing}

\affiliation[2]{organization={School of Mathematical Sciences, University of Nottingham},
            addressline={University Park}, 
            city={Nottingham},
            postcode={NG7 2RD}, 
            state={Nottinghamshire},
            country={United Kingdom}}

\cortext[1]{Corresponding author}



\begin{abstract}
We study a class of \acp{SDE} with jumps modeling multistage Michaelis--Menten enzyme kinetics, in which a substrate is sequentially transformed into a product via a cascade of intermediate complexes. These networks are typically high-dimensional and exhibit multiscale behavior with a strong coupling between different components, posing substantial analytical and computational challenges. In particular, the problem of statistical inference of reaction rates is significantly difficult, and becomes even more intricate when direct observations of system states are unavailable and only a random sample of product formation times is observed.
We address this problem in two stages. First, in a suitable scaling regime consistent with the \ac{QSSA}, we rigorously establish a stochastic averaging principle yielding a reduced model for the product-substrate dynamics. Guided by the reduced-order dynamics, we next construct a novel \ac{IPS} that approximates the product-substrate process at the particle level. This \ac{IPS} plays a pivotal role in the inference methodology, and we prove several non-asymptotic bounds and limiting results for this system. These results facilitate the construction of an estimator based on a  product-form approximate likelihood requiring only a random sample of product formation times. This approach does not need access to the system states, and we rigorously prove consistency of the estimator.  
\end{abstract}




\maketitle


\section{Introduction}
\label{sec:intro}
    Biochemical reactions can occur over an extremely wide range of time scales. Quoting Wolfenden and Snider \cite{Wolfenden2001Depth}, ``\emph{The fastest known reactions include reactions catalyzed by enzymes, but the rate enhancements that enzymes produce had not been fully appreciated until recently. In the absence of enzymes, these same reactions are among the slowest that have ever been measured, some with half-times approaching the age of the Earth.}'' It implies that some reaction rates can never be estimated using traditional methods. While this poses a challenge to the experimentalists, it also presents a unique opportunity to the mathematicians. Indeed, a variety of techniques, and postulates have been developed to deal with the multiscale nature of biochemical reactions. Michaelis and Menten studied the reaction of invertase and derived what became known as the \acf{MM} equation \cite{Michaelis1913MM}. 

    In the simplest form, the \ac{MM} enzyme kinetic reactions (\cite{Segel:1975:EK,Cornish-Bowden:2004:FEK}) consist of a reversible binding of a substrate and an enzyme into a substrate-enzyme complex, and  the conversion of the substrate-enzyme complex into a product freeing up the bound enzyme. Schematically, the \ac{MM} reaction network is represented as follows
    \begin{align}
    S + E \xrightleftharpoons[k_{-1}]{k_1} C \xrightarrow[]{k_P}    P+ E,
    \label{eq:mm_standard_det}
    \end{align}
    where $S, E$, and $P$ denote molecules of the substrate, the enzyme, and the product while $C$ denotes the intermediate complex. The nonnegative real numbers $k_1$, $k_{-1}$, and $k_P$ are reaction rate constants. Under the law of mass-action \cite{Anderson2015book}, the time evolution of the concentration of molecules of $S, E, C$, and $P$ can be described using the following system of \acp{ODE}:
\begin{align}
\begin{split}
    \timeDerivative{S} &{} = -k_1 S E + k_{-1} C \eqcomma \quad
    \timeDerivative{E} = -k_1 S E + k_{-1} C + k_P C \eqcomma \\
    \timeDerivative{C} &{} = k_1 S E - (k_P + k_{-1}) C \eqcomma \quad
    \timeDerivative{P} = k_P C \eqcomma
\end{split}
\label{eq:mm_standard_det_ode}
\end{align}
with initial conditions $S(0) = s_0, E(0)=e_0,  P(0)=0$, and $C(0)=0$. Based on empirical data, it is argued that the intermediate complex $C$ reaches a steady-state, \ie, $\timeDerivative{C}\approx 0$ very quickly while the species $S, E$, and $P$ remain in their transient states. Therefore, heuristically setting $\timeDerivative{C} = 0$ and using the conservation law $e_0 = E+C$, we get the steady-state value $C = e_0 S / (k_M+S),$
    where $k_M=(k_{-1}+k_P)/k_1$ is known as the \acl{MM} constant. The substrate concentration is then given by the  \ac{ODE}
    \begin{align}
    \timeDerivative{S} &= -\frac{k_P e_0 S}{k_M+S}\eqstop   \label{eq:det_sQSSA_S}
    \end{align}
    This approximation is known as the  \ac{sQSSA} for the deterministic \acl{MM} enzyme kinetic reaction system in  \eqref{eq:mm_standard_det}. Roughly speaking, the \ac{sQSSA} approximates the \ac{CRN} in \eqref{eq:mm_standard_det} by a single reaction of the form: $S \longrightarrow P$. 

    The original motivation of Michaelis and Menten, however, was to express the quantity $k_P C$, known as the velocity of the reaction  (the rate of change of the complex in the direction of product formation), in terms of a known quantity such as the substrate concentration $S$ rather than the unknown and unobservable complex concentration $C$. 
    The validity of this approximation, its generalizations, and its apparent misuse (the utilization of \ac{sQSSA} even when the amount of substrate is not sufficiently large compared to the amount of enzymes) have been studied extensively in the applied mathematics literature \cite{Segel:1988:VSS,Segel:1989:QSS,Stiefenhofer:1998:QSS,Tzafriri:2007:QSS,Choi2017BeyondMM,Kim2020Misuse,Schnell2013Validity,Schnell2000Enzymekinetics,Kang2019QSSA}.

    Without going into detailed biochemistry, which is outside the scope of the present paper, and can be found in standard textbooks \cite{Cornish-Bowden:2004:FEK}, we note that the simple description of the \ac{MM} reaction network in \eqref{eq:mm_standard_det} is a significant abstraction of the actual reaction network. Indeed, the \ac{MM} reaction network in reality consists of a reversible binding of a substrate and an enzyme into a substrate-enzyme complex, which, in turn, reversibly undergoes several intermediate stages, to eventually produce a product freeing up the bound enzyme. Thus, a more realistic \ac{MM} reaction network can be schematically represented as follows
    \begin{equation}
    S + E \xrightleftharpoons[k_{-1}]{k_1} C_1 \xrightleftharpoons[k_{-2}]{k_2}  C_2 \xrightleftharpoons[k_{-3}]{k_3} \cdots \xrightleftharpoons[k_{-r+1}]{k_{r-1}}C_{r-1} \xrightleftharpoons[k_{-r}]{k_{r}} C_r\xrightharpoonup[]{k_P}    P+ E,
    \label{eq:mm_det}
    \end{equation}
    where $S, E$, and $P$ denote molecules of the substrate, the enzyme, and the product while $C_1, C_2, \ldots$, and $ C_r$ denote the $r \in \setOfNaturals$ intermediate complexes. The nonnegative real numbers $k_i$, $k_{-i}$, for $i=1, 2, \ldots, r$ and $k_P$ are reaction rate constants. We refer the readers to \cite{Srinivasan2021MM} for an accessible discussion on multi-stage \ac{MM} reaction networks such as the one in \eqref{eq:mm_det}. 

    \subsection{Our contributions}
  This paper considers stochastic modeling of the multi-stage \ac{MM} system described in \eqref{eq:mm_det} via a \ac{CTMC}, whose sample paths are represented by a system of strongly coupled \acp{SDE} driven by \acp{PRM} (see \eqref{eq:X_n_traj}). The analysis of this system is mathematically challenging due to its high dimensionality, interaction among various components and the multiscale nature, where certain processes such as enzyme-substrate binding, dissociation, and intermediate complex formation occur at significantly faster time scales than product formation. An important problem in this context is the statistical inference for \ac{MM} systems, made difficult not only by high dimensionality but also a lack of data on fast, short-lived intermediate complexes. In practical scenarios, the situation is actually further complicated by the fact that even the amount of product ($P$) over time may not be observed. Instead, one often has access only to a limited dataset consisting of random samples of product formation times. 
  The core obstacle here is  {\em mathematical rather than statistical}: without observations of the system's internal states at any time point, one cannot even write down a likelihood function, and the primary goal of this paper is to develop a novel approach to this important inference problem. 
  
  We address this in two stages --- first, by studying stochastic averaging in a suitable scaling regime, which gives a reduced order \ac{ODE} model and provides a rigorous justification for the heuristic derivation of \acp{sQSSA}, and then by constructing an \ac{IPS} based on the simpler limiting model to capture particle-level dynamics.

    
    \paragraph{Stochastic averaging:} The first goal of this paper is to rigorously derive the \ac{sQSSA} for the \ac{MM} system in \eqref{eq:mm_det} directly from a stochastic \ac{CTMC}-based model. Specifically, we show that the \ac{sQSSA} follows as a consequence of the \ac{FLLN} (\Cref{thm:FLLN}) in an appropriate scaling regime. A  mathematical derivation of the \ac{sQSSA} for the simple \ac{MM} model in \eqref{eq:mm_standard_det} has been given in \cite{Kang:2013:STM}
(see also \cite{Ball:2006:AAM,Kang2019QSSA,Enger2023Unified}). While the high-level strategy of proving an averaging principle is well known and can be summarized in a few steps, and has been elegantly formalized in \cite{Kurtz1992Averaging, Kang:2013:STM} (see also \cite{Khas68,Skorokhod1989Asymptotic,Veretennikov1990Averaging,Lip96,PaVe01,PaVe02,FrWe08,FrWe12} for averaging principles in generic settings of continuous diffusion processes), establishing an averaging principle for specific complex systems still requires involved analysis. In our case, some of the tightness arguments simplify because of the conservation laws satisfied by \eqref{eq:mm_det}, but a rigorous derivation of the closed-form limiting dynamics nevertheless involves subtle arguments and careful calculations.

    \paragraph{Parameter inference via propagation of chaos:}

    
    We next come to the primary objective of this paper: developing a mathematical framework for estimating key parameters of the system from data consisting of a sample of {\em product formation of times}. The derivation of the reduced order model, an important result in its own right, is a crucial step toward this goal since it simplifies the original system in two major ways. First, it reduces the process complexity by demonstrating that the full system can be approximated by a simplified system of the form
    $$S \xrightharpoonup{\hphantom{\text{long}}} P,$$ 
    and second, it shows that the effective parameter $\theta$  driving the reduced-order \ac{ODE} lives in a lower-dimensional parameter space. However, the absence of temporal data on the amount of species of the system means that traditional estimation methods like those used in \cite{Wilkinson2018SMS,Choi2011Inference,Burnham2008inference,Ocal2019parameter} cannot be used for estimating the relevant parameters even from the simplified model. Indeed, note that  a dataset consisting only of a random sample of times of product formation does not allow reconstruction of trajectories of species copy numbers rendering any trajectory-based inference method categorically inapplicable even if the method allows partially observed trajectories. As noted earlier, the core difficulty is conceptual rather than technological, stemming from the lack of a clear formulation of the likelihood function suitable for the data type considered in the paper. This requires a completely new way of thinking about parameter inference in this case. The crux of our parameter inference method lies in the shift of our focus from \emph{population counts} to \emph{times of conversion of individual molecules}. This shift in perspective has multiple advantages, as evidenced in the \ac{DSA} approach in infectious disease epidemiology \cite{KhudaBukhsh2020Focus,KhudaBukhsh2024Howto}.
    
    
    To this end, taking cues from the reduced-order model, we first construct a weakly \acl{IPS} $\{\subst^{(n)}_i : i=1, \dots, n\}$ taking values in $\{0, 1\}$, where $1$ and $0$ represent the $S$ and the $P$ states, respectively. This system approximately captures the dynamics of the multistage \ac{MM} reaction network at the molecular level. We prove several important non-asymptotic bounds for this \ac{IPS}, including a concentration inequality in \Cref{th:conc-ineq}, the $L^p$ rate of convergence to the limiting \ac{ODE}, and a quantitative rate of convergence for the \emph{propagation of chaos} with respect to the Kolmogorov distance on bounded intervals (\Cref{cor:prop-chaos-tau}). 
    
    We next model the data, consisting of a random sample of product formation times within a fixed interval $[0,T]$, as a random subset of the conversion times of $\subst^{(n)}_i$ from state $1$ to state $0$ occurring before time $T$. While the joint distribution of these times is analytically intractable, our propagation of chaos result justifies a product-form approximation for large $n$. This facilitates the estimation of the effective parameter $\theta$  via an approximate likelihood function that bypasses the need for any state data, relying instead only on the observed product formation times (see \Cref{sec:like}). We rigorously establish the consistency of this estimator in \Cref{sec:const}.

    \subsection{Outline and notations} \label{sec:out-not}

    The rest of the paper is structured as follows. In \Cref{sec:stoch_model}, we describe the \ac{CTMC} model of the multi-stage \ac{MM} enzyme kinetic reaction network. In \Cref{sec:sQSSA}, we describe the \ac{sQSSA} and prove the \ac{FLLN} for the scaled process. In \Cref{sec:IPS}, we construct the \ac{IPS} and provide the necessary limit theorems for the propagation of chaos and the parameter inference. We make some concluding remarks in \Cref{sec:conclusion}. 
    Additional mathematical derivations are provided in the appendices. 

    The sets of natural numbers and real numbers are denoted by $\setOfNaturals$, and $\setOfReals$ respectively. We use $\setOfNonnegativeIntegers$ to denote the set of nonnegative integers. The set $\R^{d\times d'}$ will denote the space of all $d\times d'$ real matrices. We will  denote the canonical unit vector in $\R^d$ whose $i$-th component is 1 and all other components are 0 by $e^{(d)}_i$. For a column vector $v \in \R^d$, we write $v^{\ot 2} \defeq vv^\top \in \R^{d\times d}.$
    We will use $\norm{\cdot}_1$ and $\norm{\cdot}_{\infty}$ to denote $l_1$ and $l_\infty$ norms. 
     
    The spaces of continuous functions and Lipschitz continuous functions from a metric space $E$ to a metric space $F$ are denoted by $C(E, F)$ and $\mrm{Lip}(E,F)$, respectively. The space $C([0,T], F)$ is a subset of the space $D([0,T], F)$, the space of \cadlag functions from  $[0,T]$ to a metric space $F$, \ie, functions that are right continuous and have finite left-hand limits. For a function $x$, we will denote the left-hand limit of $x$ at $t$ by $x(t-)$. The symbol $\indicator{A}{\cdot}$ denotes the indicator function of a set $A$. 
    
    For a measurable space $(E, \mathcal{E})$, the symbols $\mathcal{P}(E)$ and $\mathcal{M}(E)$ denote the space of probability measures and finite non-negative measures on $E$, respectively. If $\mu$ and $\nu$ are two probability measures on $E$ with $\mu \ll \nu$, the {\em relative entropy} of $\mu$ with respect to $\nu$ is defined as
    $$\mrm{RE}(\mu \| \nu) \defeq \int_E \ln \lf(\f{\differential{\mu}}{\differential{\nu}}\ri)\differential{\mu}.$$
    Suppose for a reference measure $\lambda$ on $E$, $\mu \ll \lambda$ and $\nu \ll \lambda$.
    Then, the {\em entropy} of $\mu$ and the {\em cross entropy} of $\mu$ relative to $\nu$ are defined as
    \begin{align*}
        \mrm{Ent}(\mu) \defeq -\int_E \f{\differential{\mu}}{\differential{\lambda}} \ln\lf(\f{\differential{\mu}}{\differential{\lambda}} \ri) \differential{\lambda}, \qquad 
        \mrm{CE}(\mu, \nu) \defeq -\int_E \ln\lf(\f{\differential{\nu}}{\differential{\lambda}}\ri)\f{\differential{\mu}}{\differential{\lambda}} \differential{\lambda}.
    \end{align*}
    Clearly, $\mrm{CE}(\mu, \nu) = \mrm{Ent}(\mu)+ \mrm{RE}(\mu \| \nu).$ When $E = \R^d,$ we typically take $\lambda$ as the Lebesgue measure, denoted by $\leb$.



    \section{Stochastic model}
    \label{sec:stoch_model}
    In this section, we describe the stochastic model of the multi-stage \ac{MM} enzyme kinetic reaction network in \eqref{eq:mm_det} in terms of a \ac{CTMC}. 
    We will use the symbol $n$ as a scaling parameter that will encode reaction speeds and species abundance. 
    Let $\kappa_i^{(n)}, \kappa_{-i}^{(n)}$, for $i=1, 2, \ldots, r$, and $\kappa_P^{(n)} $
    be the (stochastic) reaction rate constants.
    Denote by $X_S^{(n)}(t)$, $X_E^{(n)}(t)$, $X_{C, 1}^{(n)}(t)$, $X_{C, 2}^{(n)}(t), \ldots, X_{C, r}^{(n)}(t)$,  and $X_P^{(n)}(t)$ the species copy numbers of  $S, E, C_1, C_2, \ldots, C_r$, and $P$ at time $t\ge 0$. The dynamics of the system can be described by the following system of \acp{SDE} driven by \acp{PRM} 
     \begin{align}
        \begin{aligned}
            X_E^{(n)}(t) &{} = X_E^{(n)}(0) - \int_{[0,t]\times [0,\infty)}\indicator{[0, \lambda_1^{(n)}(X^{(n)} (s-))]}{v} \xi_1(\differential{v}\times \differential{s}) 
            + \int_{[0,t]\times [0,\infty)}\indicator{[0, \lambda_{-1}^{(n)}(X^{(n)} (s-))]}{v} \xi_{-1}(\differential{v}\times \differential{s}) \\
            &{} \quad + \int_{[0,t]\times [0,\infty)}\indicator{[0, \lambda_{P}^{(n)}(X^{(n)} (s-))]}{v} \xi_{P}(\differential{v}\times \differential{s}) \eqcomma \\
            X_{C, i}^{(n)}(t) &{} = X_{C, i}^{(n)}(0) + \int_{[0,t]\times [0,\infty)}\indicator{[0, \lambda_i^{(n)}(X^{(n)} (s-))]}{v} \xi_i(\differential{v}\times \differential{s}) 
            + \int_{[0,t]\times [0,\infty)} \indicator{[0, \lambda_{-(i+1)}^{(n)}(X^{(n)} (s-))]}{v} \xi_{-(i+1)}(\differential{v}\times \differential{s}) \\
            &{} \quad 
            - \int_{[0,t]\times [0,\infty)}\indicator{[0, \lambda_{-i}^{(n)}(X^{(n)} (s-))]}{v} \xi_{-i}(\differential{v}\times \differential{s}) 
            -  \int_{[0,t]\times [0,\infty)}\indicator{[0, \lambda_{i+1}^{(n)}(X^{(n)} (s-))]}{v} \xi_{i+1}(\differential{v}\times \differential{s})\eqcomma\\
            & \hs{.3cm} \text{ for } i=1, 2, \ldots, r-1\eqcomma  \\
            X_{C, r}^{(n)}(t) &{} = X_{C, r}^{(n)}(0) + \int_{[0,t]\times [0,\infty)}\indicator{[0, \lambda_r^{(n)}(X^{(n)} (s-))]}{v} \xi_r(\differential{v}\times \differential{s}) 
            - \int_{[0,t]\times [0,\infty)}\indicator{[0, \lambda_{-r}^{(n)}(X^{(n)} (s-))]}{v} \xi_{-r}(\differential{v}\times \differential{s}) \\
            &{} \quad -  \int_{[0,t]\times [0,\infty)}\indicator{[0, \lambda_{P}^{(n)}(X^{(n)} (s-))]}{v} \xi_{P}(\differential{v}\times \differential{s})\eqcomma  \\
            X_S^{(n)}(t) &{} = X_S^{(n)}(0) - \int_{[0,t]\times [0,\infty)}\indicator{[0, \lambda_1^{(n)}(X^{(n)} (s-))]}{v} \xi_1(\differential{v}\times \differential{s}) 
            + \int_{[0,t]\times [0,\infty)}\indicator{[0, \lambda_{-1}^{(n)}(X^{(n)} (s-))]}{v} \xi_{-1}(\differential{v}\times \differential{s})\eqcomma \\
            X_P^{(n)}(t) &{} =  X_P^{(n)}(0) + \int_{[0,t]\times [0,\infty)}\indicator{[0, \lambda_{P}^{(n)}(X^{(n)} (s-))]}{v} \xi_{P}(\differential{v}\times \differential{s})\eqcomma 
        \end{aligned}
        \label{eq:X_n_traj}
    \end{align}
    where $\lambda_{k}^{(n)}: \setOfNonnegativeIntegers^{r+2} \mapsto \setOfPositiveReals \eqcomma$  for $k \in \{\pm 1, \pm 2, \ldots, \pm r, P\} $ denotes the propensity function of the $k$-th reaction,  $\xi_i, \xi_{-i}$, for $i=1, 2, \ldots, r$, and $\xi_P$ are independent \acp{PRM}  on $\setOfPositiveReals\times \setOfPositiveReals$ with intensity measure $\leb \ot \leb$. The \acp{PRM} $\xi_i, \xi_{-i}$, for $i=1, 2, \ldots, r$, and $\xi_P$ are defined on the same probability space $(\Omega, \history{}, \prob)$, and are independent of $X^{(n)}(0)$.  Assume $\history{}$ is $\prob$-complete and  associate to $(\Omega, \history{}, \prob)$ the filtration $(\history{t})_{t\ge 0}$ is generated by $X^{(n)}(0)$ and the \acp{PRM}, $\xi_{i}, \xi_{-i}, i=1,  \ldots, r $ and $\xi_P$,
  with $\history{0}$ containing all $\prob$-null sets in $\history{}$. The filtration $(\history{t})_{t\ge 0}$ is right continuous in the sense that 
    \begin{align*}
        \history{t+} \defeq \cap_{s>0}\history{t+s} = \history{t}\eqstop 
    \end{align*}
    Therefore, the stochastic basis $(\Omega, \history{}, (\history{t})_{t\ge 0}, \prob) $, also known as the filtered probability space, is complete or the \emph{usual  conditions} (also referred to as the Dellacherie conditions; see 
    \cite[Definition 1.3]{jacod2003limit})  are satisfied.

    Notice that the following {\em conservation law} holds at all times $t\geq 0$
    \begin{align} \label{eq:e_conserv_law}
    \begin{aligned}
        X_E^{(n)}(t) +\sum_{i=1}^{r} X_{C, i}^{(n)}(t) &\ = X_E^{(n)}(0) +\sum_{i=1}^{r} X_{C, i}^{(n)}(0) \equiv \conscon \eqcomma \\  
         X_S^{(n)}(t) +  X_P^{(n)}(t) +\sum_{i=1}^{r} X_{C, i}^{(n)}(t) &\ = X_S^{(n)}(0) +  X_P^{(n)}(0) +\sum_{i=1}^{r} X_{C, i}^{(n)}(0)\eqcomma 
     \end{aligned}     
    \end{align}
which means that we do not need to keep track of $X_E^{(n)}$. As such, we define the state vector $X^{(n)}$ as
    $$X^{(n)} \defeq ( X_{C, 1}^{(n)}, X_{C, 2}^{(n)}, \ldots, X_{C, r}^{(n)}, X_S^{(n)}, X_P^{(n)}).$$
    Then, $X^{(n)} $ is  a jump Markov process (see \cite[Chapter~4]{Ethier:1986:MPC}) with 
    paths in $D([0, \infty), \setOfNonnegativeIntegers^{r+2})$ and  generator 
    \begin{align*}
        \gen^{(n)} f(x) \defeq &{} \lambda_1^{(n)}(x)  \left(f(x+e^{(r+2)}_1-e^{(r+2)}_{r+1}) - f(x) \right) 
        + \lambda_{-1}^{(n)}(x) \left( f(x-e^{(r+2)}_1+e^{(r+2)}_{r+1}) -f(x)  \right) \\
        &{} +\sum_{i=2}^{r} \lambda_{i}^{(n)}(x) \left( f(x-e^{(r+2)}_{i-1} +e^{(r+2)}_{i}) - f(x) \right) 
        +\sum_{i=2}^{r} \lambda_{-i}^{(n)}(x) \left( f(x+e^{(r+2)}_{i-1} -e^{(r+2)}_{i}) - f(x) \right)\\
        &{} 
        + \lambda_{P}^{(n)}(x)  \left(f(x-e^{(r+2)}_{r}+e^{(r+2)}_{r+2}) - f(x) \right)\eqcomma 
    \end{align*}
    for bounded functions $f: \setOfNonnegativeIntegers^{r+2} \mapsto \setOfReals$, $x\in \setOfNonnegativeIntegers^{r+2}$. It will be convenient to write a typical state $x \in \setOfNonnegativeIntegers^{r+2}$ of the process $X^{(n)}$ as
    $$x =(\underbrace{x_{C,1},x_{C,2}, \ldots, x_{C,r}}_{x_C}, x_S, x_P),$$
    rather than the more conventional notation, $x= (x_1,x_2,\hdots,x_{r+2})$.
    
    Under the stochastic law of mass action (see \cite{Anderson2015book} and also \cite{Wilkinson2018SMS}), we  take the propensity functions to be of the following forms: for $x\in \setOfNonnegativeIntegers^{r+2}$, 
    \begin{align}
        \begin{aligned}
            \lambda_1^{(n)}(x) &{}\defeq \kappa_1^{(n)} x_{S} \Big(\conscon - \sum_{i=1}^{r}x_{C,i}\Big)\eqcomma \\
            \lambda_{-1}^{(n)}(x) &{} \defeq  \kappa_{-1}^{(n)}x_{C,1} \eqcomma \ \ \lambda_{P}^{(n)}(x)  \defeq \kappa^{(n)}_P x_{C,r}\eqcomma   \\
            \lambda_{i}^{(n)}(x) &{}\defeq  \kappa_{i}^{(n)}x_{C,i-1}, \quad \lambda_{-i}^{(n)}(x) \defeq  \kappa_{-i}^{(n)}x_{C,i} \eqcomma \text{ for } i =2, 3, \ldots, r\eqstop  \\
          \end{aligned}
        \label{eq:propensity_functions}
    \end{align}
   

    \section{The \acl{sQSSA}}
    \label{sec:sQSSA}

    In order to study various averaging phenomena, we will consider the scaled process $$Z^{(n)} \defeq (Z_{C, 1}^{(n)}, Z_{C, 2}^{(n)}, \ldots, Z_{C, r}^{(n)}, Z_S^{(n)}, Z_P^{(n)}),$$ where 
    \begin{align}
    \begin{aligned}
         Z_{C, i}^{(n)}(t) &{} = n^{-\alpha_{C, i}} X_{C, i}^{(n)} (n^\gamma t)\eqcomma \text{ for } i =1, 2, \ldots, r\eqcomma \quad 
         Z_S^{(n)}(t) 
         = n^{-\alpha_S} X_S^{(n)}(n^\gamma t)\eqcomma \quad
         Z_P^{(n)}(t) = n^{-\alpha_P} X_P^{(n)}(n^\gamma t) \eqcomma 
    \end{aligned}
        \label{eq:scaled_process}
    \end{align}
    where the exponents $\alpha_S, \alpha_P, \alpha_E, \alpha_{C, 1}, \alpha_{C, 2}, \ldots, \alpha_{C, r}$ capture the species abundance, and  the exponent $\gamma $ is used to speed up or slow down time. Also, define 
    \begin{align}
        Z_{E}^{(n)}(t) \defeq n^{-\alpha_E} X_E^{(n)}(n^\gamma t) =  n^{-\alpha_E} \conscon - \sum_{i=1}^{r} n^{\alpha_{C,i} - \alpha_E} Z_{C, i}^{(n)}(t) \eqstop 
        \label{eq:scaled_E}
    \end{align}
    In addition to the above parameters, we will also consider scaling exponents $\beta_i, \beta_{-i}$ for $i = 1, 2, \ldots, r$, and $\beta_P$ to describe the reaction-speeds so that 
    \begin{align} \label{eq:scaled_rate}
        \kappa_i^{(n)} = n^{\beta_i} \kappa_i\eqcomma \quad \kappa_{-i}^{(n)} = n^{\beta_{-i}} \kappa_{-i}\eqcomma \text{ for } i =1, 2, \ldots, r\eqcomma \quad \kappa_P^{(n)} = n^{\beta_P} \kappa_P,
    \end{align}
    for constants $\kappa_i, \kappa_{-i}$, for $i=1, 2, \ldots, r$, and $\kappa_{P}$ not depending on $n$. Note that these scaling parameters can take the value zero, and both positive,  or negative values. Such parametrizations, standard in the stochastic multiscale literature \citep{Kang:2013:STM,Kang:2014:CLT,Ball:2006:AAM,Kang2019QSSA,Ganguly2025tQSSA}, are useful as a means to describe the differences in the species abundances, and the reaction rates.  In the next section, we will choose a particular scaling regime and derive the \ac{sQSSA} as a consequence of the \ac{FLLN} for the sequence of scaled stochastic processes $Z^{(n)}$ as $n \to \infty$.

    In line with the literature on the \ac{sQSSA} for the standard single-stage \ac{MM} model \cite{Ball:2006:AAM,Kang2019QSSA}, we adopt the following scaling regime for the stochastic \ac{sQSSA} for the multistage \ac{MM} model:
    \begin{align}
    \begin{aligned}
        & \alpha_S=\alpha_P=1,\quad \alpha_E=\alpha_{C, i}=0, \text{ for } i =1, 2, \ldots, r\eqcomma \\
       & \beta_1=0,\quad \beta_{-1}=\beta_2=\cdots=\beta_{-r} = \beta_{r} = \beta_{P}=  1,\\
       &\gamma = 0\eqstop 
    \end{aligned}
        \label{eq:sQSSA_scalings}
    \end{align}
    The interpretation of the above choice is as follows: 
    \begin{itemize}
        \item The species $S$ and $P$ are abundant, \ie, $\BigO{n}$, and the species $E$ and $C_i$ are $\BigO{1}$.
        \item The reaction $S+E \longrightarrow C_1$ is slow, and all other reactions are fast. 
        \item We do not speed up or slow down time.
    \end{itemize}

We will assume the following condition throughout the paper. Write 
    \begin{align}
        Z_{V}^{(n)} \defeq (Z_S^{(n)}, Z_P^{(n)})\eqcomma \nonumber \quad 
        Z_C^{(n)}  \defeq (Z_{C, 1}^{(n)}, Z_{C, 2}^{(n)}, \ldots, Z_{C, r}^{(n)})\eqstop  \nonumber
    \end{align}

\begin{myAssumption}\label{assum:initial}
\begin{enumerate}[(i)]
    \item There exists an $\setOfPositiveReals^2$-valued random variable $Z_V(0)$ such that $Z^{(n)}_V(0) \rt Z_V(0)$ in $L^1(\Omega, \PP)$  as $n \rt \infty$.
    \item The conservation constant $\conscon$, introduced in \eqref{eq:e_conserv_law}, is deterministic and independent of $n$. 
\end{enumerate}

\end{myAssumption}

 \begin{myRemark}\label{rem:cons-const}        
 For convenience, we assume that the conservation constant $ \conscon$ in \eqref{eq:e_conserv_law} is deterministic and independent of $n$, although the proofs extend with minor changes for $n$-dependent random variables $\conscon^{(n)}$ under the mild assumption that $\conscon^{(n)} \stackrel{\PP} \Rt \conscon$ as $n\rt \infty$. For example, see \cite{Ganguly2025tQSSA} where the \ac{tQSSA} was derived for a simpler model of \ac{MM} enzyme kinetic reactions. 
\end{myRemark}

    As before, we write a typical state $z \in \setOfPositiveReals^{r+2}$ of the process $Z^{(n)}$ as
    $$z =(\underbrace{z_{C,1},z_{C,2}, \ldots, z_{C,r}}_{z_C}, z_S, z_P)$$
   instead of the conventional notation, $z= (z_1,z_2,\hdots,z_{r+2})$.
   We introduce the `scaled' propensity functions $\lambda_k$, which are independent of the scaling parameter $n$, as  
    \begin{align}
        \begin{aligned}
             \lambda_1(z) \defeq&\  \kappa_1 z_{S} \Big(\conscon - \sum_{i=1}^{r}z_{C,i}\Big) \eqcomma \ \ \lambda_{-1}(z) \defeq\ \kappa_{-1} z_{C,1}, \ \ \lambda_{P}(z) \defeq \  \kappa_P z_{C,r}\eqcomma  \\
                     \lambda_{i}(z) \defeq& \ \kappa_{i} z_{C,i-1}, \quad \lambda_{-i}(z) \defeq \  \kappa_{-i} z_{C,i}\eqcomma \text{ for } i =2, 3, \ldots, r\eqstop
                    \end{aligned}
        \label{eq:scaled_propensities}
    \end{align}
    Note that the $\l_k(z)$ do not depend on the coordinate $z_P$, that is, for $z = (z_C, z_S, z_P)$
   \begin{align}\label{eq:prop-nodep-zp}
        \lambda_k(z) \equiv \lambda_k(z_C, z_S, z_P) \equiv \lambda_k(z_C,z_S), \quad k \in \{\pm 1, \pm 2, \hdots \pm r, P\}\eqstop 
   \end{align}
   Because of the scaling described in \eqref{eq:scaled_process}, \eqref{eq:scaled_rate} and \eqref{eq:sQSSA_scalings}, the relation between the original propensity functions $\lambda_{k}^{(n)}$ and their scaled versions $\lambda_{k}$ is as follows:  for $x =(x_{C}, x_S, x_P) \in \setOfNonnegativeIntegers^{r+2}$ setting $z = (x_{C}, x_S/n, x_P/n)$ gives 
   $$\lambda^{(n)}_k(x) = n\lambda_k(z), \quad k \in \{\pm 1, \pm 2, \hdots r, P\}\eqstop $$
 
   It now readily follows from \eqref{eq:X_n_traj} that the trajectory equations for the scaled process $Z^{(n)}$ are given by  
   
    \begin{align}
        \begin{aligned}
            Z_{C, i}^{(n)}(t) &{} = Z_{C, i}^{(n)}(0) + \int_{[0,t]\times [0,\infty)}\indicator{[0, n \lambda_i(Z^{(n)} (s-))]}{v} \xi_i(\differential{v}\times \differential{s}) 
            + \int_{[0,t]\times [0,\infty)}\indicator{[0, n \lambda_{-(i+1)}(Z^{(n)} (s-))]}{v} \xi_{-(i+1)}(\differential{v}\times \differential{s}) \\
            &{} \quad - \int_{[0,t]\times [0,\infty)}\indicator{[0, n \lambda_{-i}(Z^{(n)} (s-))]}{v} \xi_{-i}(\differential{v}\times \differential{s}) 
            -  \int_{[0,t]\times [0,\infty)}\indicator{[0, n \lambda_{i+1}(Z^{(n)} (s-))]}{v} \xi_{i+1}(\differential{v}\times \differential{s})\eqcomma\\
            &\quad 
            \text{ for } i=1, 2, \ldots, r-1\eqcomma  \\
            Z_{C, r}^{(n)}(t) &{} = Z_{C, r}^{(n)}(0) + \int_{[0,t]\times [0,\infty)}\indicator{[0, n \lambda_r(Z^{(n)} (s-))]}{v} \xi_r(\differential{v}\times \differential{s}) 
            - \int_{[0,t]\times [0,\infty)}\indicator{[0, n \lambda_{-r}(Z^{(n)} (s-))]}{v} \xi_{-r}(\differential{v}\times \differential{s}) \\
            &{} \quad - \int_{[0,t]\times [0,\infty)}\indicator{[0, n \lambda_{P}(Z^{(n)} (s-))]}{v} \xi_{P}(\differential{v}\times \differential{s})\eqcomma  \\
            Z_S^{(n)}(t) &{} = Z_S^{(n)}(0) - \frac{1}{n}\int_{[0,t]\times [0,\infty)}\indicator{[0, n \lambda_1(Z^{(n)} (s-))]}{v} \xi_1(\differential{v}\times \differential{s}) 
            + \frac{1}{n} \int_{[0,t]\times [0,\infty)}\indicator{[0, n \lambda_{-1}(Z^{(n)} (s-))]}{v} \xi_{-1}(\differential{v}\times \differential{s})\eqcomma \\
            Z_P^{(n)}(t) &{} =  Z_P^{(n)}(0) + \frac{1}{n}\int_{[0,t]\times [0,\infty)}\indicator{[0, n \lambda_{P}(Z^{(n)} (s-))]}{v} \xi_{P}(\differential{v}\times \differential{s})\eqstop 
        \end{aligned}
        \label{eq:Z_n_traj}
    \end{align}
    Similarly, from \eqref{eq:scaled_E} we can write the trajectory equation for $Z_E^{(n)}$ as
    \begin{align}
        \begin{aligned}
            Z_E^{(n)}(t) &{} = Z_E^{(n)}(0) - \int_{[0,t]\times [0,\infty)}\indicator{[0, n \lambda_1(Z^{(n)} (s-))]}{v} \xi_1(\differential{v}\times \differential{s}) 
            +\int_{[0,t]\times [0,\infty)}\indicator{[0, n \lambda_{-1}(Z^{(n)} (s-))]}{v} \xi_{-1}(\differential{v}\times \differential{s}) \\
            &{} \quad +  \int_{[0,t]\times [0,\infty)}\indicator{[0, n \lambda_{P}(Z^{(n)} (s-))]}{v} \xi_{P}(\differential{v}\times \differential{s})\eqstop
        \end{aligned}
        \label{eq:Z_En_traj}
    \end{align}
    In the next section, we show that the deterministic \ac{sQSSA} can be derived as an \ac{FLLN} under the above scaling regime. 

    \subsection{\acl{FLLN}}
    \label{sec:FLLN}
    From \eqref{eq:Z_n_traj}, and \eqref{eq:Z_En_traj}, it is clear that the variables $Z_{C, i}^{(n)}$ for $i=1, 2, \ldots, r$ and $Z_E^{(n)}$ are fast variables whereas the variables $Z_S^{(n)}$ and $Z_P^{(n)}$ are slow variables. In order to derive the \ac{sQSSA}, we will consider the occupation measure of the fast variables. To this end, notice that because of the conservation law in \eqref{eq:e_conserv_law} the process $Z_C^{(n)}$ takes values in the finite state space,
\begin{align}\label{eq:bdd-Z-M}
\BddZ_{\conscon,+}^r \defeq \lf\{u=(u_1,u_2,\hdots, u_r) \in \setOfNonnegativeIntegers^r: \sum_{i=1}^r u_i \leq \conscon\ri\}.
\end{align}

  Note that $Z_E^{(n)} = \conscon - \norm{Z_C^{(n)}}_1$, where recall that $\norm{\cdot}_1$ denotes the $\ell_1$-norm on $\setOfReals^{r}$, \ie,  $\norm{u}_1 \defeq \sum_{i=1}^{r} \absolute{u_i}$ for $u=(u_1, u_2, \ldots, u_r) \in \setOfReals^r$.    
 Define the  occupation measure $\occ_n$ of $Z_C^{(n)}$ by
    \begin{align}
        \occ_n(z_C \times [0, t])\defeq  \int_0^t \indic{\{Z_C^{(n)}(s) =z_C\}} \differential{s}, \quad z_C \in\BddZ_{\conscon,+}^r \eqstop 
        \label{eq:occupation_measure-0}
    \end{align}
  For $ A \subset \BddZ_{\conscon,+}^r$, we write
     \begin{align}
        \occ_n(A \times [0, t])\defeq \int_0^t \indicator{A}{Z_C^{(n)}(s)} \differential{s} =\sum_{z_C \in A} \int_0^t \indic{\{Z_C^{(n)}(s) =z_C\}} \differential{s}\eqstop  
        \label{eq:occupation_measure}
    \end{align}
   For a fixed $T$ such that $0< T< \infty$,  the occupation measures $\Gamma_n$ are random measures on $ G_T \defeq \BddZ_{\conscon,+}^r\times [0,T]$. 
In other words, $\occ_n$ are $\mathcal{M}(G_T)$-valued random variables, where  $\mathcal{M}(G_T)$ denotes the space of finite (non-negative) measures on $G_T$ endowed with the topology of weak convergence  \cite{Billingsley1999Convergence,Kallenberg2017RandomMeasures,jacod2003limit}.

The following lemma is a simple consequence of  \Cref{assum:initial}.

\begin{myLemma}\label{lem:Zv_mmt-bd}
For any $T>0$,  $$\sup_n\EE\lf[\sup_{0\leq t\leq T}\norm{Z_{V}^{(n)}(t)}_1\ri] \leq  C \equiv \sup_n \EE\lf[\norm{Z_{V}^{(n)}(0)}_1\ri] +J.$$
\end{myLemma}

\begin{proof}[Proof of \Cref{lem:Zv_mmt-bd}]
The proof follows readily from \eqref{eq:e_conserv_law}, which --- when expressed in terms of the scaled process $Z^{(n)}$ --- shows that for any $t>0$, 
\begin{align*}
\norm{Z_{V}^{(n)}(t)}_1+n^{-1}\norm{Z_{C}^{(n)}(t)}_1  = &\
\norm{Z_{V}^{(n)}(0)}_1+n^{-1}\norm{Z_{C}^{(n)}(0)}_1 
\leq 
\norm{Z_{V}^{(n)}(0)}_1+n^{-1}\conscon.
\end{align*}
\end{proof}


The following proposition establishes the tightness of the sequence of random variables $\{(\occ_n, Z_V^{(n)}): n\in \setOfNaturals\}$.

    \begin{myProposition}
        The collection $\{(\occ_n, Z_V^{(n)}): n\in \setOfNaturals\}$ is tight as a sequence of $\mathcal{M}(G_T)\times D([0,T], {\R_+^2})$-valued random variables. Furthermore, if $Z_V$ is a limit point of $Z^{(n)}_V$, then the paths of $Z_V$ are almost surely in $C([0,T], \R^2)$.
        \label{proposition:occupation_tightness}
    \end{myProposition}
    \begin{proof}[Proof of \Cref{proposition:occupation_tightness}]
       
      By \cite[Theorem~2.11]{Budhiraja2019Analysis} the tightness of $\{\occ_n\}$ is equivalent to the tightness of the  sequence of corresponding  deterministic mean measures $\{\nu_n \}$ in $\mathcal{M}(G_T)$, where the $\nu_n$  are defined by
    \begin{align}
        \nu_n(A\times [0, t]) \defeq \Eof{\Gamma_n(A\times [0, t])} = \int_0^t \probOf{Z_C^{(n)}(s) \in A} \differential{s}\eqcomma \quad A \subset \BddZ_{\conscon,+}^r\eqstop 
    \end{align} 
   But this is obvious as $G_T$ is compact.  

   We now establish $C$-tightness of $Z_V^{(n)}$ in $D([0,T], {\R_+^2})$. To this end, notice that 
   \begin{align}
           \begin{aligned}
               \label{eq:Z_V_n_traj}
           Z_V^{(n)}(t) = &\ Z_V^{(n)}(0) + \int_{[0,t]\times [0,\infty)} \begin{pmatrix} -n^{-1} \indicator{[0, n \lambda_1(Z^{(n)} (s-))]}{v} \\
               0 
           \end{pmatrix} \xi_1(\differential{v}\times \differential{s}) \\  
           &{} + \int_{[0,t]\times [0,\infty)} \begin{pmatrix} n^{-1} \indicator{[0, n \lambda_{-1}(Z^{(n)} (s-))]}{v} \\
               0
           \end{pmatrix} \xi_{-1}(\differential{v}\times \differential{s})\\
           &{} + \int_{[0,t]\times [0,\infty)} \begin{pmatrix} 0 \\
               n^{-1} \indicator{[0, n \lambda_{P}(Z^{(n)} (s-))]}{v}
           \end{pmatrix} \xi_{P}(\differential{v}\times \differential{s})\eqstop \\
         \end{aligned}
       \end{align}
   Rewrite the trajectory equation for $Z_V^{(n)}$ from \eqref{eq:Z_V_n_traj} as
  \begin{align}
            \label{eq:Z-split-eq}
        Z_V^{(n)}(t) \equiv \Upsilon_V^{(n)}(t)+n^{-1}M_V^{(n)}(t),
    \end{align}
 where the process $\Upsilon_V^{(n)}$ is given by
\begin{align} \label{eq:Z-split-det}
\begin{aligned}
   \Upsilon_V^{(n)}(t) =&\ Z_V^{(n)}(0) - e_1^{(2)} \int_0^t  \lambda_1(Z^{(n)}(s)) \differential{s} + e_1^{(2)} \int_0^t \lambda_{-1}(Z^{(n)}(s)) \differential{s} 
        + e_2^{(2)} \int_0^t  \lambda_{P}(Z^{(n)}(s)) \differential{s}\eqcomma 
\end{aligned}      
    \end{align}
    and $M_V^{(n)}$ is a zero-mean martingale given by
    \begin{align}\label{eq:M_n_defn}
        \begin{aligned}
        M_V^{(n)}(t) \defeq& -e_1^{(2)} \int_{[0,t]\times [0,\infty)} \indicator{[0, n \lambda_1(Z^{(n)} (s-))]}{v} \tilde{\xi}_1(\differential{v}\times \differential{s}) 
        + e_1^{(2)} \int_{[0,t]\times [0,\infty)} \indicator{[0, n \lambda_{-1}(Z^{(n)} (s-))]}{v} \tilde{\xi}_{-1}(\differential{v}\times \differential{s}) \\
        &{} + e_2^{(2)} \int_{[0,t]\times [0,\infty)} \indicator{[0, n \lambda_{P}(Z^{(n)} (s-))]}{v} \tilde{\xi}_{P}(\differential{v}\times \differential{s}).
        \end{aligned} 
    \end{align}
     Here $e_1^{(2)}\defeq {(1, 0)}\trans, e_2^{(2)} \defeq (0, 1)\trans$ are the two unit canonical basis vectors in $\setOfReals^2$, $\tilde{\xi}_1, \tilde{\xi}_{-1}$, and $\tilde{\xi}_{P}$ are the compensated (centered) \acp{PRM} corresponding to the \acp{PRM} $\xi_1, \xi_{-1}$, and $\xi_{P}$, respectively. Observe that by the \ac{BDG} inequality 
     \begin{align*}
  \Eof{\sup_{t\leq T} \norm{M_V^{(n)}(t)}^2} &\leq \ \const_0 \Eof{[M_V^{(n)}]_T} = \const_0\EE\Bigg(\int_{[0,T]\times [0,\infty)}  \indicator{[0, n \lambda_1(Z^{(n)} (s-))]}{v} \xi_1(\differential{v}\times \differential{s}) \\
  &\ + \int_{[0,T]\times [0,\infty)} \indicator{[0, n \lambda_{-1}(Z^{(n)} (s-))]}{v} \xi_{-1}(\differential{v}\times \differential{s}) \\ 
  &\ + \int_{[0,T]\times [0,\infty)} \indicator{[0, n \lambda_{P}(Z^{(n)} (s-))]}{v} \xi_{P}(\differential{v}\times \differential{s})\Bigg)\\
  =&\ \const_0 n \Eof{ \int_0^T \left(\lambda_1(Z^{(n)}(s)) + \lambda_{-1}(Z^{(n)}(s)) + \lambda_P(Z^{(n)}(s))  \right)\differential{s} }\\
  \leq&\ \const_0 n \conscon\lf((\kappa_{-1}+\kappa_P)T + \kappa_1\Eof{\sup_{0\leq t\leq T}\norm{Z_{V}^{(n)}(t)}_1} T\ri)
  \equiv \const_1(T)n \eqcomma 
     \end{align*}
where for the last inequality we used \Cref{lem:Zv_mmt-bd}. It follows that
\begin{align}
    \label{eq:mart-L2-conv}
 n^{-1}\sup_{t\leq T}\norm{M_V^{(n)}(t)} \nrt 0 \text { in } L^2(\PP).   
\end{align}
Thus, for establishing tightness of $Z_{V}^{(n)}$ in $D([0,T], {\R_+^2})$, it is enough to prove tightness of $\Upsilon_V^{(n)}$.
Notice that
\begin{align*}
\|\Upsilon_V^{(n)}(t) - \Upsilon_V^{(n)}(t')\| \leq \int_{t\wedge t'}^{t\vee t'} \left(\lambda_1(Z^{(n)}(s)) + \lambda_{-1}(Z^{(n)}(s)) + \lambda_P(Z^{(n)}(s))  \right)\differential{s}.
\end{align*}
Hence, taking supremum and then expectation, because of   \Cref{lem:Zv_mmt-bd},
\begin{align*}
\Eof{\sup\limits_{|t-t'|\leq h}\|\Upsilon_V^{(n)}(t) - \Upsilon_V^{(n)}(t')\|} &{} \leq \conscon(\kappa_{-1}+\kappa_P)h + \conscon\kappa_1\Eof{\sup_{0\leq t\leq T}\norm{Z_{V}^{(n)}(t)}_1} h 
\leq \const_3(T)h\eqcomma 
\end{align*}
which, in fact, shows that $\{\Upsilon_V^{(n)}\}$ is tight in $C([0,T], {\R^2})$ (see the discussion following Theorem 3.2 in \cite{Whitt2007MCLT}). This readily implies the second assertion.
\end{proof}

For every fixed $z_S \in \R_+$, and for a measurable $f:\BddZ^r_{\conscon,+} \rt \R$,  define the operator $\fgen_{z_S}$ as follows: 
   \begin{align}
            \begin{aligned}
            \fgen_{z_S}f(z_C)\defeq &\  \kappa_1 z_S (\conscon - \norm{z_C}_1)\left(f(z_C + e_1^{(r)}) - f(z_C)\right) 
            + \sum_{i=2}^{r}\kappa_{i}z_{C,i-1} \left(f(z_C - e_{i-1}^{(r)} + e_i^{(r)}) - f(z_C)\right) \\
            &{} +  \kappa_P z_{C,r} \left(f(z_C - e_r^{(r)}) - f(z_C)\right) 
            + \sum_{i=1}^{r}  \kappa_{-i}z_{C,i} \left(f(z_C + e_{i-1}^{(r)} - e_{i}^{(r)}) - f(z_C)\right) \eqcomma 
           \end{aligned}
            \label{eq:fast_generator}
        \end{align}    
   where $z_{C} = (z_{C,1},z_{C,2},\hdots, z_{C,r}) \in \BddZ^r_{\conscon,+}$, $e^{(r)}_0 \equiv \bm{0}_r$, $e^{(r)}_i, \ i=1,2,\hdots,r$ are the canonical basis vectors of $\R^r$. 
Notice that  for every $z_S \in \R_{+}$, the operator $\fgen_{z_S}$ is the generator of a  $\BddZ_{\conscon,+}^r$-valued \ac{CTMC} $\tilde Z_C \equiv \tilde Z^{z_S}_{C}$ whose paths can be described as follows: for $i=1,2,\hdots, r,$ 
   \begin{align}\label{eq:fast-proc-frozen-slow}
        \begin{aligned}
           \tilde Z_{C, i}(t) &{} = Z_{C, i}(0) + \indic{\{i= 1\}}\int_{[0,t]\times [0,\infty)}\indicator{\lf[0,  \kappa_1 z_S \lf(\conscon - \norm{\tilde Z_C(s-)}_1\ri)\ri]}{v} \xi_1(\differential{v}\times \differential{s}) \\
           &{}\quad + \indic{\{i\neq 1\}}\int_{[0,t]\times [0,\infty)}\indicator{\lf[0,  \kappa_i \tilde Z_{C,i-1}(s-)\ri]}{v} \xi_i(\differential{v}\times \differential{s}) 
            + \indic{\{i\neq r\}}\int_{[0,t]\times [0,\infty)}\indicator{[0,  \kappa_{-(i+1)}\tilde Z_{C,i+1} (s-)]}{v} \xi_{-(i+1)}(\differential{v}\times \differential{s}) \\
            &{} \quad - \int_{[0,t]\times [0,\infty)}\indicator{[0,\kappa_{-i}\tilde Z_{C,i}(s-)]}{v} \xi_{-i}(\differential{v}\times \differential{s}) 
            -  \indic{\{i\neq r\}}\int_{[0,t]\times [0,\infty)}\indicator{[0,  \kappa_{i+1}\tilde Z_{C,i} (s-)]}{v} \xi_{i+1}(\differential{v}\times \differential{s})\\
            &{} \quad - \indic{\{i=r\}}\int_{[0,t]\times [0,\infty)}\indicator{[0, \kappa_{P}\tilde Z_{C,r} (s-)]}{v} \xi_{P}(\differential{v}\times \differential{s}), 
            \end{aligned}
    \end{align} 
where we reused the \acp{PRM} from \eqref{eq:Z_n_traj}.
  The stochastic process $\tilde Z_C$ can be thought of as a process capturing the dynamics of the fast process $Z^{(n)}_C$ with the slow component $Z^{(n)}_S$ frozen at $z_S$.  

\begin{myLemma}\label{lem:fast-stat-dist}
For every fixed $z_S \in \R_+$, the \ac{CTMC} corresponding to the generator $\fgen_{z_S}$ has the unique stationary distribution 
$$\pi_{z_S}\equiv \sys{Multinomial}(\conscon, p_1(z_S), p_2(z_S), \hdots, p_r(z_S))$$
on $\BddZ_{\conscon,+}^r$, that is, 
\begin{align*}
            \pi_{z_S}(z_C) = \frac{\conscon!}{ (\conscon-\sum_{i=1}^{r}z_{C,i})! \prod_{i=1}^{r} z_{C,i}!} \lf(1- \sum_{i=1}^{r}p_i(z_S)\ri)^{(\conscon-\sum_{i=1}^{r}z_{C,i})}\prod_{i=1}^{r} p_i(z_S)^{z_{C, i}} \eqcomma 
        \end{align*}
for $z_C = (z_{C,1}, z_{C, 2}, \ldots, z_{C, r}) \in \BddZ_{\conscon,+}^r$, with 
        \begin{align*}
            \begin{aligned}
                p_1(z_S) &= \left(1+ a_1(z_S) + \sum_{i=2}^{r}  \frac{1}{\prod_{j=2}^{i} a_j} \right)^{-1} \eqcomma\;\quad  
                p_{i}(z_S) 
                = \frac{p_1(z_S)}{\prod_{j=2}^{i} a_j}, \quad \text{ for } i=2, 3, \ldots, r\eqcomma 
            \end{aligned}
        \end{align*}
        where the numbers $a_1(z_S), a_2, a_3, \ldots, a_{r}$  satisfy the following recursive relations 
        \begin{align}\label{eq:ai-def}
            \begin{aligned}
                a_{r} &= \frac{ (\kappa_{-r}+\kappa_{P}) }{ \kappa_r}\eqcomma \quad   a_1(z_S)  = \frac{\kappa_{-1}}{\kappa_1 z_S} + \frac{\kappa_{P}}{\kappa_1 z_S} \prod_{i=2}^{r} \frac{1}{a_i}\eqcomma \quad 
                a_{i} 
                = \frac{(\kappa_{-i} + \kappa_{i+1})}{ \kappa_i } - \frac{\kappa_{-(i+1)}}{a_{i+1} \kappa_i}  \text{ for } i=2, 3,  \ldots, r-1\eqstop 
            \end{aligned}
        \end{align}
\end{myLemma}

\begin{proof}[Proof of \Cref{lem:fast-stat-dist}]
This follows from  \Cref{lemma:stationary_distribution} in Appendix \ref{sec:stationary_distribution} by putting $l_1 = \kappa_1 z_S,\ l_{-1} = \kappa_{-1},\ l_i = \kappa_i,\ l_{-i} = \kappa_{-i},\ l_{r+1} = \kappa_P$.
\end{proof}

It is easy to see that 
\begin{align}
    \label{eq:deriv-bd-p1}
    |\partial p_1(z_S)| \vee |\partial^2p_1(z_S)| \leq \const_{p_1}\eqcomma 
\end{align}
for some constant $\const_{p_1}$. 
Next, for $k=-1,1,P,$ define  the averaged propensity functions $\lambda_k^{\mrm{avg}}$ as
\begin{align}\label{eq:def-avg-prop}
     \lambda_k^{\mrm{avg}}(z_S) &\defeq  \sum_{z_C \in \BddZ_{\conscon,+}^r } \lambda_k(z_C, z_S) \pi_{z_S}(z_C) = \begin{cases}
         \kappa_1 \conscon z_S \lf(1- \sum_{i=1}^{r} p_i(z_S)\ri), &\ k=1\eqcomma \\
         \kappa_{-1}\conscon p_{1}(z_S), &\ k=-1\eqcomma \\
         \kappa_P \conscon p_r(z_S), &\ k=P.
     \end{cases}
\end{align}


We now prove our main result of this section. 

\begin{myTheorem}
        \label{thm:FLLN}
        Assume that $Z_V^{(n)}(0)$ is non-random and $Z_V^{(n)}(0) \to Z_V(0)$ as $n\to \infty$, where $Z_V(0)$ is a vector in $\setOfPositiveReals^2$. Then, the sequence $\{(\occ_n, Z_V^{(n)})\}$ converges in probability to $(\pi_{Z_S}\star \leb, Z_V = (Z_S, Z_P))$ in $\mathcal{M}(G_T)\times D([0,T], \R_+^2)$  where the measure $\pi_{Z_S}\star \leb$ is given by 
    \begin{align}
    \label{eq:def-lim-occ-meas}
    \pi_{Z_S}\star \leb (A \times [0,t]) \defeq  \int_0^t  \pi_{Z_S(s)}(A)\differential{s}, \quad A \subset \BddZ_{\conscon,+}^r,\ t\leq T\eqcomma 
    \end{align}
  and $Z_V = (Z_S, Z_P)$ is the solution of the system of \acp{ODE}
        \begin{align}\label{eq:ZV-limit-ODE}
        \begin{aligned}
            \timeDerivative{Z_S(t)} &= -\lambda_1^{\mrm{avg}}(Z_S(t)) + \lambda_{-1}^{\mrm{avg}}(Z_S(t)) 
            =-h(Z_S(t))\eqcomma \qquad 
            \timeDerivative{Z_P(t)} = \lambda_P^{\mrm{avg}}(Z_S(t)) = \kappa_{P} \conscon p_r(Z_S(t)) \eqcomma 
        \end{aligned}    
        \end{align}
        where the effective propensity function $h$ is given by 
        \begin{align*}
    h(y) \defeq& \kappa_1 \conscon \left(1 - \sum_{i=1}^{r}p_i(y) \right)y - \kappa_{-1} \conscon  p_1(y) 
    =
    \f{\conscon (\kappa_1 c_2-\kappa_{-1})y}{c_1y+c_2} 
    = \f{\conscon c_3y}{c_1y+c_2},
\end{align*}
with the constants
\begin{align*} 
    c_1=1+\sum_{i=2}^r \prod_{j=2}^{i}\frac{1}{a_j}, 
    \quad 
    c_2 = \frac{\kappa_{-1}}{\kappa_1} + \frac{\kappa_{P}}{\kappa_1} \prod_{i=2}^{r} \frac{1}{a_i}, 
    \quad 
    c_3=\kappa_{P}\prod_{i=2}^{r}\f{1}{a_i},
\end{align*}
and $a_i$, $\pi_{Z_S}$, $p_i, \ i=1,2,\hdots, r$ are as in \Cref{lem:fast-stat-dist}.
\end{myTheorem}     

\begin{proof}[Proof of \Cref{thm:FLLN}]
 By Proposition \ref{proposition:occupation_tightness} and Prokhorov's theorem, $\{(\occ_n, Z_V^{(n)})\}$ is relatively compact in $\mathcal{M}(G_T)\times D([0,T], {\R_+^2})$. Let $(\occ, Z_V)$ be a limit point of $\{(\occ_n, Z_V^{(n)})\}$, that is, there exists a subsequence along which $(\occ_n, Z_V^{(n)}) \RT (\occ, Z_V)$.
By Skorokhod theorem, we can assume that $(\occ_n, Z_V^{(n)}) \rt (\occ, Z_V)$ a.s.  along this subsequence.  Notice that by \Cref{proposition:occupation_tightness}, $Z_V$ is continuous, and from the proof of \Cref{proposition:occupation_tightness}, we have 
\begin{align}\label{eq:unif-conv-Z_V}
\sup_{t\leq T}\|Z_V^{(n)}(t) - Z_V(t)\| \stackrel{n\rt\infty} \Rt 0, \quad \text{a.s.}
\end{align}
Here, by a slight abuse of notation we continued to denote the subsequence by $(\occ_n, Z_V^{(n)})$. We will show the limit is unique and is independent of the subsequence. This ensures that the convergence holds along the entire sequence.

As before, we denote a typical element in the state space, $\BddZ_{\conscon,+}^r$, of the process $Z^{(n)}_C$ by $z_C =(z_{C,1},z_{C,2},\hdots, z_{C,r})$. For any measurable function $f: \BddZ_{\conscon,+}^r \rt \setOfReals$, by  It\^o's formula for \acp{SDE} driven by \acp{PRM} (see 
\cite[Theorem~5.1]{IkedaWatanabe2014Stochastic}), we have
\begin{align} \label{eq:ZC-ito}
\begin{aligned}
&f(Z_C^{(n)}(t)) - f(Z_C^{(n)}(0))\\
=&\ \int_{[0,t]\times [0,\infty)} \left(f(Z_C^{(n)}(s-) + e_1^{(r)}) - f(Z_C^{(n)}(s-) )\right) \indicator{\lf[0, n \lambda_1(Z_C^{(n)}(s-), Z_S^{(n)}(s-))\ri]}{v} \xi_1(\differential{v}\times \differential{s}) \\
        &{} + \sum_{i=2}^{r}\int_{[0,t]\times [0,\infty)} \left(f(Z_C^{(n)}(s-) - e_{i-1}^{(r)} + e_i^{(r)}) - f(Z_C^{(n)}(s-) )\right) \indicator{[0, n \lambda_{i}(Z^{(n)}_C (s-))]}{v} \xi_{-(i+1)}(\differential{v}\times \differential{s})\\
         &{} + \sum_{i=1}^{r} \int_{[0,t]\times [0,\infty)} \left(f(Z_C^{(n)}(s-) + e_{i-1}^{(r)} - e_{i}^{(r)}) - f(Z_C^{(n)}(s-) )\right) \indicator{[0, n \lambda_{-i}(Z^{(n)}_C(s-))] }{v} \xi_{i}(\differential{v}\times \differential{s})\\
        &{} + \int_{[0,t]\times [0,\infty)} \left(f(Z_C^{(n)}(s-) - e_r^{(r)}) - f(Z_C^{(n)}(s-) )\right) \indicator{[0, n \lambda_{P}(Z^{(n)}_C(s-))] }{v} \xi_{P}(\differential{v}\times \differential{s})\\
        =&\  n \int_0^t \fgen_{Z^{(n)}_S(s)}f(Z^{n)}_C(s)) \differential{s}+ M^{(n)}_{C, f}(t)\\
        =&\  n\sum_{z_C \in \BddZ_{\conscon,+}^r }\int_0^t \fgen_{Z^{(n)}_S(s)}f(z_C) \occ_n(z_C\times \differential{s})+ M^{(n)}_{C,f}(t).
 \end{aligned}       
\end{align}
Here, the stochastic process $M^{(n)}_{C,f}$ is a martingale given by
\begin{align*}
 M^{(n)}_{C,f}(t) 
  \defeq & \int_{[0,t]\times [0,\infty)} \left(f(Z_C^{(n)}(s-) + e_1^{(r)}) - f(Z_C^{(n)}(s-) )\right) \indicator{[0, n \lambda_1(Z^{(n)}_S (s-),Z^{(n)}_C (s-))]}{v} \tilde\xi_1(\differential{v}\times \differential{s}) \\
        &\quad + \sum_{i=2}^{r}\int_{[0,t]\times [0,\infty)} \left(f(Z_C^{(n)}(s-) -  e_{i-1}^{(r)} + e_i^{(r)}) - f(Z_C^{(n)}(s-) )\right) \indicator{[0, n \lambda_{i}(Z^{(n)}_C (s-))]}{v} \tilde\xi_{i}(\differential{v}\times \differential{s})\\
         &\quad + \sum_{i=1}^{r} \int_{[0,t]\times [0,\infty)} \left(f(Z_C^{(n)}(s-) + e_{i-1}^{(r)} - e_{i}^{(r)}) - f(Z_C^{(n)}(s-) )\right) \indicator{[0, n \lambda_{-i}(Z^{(n)}_C(s-))] }{v} \tilde\xi_{-i}(\differential{v}\times \differential{s})\\
        &\quad + \int_{[0,t]\times [0,\infty)} \left(f(Z_C^{(n)}(s-) - e_r^{(r)}) - f(Z_C^{(n)}(s-) )\right) \indicator{[0, n \lambda_{P}(Z^{(n)}_C(s-))] }{v} \tilde\xi_{P}(\differential{v}\times \differential{s}).
\end{align*}
Since $\BddZ_{\conscon,+}^r$ is a finite set, any function $f$ on $\BddZ_{\conscon,+}^r$ is continuous and  bounded, and it follows 
\begin{align*}
n^{-1}\lf(f(Z_C^{(n)}(t))-f(Z_C^{(n)}(0))\ri) \stackrel{n \rt \infty}\Rt 0, \quad n^{-2}\Eof{\sup_{t\leq T} \norm{M^{(n)}_{C,f}(t)}^2} \leq \const_f(T)n^{-1} \stackrel{n \rt \infty}\Rt 0.
\end{align*}
Furthermore, notice by the assumption that $(\occ_n, Z_V^{(n)}) \stackrel{n \rt \infty}\Rt (\occ, Z_V = (Z_S, Z_P))$ a.s., which means 
\begin{align*}
\sum_{z_C \in \BddZ_{\conscon,+}^r }\int_0^t \fgen_{Z^{(n)}_S(s)}f(z_C) \occ_n(z_C \times \differential{s}) \stackrel{n \rt \infty}\Rt  \sum_{z_C \in \BddZ_{\conscon,+}^r }\int_0^t \fgen_{Z_S(s)}f(z_C) \occ(z_C \times \differential{s}).
\end{align*}
Now rearranging the terms in \eqref{eq:ZC-ito} gives
\begin{align*}
\sum_{z_C \in \BddZ_{\conscon,+}^r }\int_0^t \fgen_{Z^{(n)}_S(s)}f(z_C) \occ_n(z_c\times  \differential{s}) = n^{-1}\lf(f(Z_C^{(n)}(t))-f(Z_C^{(n)}(0))\ri) -n^{-1}M^{(n)}_{C,f}(t).
\end{align*}
Taking $n\rt \infty$, it follows that there exists an $\Omega_0 \subset \Omega$ such that $\PP(\Omega_0)=1$ and for all $\omega \in \Omega_0$,
$$\sum_{z_C \in \BddZ_{\conscon,+}^r }\int_0^t  \fgen_{Z_S(s, \om)}f(z_C) \occ(z_C\times  \differential{s})(\om) =0.$$
Notice the probability one set $\Omega_0$ where the above equality holds can depend on $f$. But since the set of functions from $\BddZ_{\conscon,+}^r \rt \R$ is separable (it is isomorphic to $\R^{(J+1)^r}$), the $\Omega_0$ set can be taken independent of $f$. 
Since for any $t>0$, $\occ_n(\BddZ_{\conscon,+}^r \times [0,t]) = t$, we have $\occ(\BddZ_{\conscon,+}^r \times [0,t]) = t$. Thus, splitting $\occ(z_C \times \differential{s}) \equiv \occ_{(1|2)}(z_C|s)\differential{s}$ we see that for a.a. $s$ 
$$\sum_{z_C \in \BddZ_{\conscon,+}^r } \fgen_{Z_S(s)}f(z_C) \occ_{(1|2)}(z_C|s) =0.$$
In other words, for a.a. $s$, $\occ_{(1|2)}(\cdot|s)$ is a stationary distribution of the generator $\fgen_{Z_S(s)}$, and by Lemma \ref{lem:fast-stat-dist}, $\occ_{(1|2)}(\cdot|s) \equiv \pi_{Z_S(s)}(\cdot),$ which shows that  $\occ = \pi_{Z_S} \star \leb$ (defined in \eqref{eq:def-lim-occ-meas}).

Recalling that the  functions $\lambda_k$ depends on $Z^{(n)} =(Z^{(n)}_C, Z^{(n)}_S, Z^{(n)}_P)$ only through $(Z^{(n)}_C, Z^{(n)}_S)$ (see \eqref{eq:prop-nodep-zp}), we now have by \Cref{lem:conv-int} that for $k =P, -1, 1$,
\begin{align}\label{eq:weak-conv-occ-meas}
\begin{aligned}
   \int_0^t \lambda_k(Z^{(n)}(s)) \differential{s} =&\ \sum_{z_C \in \BddZ_{\conscon,+}^r } \int_0^t \lambda_k(z_C, Z^{(n)}_S(s))  \occ_n(z_C\times  \differential{s}) \\
    &\quad \nrt 
    \sum_{z_C \in \BddZ_{\conscon,+}^r }\int_0^t \lambda_k(z_C, Z_S(s))  \pi_{Z_S(s)}(z_C) \differential{s} \equiv \int_0^t\lambda_k^{\mrm{avg}}(Z_S(s)) \differential{s}.
\end{aligned}
\end{align}
It  follows from \eqref{eq:Z-split-eq}, \eqref{eq:Z-split-det}, together with \eqref{eq:mart-L2-conv} and \eqref{eq:weak-conv-occ-meas}, that the limit point $Z_V$ is a solution of the \ac{ODE} \eqref{eq:ZV-limit-ODE}.

Finally, note that the mapping $z_S \in [0,\infty) \mapsto p_1(z_S) \in [0,1]$ is Lipschitz continuous, and hence, so are all the mappings $z_S \in [0,\infty) \mapsto p_i(z_S) \in [0,1]$, $i=1,2,3\hdots,r$ (see \Cref{lem:fast-stat-dist}). Thus, the \ac{ODE} for $Z_V$ admits a unique solution, and hence the limit point $(\pi_{Z_S}\star \leb, Z_V)$ is unique.

\end{proof}

    One can verify directly from the first equation in \eqref{eq:_pr_recursions} in \Cref{sec:stationary_distribution} that $\timeDerivative{Z_S} + \timeDerivative{Z_P} =0$. We remark that an \ac{sQSSA} of the type of \Cref{thm:FLLN} can indeed be seen as a model reduction of the original multi-stage \ac{MM} system in \eqref{eq:mm_det} into a simple \ac{CRN} of the form $S \xrightharpoonup[]{} P$. Interestingly, the effective propensity function $h$ has the same form as the standard \ac{MM} propensity function ($r=1$ case), although the explicit expression cannot be intuitively guessed from the simpler \ac{MM} model.

\begin{figure}[t]
        \centering
        \includegraphics[width=0.49\textwidth]{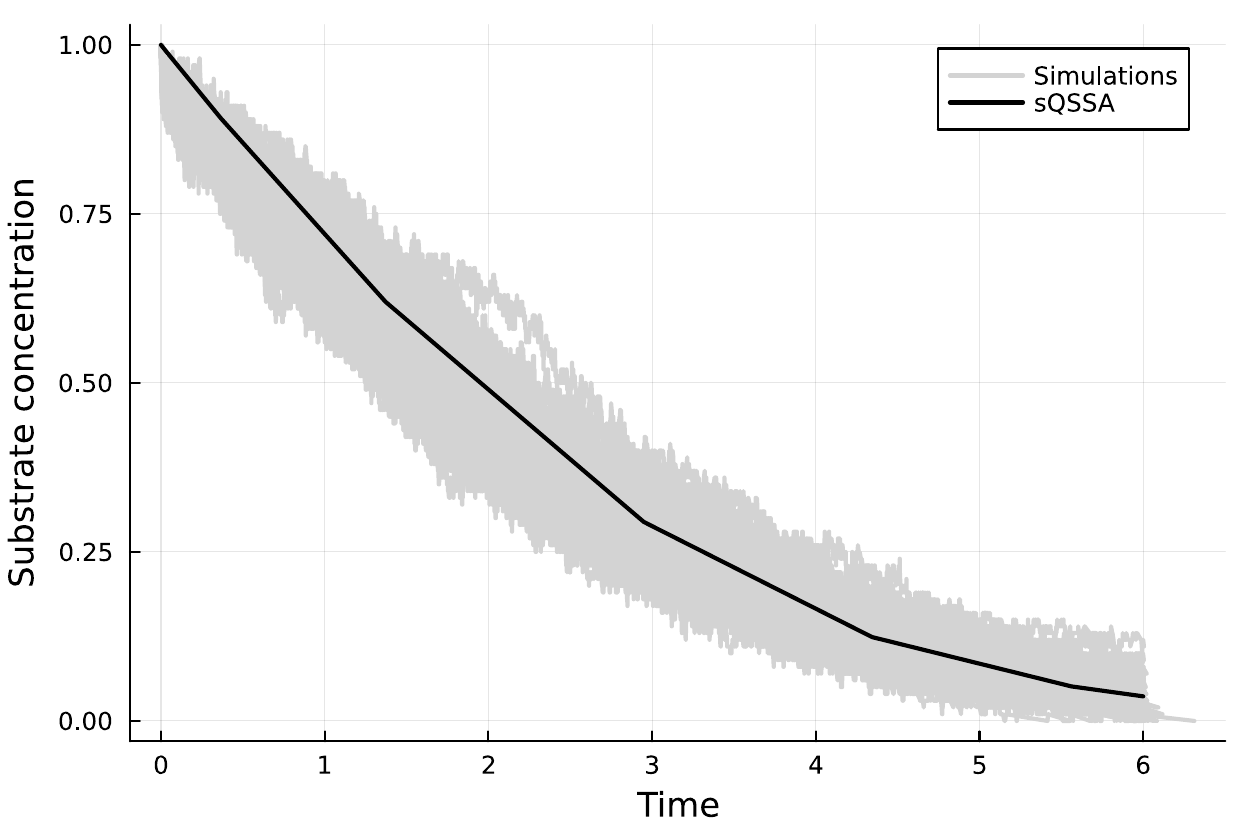}
        \includegraphics[width=0.49\textwidth]{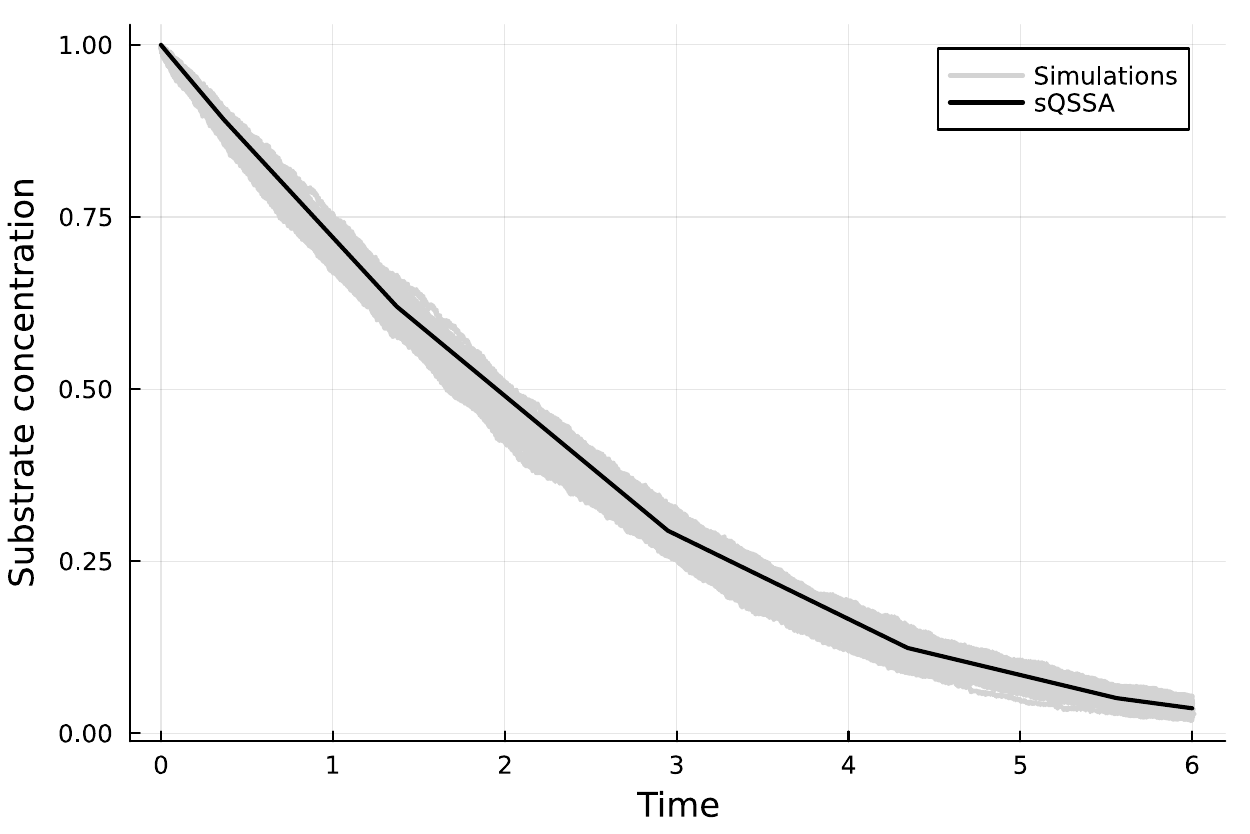}

        \caption{The accuracy of the \ac{sQSSA} for the multi-stage \ac{MM} reaction system in \Cref{example:two_species}. We compare the deterministic \ac{ODE} with 100 trajectories of Doob--Gillespie simulations of the original stochastic model.  (Left) $n=100$. (Right) $n=1000$. Other parameters are $J=10, \kappa_1 = 1, \kappa_{-1} = 1, \kappa_2 = 1, \kappa_{-2} = 1, \kappa_P = 0.1$. 
        }
        \label{fig:2_species_sQSSA}
    \end{figure}

    \begin{myExample}
        \label{example:two_species}
        Consider the case of two intermediate complex species $C_1, C_2$, \ie, when $r=2$. In this case, we have the following reactions
        \begin{equation}
            S + E \xrightleftharpoons[k_{-1}]{k_1} C_1 \xrightleftharpoons[k_{-2}]{k_2}  C_2 \xrightharpoonup[]{k_P}    P+ E\eqcomma  
            \label{eq:two_species_reactions}
        \end{equation}
        with the corresponding stochastic rates denoted by $\kappa_1, \kappa_{-1}, \kappa_2, \kappa_{-2}, \kappa_P$. This system has been referred to as ``Scheme 2'' in \cite{Srinivasan2021MM}. Then, the probabilities of the $\sys{Multinomial}(\conscon, p_1(z_S), p_2(z_S))$ distribution are given by
        \begin{align}
            \begin{aligned}
                p_1(z_S) &= \frac{z_S \kappa_1 (\kappa_{-2}+\kappa_P)}{ (\kappa_{-2}+\kappa_P)(\kappa_{-1} + \kappa_1 z_S) + \kappa_2 (\kappa_1 z_S + \kappa_P) }\eqcomma \quad  
             p_{2}(z_S) 
             = \frac{z_S \kappa_1 \kappa_2}{(\kappa_{-2}+\kappa_P)(\kappa_{-1} + \kappa_1 z_S) + \kappa_2 (\kappa_1 z_S + \kappa_P) } \eqstop 
            \end{aligned}
            \label{eq:two_species_probabilities}
        \end{align}
        Note that the probabilities of the $\sys{Multinomial}$ distribution are not easy to generalize from the standard  case of \ac{MM} reaction system with $r=1$ intermediate complex (\ie, when $C_1$ is the only intermediate complex), which has been considered previously in the literature, \eg, \cite{Kang:2014:CLT,Ball:2006:AAM,Enger2023Unified}. See \Cref{fig:2_species_sQSSA} for a comparison of the \ac{sQSSA} with Doob--Gillespie trajectories.  
    \end{myExample}


    \section{\acl{IPS} and Statistical Inference}
    \label{sec:IPS}

   One of the main practical advantages of \acp{QSSA} for \ac{MM} enzyme kinetic reaction systems (see \Cref{thm:FLLN}) is the resulting reduction in model complexity through averaging of the fast intermediate species. In particular, \Cref{thm:FLLN} shows that the original multi-stage \ac{MM} reaction system in \eqref{eq:mm_det} can be approximated by an effective single conversion reaction of the form
\begin{align}\label{eq:h_defn}
    S \xrightharpoonup[]{h_{\theta}(S)} P,\quad \text{where} \quad h_{\theta}(y) =  \f{\conscon (\kappa_1 c_2-\kappa_{-1})y}{c_1y+c_2} 
    = \f{\conscon c_3y}{c_1y+c_2} \equiv \f{\theta_1 y}{\theta_2+y},
\end{align}
with the effective parameter $\theta$ of the reduced model given by
\begin{align}
    \label{eq:para-red}
    \theta \equiv (\theta_1,\theta_2)  \defeq (\conscon c_3/c_1,\ c_2/c_1). 
\end{align}
This simultaneous reduction in both process complexity and effective parameter dimension yields substantial computational advantages. 
In this section, we  develop a rigorous mathematical framework for inference of the parameter $\theta$ from {\em data on product formation times}.

    This estimation problem is challenging and  different from traditional statistical inference of dynamical systems, which requires data on the states of the system at different time points --- in this case, paths of the process $Z_V = (Z_S, Z_P)$.
     When complete trajectories are available, one could construct a likelihood function in a  straightforward manner from the Doob--Gillespie simulation algorithm (see \cite{Wilkinson2018SMS,Anderson2015book} for example) or as an appropriate likelihood ratio (as a Radon--Nikodym derivative; see \cite[Appendix 1, Proposition 2.6, p. 320]{Kipnis1999ScalingLimits}, \cite{Kutoyants2023Poisson}). However, the assumption of having access to complete trajectories of $Z_V$, or even  high-frequency observations on $Z_V$ is often unrealistic. Instead, often in practice, we only have access to a limited amount of data in the form of a random sample $t_1, t_2, \ldots, t_K$ of times of product formation over a fixed time interval $[0,T]$. This is similar to the situation in infectious disease epidemiology when one has access to only a random sample of infection or removal times (death or recovery). In the context of compartmental epidemic models, such times have been termed \emph{transfer times} \cite{DiLauro2022Interface,KhudaBukhsh2020Focus,KhudaBukhsh2023JTB,KhudaBukhsh2024Howto,KhudaBukhsh2023Prison} since they refer to a transfer of an entity from one compartment to another. The lack of information on the system's state, specifically the values of $Z_V$ at $t_1, t_2, \ldots, t_K$, makes it impossible to estimate $\theta$ directly from \eqref{eq:ZV-limit-ODE} using conventional methods like least squares. It is, therefore, crucial to develop a novel statistical inference methodology that is designed to work with datasets consisting of only a sample of product formation times. 
     
     We achieve this through the construction of a suitable  \ac{IPS}. 
          This part of the paper develops several important asymptotic results for this \ac{IPS} providing the foundation to construct and prove the consistency of the estimator of $\theta$.
     

\subsection{Construction of the \ac{IPS} and its asymptotic properties}
\label{sec:prop_chaos}

Consider a particle-system of $n$ substrate molecules $S_1, S_2, \hdots, S_n$ with reactions
$
        S_i \rightharpoonup P_i\eqcomma \quad i=1, 2, \ldots, n\eqstop 
        $
The propensity of these reactions is determined by the total number of substrate molecules present in the system.  The evolution of the system will be captured by the process $\subst^{(n)} = (\subst^{(n)}_1,\subst^{(n)}_2,\hdots,\subst^{(n)}_n)$, where $\subst^{(n)}_i(t) =1$ if the $S_i$-molecule has not undergone a conversion into a product molecule $P_i$ molecule by time $t$ (that is, the  $i$-th particle is still in $S$ state at time $t$), and $\subst^{(n)}_i(t) =0$, otherwise. Once the particle has undergone a conversion into a product molecule ($P$ state), the particle stays in that state forever. Denote by
  \begin{align}
        \bar{\subst}^{(n)}(t) \defeq \frac{1}{n} \sum_{i=1}^{n} \subst_i^{(n)}(t),  \label{eq:barQ_defn} 
    \end{align}
 the proportion of substrate molecules ($S$-state) at time $t$. Theorem \ref{thm:FLLN} suggests that the macro-level propensity function of the reaction system is given by $n h_{\theta}(\bar\subst^{(n)}(t))$, where $h_{\theta}$ is given by \eqref{eq:h_defn}. In other words, borrowing language from infectious disease epidemiology, $n h_{\theta}(\bar \subst^{(n)}(t))$ is interpreted as the total instantaneous pressure on $S$ molecules to convert into $P$ molecules at time $t$. Thus, the quantity, 
 $h_{\theta}(\bar \subst^{(n)}(t))/\bar \subst^{(n)}(t)$ denotes the ``per molecule'' instantaneous pressure on any specific $i$-th reaction. Consequently, the evolution of $\subst^{(n)}$ can be described by the following weakly \acl{IPS} of \acp{SDE} driven by \acp{PRM} (written in integral form)   
    \begin{align}
        \subst_i^{(n)}(t) = 1- \int_{[0, t]\times [0, \infty)} \indicator{[0, \subst_i^{(n)}(s-) g_{\theta}(\bar{\subst}^{(n)}(s-)) ]}{v} \eta_i\left(\differential{v}\times \differential{s}\right)\eqcomma 
        \label{eq:particle_i_trajectory}
    \end{align}
    where $\eta_1, \eta_2, \ldots, $ are independent \acp{PRM} on $\setOfPositiveReals\times \setOfPositiveReals$ with intensity measure $\leb \ot \leb$, and the function $g_\theta: [0,\infty) \rt [0,\infty)$ is defined by 
   \begin{align}
       \label{eq:avg-prop}
     g_\theta(y) \defeq 
        h_\theta(y)/y \equiv \theta_1 /(\theta_2+y).  
   \end{align}
For $i=1, 2, \ldots, n$, let 
    \begin{align}
        \tau_i^{(n)} \defeq \inf\{t \ge 0 \mid \subst^{(n)}_i(t) = 0\}\eqcomma \label{eq:tau_i_defn}
    \end{align}
   denote the time of conversion of the $i$-th particle from $S$ to $P$ state. Since
\begin{align} \label{eq:tau-subst-rel}
  \{\tau_i^{(n)} >t\} = \{\subst^{(n)}_i(t)=1\}  \eqcomma 
\end{align}
the system can be equivalently described by the collection of stopping times $\{\tau_i^{(n)}: i=1,2,\hdots,n\}$. The random variables $\tau^{(n)}_i$, which track the product formation times, are important for formulating a suitable likelihood function for the parameter $\theta$.

The distribution of the $\tau^{(n)}_i$ does not admit a tractable expression; so, we consider an approximate likelihood function based on its limiting form. Before introducing this likelihood, we first establish the necessary convergence results for the $\subst^{(n)}$-system, which are not only interesting in their own right but also play a crucial role in proving the consistency of the estimator.
 
It is intuitively clear that $\subst^{(n)}$ approximates the reduced order model of \Cref{thm:FLLN} in the sense that  as the number of substrate molecules, $n$, goes to infinity,  the mean process $\bar{\subst}^{(n)}$ converges to $Z_S$, the solution of the system of \acp{ODE} in \eqref{eq:ZV-limit-ODE}, with initial condition $Z_S(0)=1$. \Cref{prop:FLLN_Tbar} provides the $L^p$-convergence rate for this process, while the following result, \Cref{th:conc-ineq}, establishes a non-asymptotic concentration inequality.  

For the parameter estimation, it is necessary to denote the  probability measure $\PP$ on  $(\Om, \cal{F})$ by $\PP_{\theta}$ to emphasize its dependence on the  parameter $\theta$. Thus, under $\PP_{\theta}$, the process $\subst^{(n)}$ satisfies \eqref{eq:particle_i_trajectory} with parameter $\theta$.

\begin{myProposition}
        For any $p\geq 2$,
          \begin{align}
            \EE_\theta\lf[\sup_{t\le T} \absolute{{\bar{\subst}^{(n)}}(t) - Z_{S,\theta}(t)}^p\ri] \leq \cnst_p  \lf(\f{\theta_1 T}{1+\theta_2} \wedge 1\ri)\exp\lf(p\ \theta_1 T/\theta_2\ri) n^{-p/2}, 
        \end{align}
        where $Z_S \equiv Z_{S,\theta}$ is the unique solution of the \ac{ODE}  
        \begin{align}
            \timeDerivative{Z_{S,\theta}(t)} = -h_{\theta}(Z_{S,\theta}(t))\eqcomma \quad Z_{S,\theta}(0)=1,
            \label{eq:barQ_ODE}
        \end{align}  
        and  $\cnst_p $ is the \ac{BDG}-constant.
        \label{prop:FLLN_Tbar}
    \end{myProposition}
    \begin{proof}[Proof of \Cref{prop:FLLN_Tbar}]
        Since $g_\theta(y) = h_\theta(y)/y$,  the process $\bar{\subst}^{(n)}$ satisfies 
        \begin{align*}
            {\bar{\subst}^{(n)}}(t) = 1 - \frac{1}{n}\sum_{i=1}^{n}\int_{[0, t]\times [0, \infty)} \indicator{[0, \subst_i^{(n)}(s-) g_{\theta}(\bar{\subst}^{(n)}(s-)) ]}{v} \bar{\eta}_i\left(\differential{v}\times \differential{s}\right) -  \int_0^t h_{\theta}(\bar{\subst}^{(n)}(s))\differential{s} \eqcomma 
        \end{align*}
        where $\bar{\eta}_i$ is the centered \ac{PRM} corresponding to the \ac{PRM} $\eta_i$. Therefore, 
        \begin{align*}
            \absolute{{\bar{\subst}^{(n)}}(t) - Z_S(t)} &{} \le \frac{1}{n}\absolute{\sum_{i=1}^{n}\int_{[0, t]\times [0, \infty)} \indicator{[0, \subst_i^{(n)}(s-) g_{\theta}(\bar{\subst}^{(n)}(s-)) ]}{v} \bar{\eta}_i\left(\differential{v}\times \differential{s}\right)   } 
            + \int_0^t \absolute{h_{\theta}(\bar{\subst}^{(n)}(s)) - h_{\theta}(Z_S(s)) }\differential{s}\\
            &{} \le \frac{1}{n}\absolute{\sum_{i=1}^{n}\int_{[0, t]\times [0, \infty)} \indicator{[0, \subst_i^{(n)}(s-) g_{\theta}(\bar{\subst}^{(n)}(s-)) ]}{v} \bar{\eta}_i\left(\differential{v}\times \differential{s}\right)   } 
            + (\theta_1/\theta_2) \int_0^t \absolute{\bar{\subst}^{(n)}(s) - Z_S(s) } \differential{s}\eqcomma 
        \end{align*}
       where we used the fact that $h_\theta$ is Lipschitz continuous with Lipschitz constant, $\|h\|_{\mrm{Lip}} = \sup_{y\geq 0}|h'(y)| \leq \theta_1/ \theta_2$. Then, by Gr\"onwall's inequality, and taking supremum on both sides, we have 
        \begin{align}
            \sup_{t\le T} \absolute{{\bar{\subst}^{(n)}}(t) - Z_S(t)} \le \frac{1}{n}\left( \sup_{t\le T }  \absolute{\SC{E}_n(t)} \right) \exp\lf((\theta_1/\theta_2) T\ri)\eqcomma \label{eq:Qbar_diff_bound}
        \end{align}
        where $\SC{E}_n$ defined by
        \begin{align*}
            \SC{E}_n(t) \defeq {\sum_{i=1}^{n}\int_{[0, t]\times [0, \infty)} \indicator{[0, \subst_i^{(n)}(s-) g_{\theta}(\bar{\subst}^{(n)}(s-)) ]}{v} \bar{\eta}_i\left(\differential{v}\times \differential{s}\right)   }\eqstop 
        \end{align*}
        is a martingale. Notice that the quadratic variation $[\SC{E}_n]$ satisfies
        \begin{align*}
         [\SC{E}_n]_t = n(1-\bar{\subst}^{(n)}(t)) \leq n, \quad \EE([\SC{E}_n]_t) \leq n \EE\int_0^t h_{\theta}(\bar{\subst}^{(n)}(s))\differential{s} \leq \lf(\f{\theta_1 t}{1+\theta_2} \wedge 1\ri) n   
        \end{align*}
      The above inequality used the fact that $\bar{\subst}^{(n)}$ takes values in $[0,1]$, and $0 \leq h(y) \leq \theta_1/(1+\theta_2)$ for $y \in [0,1]$.   
        
        By the \ac{BDG} inequality, we have for any $p\geq 2$,
        \begin{align*}
            \EE\left[\sup_{t\le T} \absolute{\SC{E}_n(t)}^p\right] &{}\le \cnst_p {\EE\left[[\SC{E}_n]^{p/2}_T\right]} \leq \cnst_p n^{p/2-1} \EE\left([\SC{E}_n]_T \right)           \leq  \cnst_p  \lf(\f{\theta_1 T}{1+\theta_2} \wedge 1\ri) n^{p/2}. 
        \end{align*}
         The assertion now easily follows from \eqref{eq:Qbar_diff_bound}.
        \end{proof}

The following concentration inequality shows that the probability of the fluctuations between $\bar{\subst}^{(n)}(t)$ and its limit $Z_{S,\theta}(t)$ exhibits an exponential decay with an $n^{1/2}$ pre-factor, which, in particular, leads to a large deviation upper bound (see \Cref{rem:FLLN_Tbar}). 

      \begin{myTheorem}\label{th:conc-ineq}
 Let $\theta \in \Theta \subset (0, \infty)^2$ and $T>0$. Then, for any $t\geq 0$, $0<\ep<\|h_\theta\|_{\infty,1} t$, 
 $$\PP_{\theta}\lf(\absolute{{\bar{\subst}^{(n)}}(t) - Z_{S,\theta}(t)}>\ep\ri) \leq 2\lf(1+ \frac{\ep }{\|h\|_{\infty,1}}\const_{\mrm{conc}}(\theta,t) n^{1/2}\ri) \myExp{-I_{t,\theta}(\ep)n},$$
 where $Z_{S,\theta}$ is the unique solution of \eqref{eq:barQ_ODE},  $$\|h_\theta\|_{\infty,1} = \sup_{ y \in [0,1]}h_\theta(y) = \theta_1/(\theta_2+1),$$
 and for $u \in(-\|h_\theta\|_{\infty,1} t, \infty), \ u \neq 0$,
\begin{align}
    \label{eq:barY-conc-rate}
    \begin{aligned}
I_{t,\theta}(u) & \defeq \|h_\theta\|_{\infty,1} t\lf(\Big(1+\frac{u}{\|h_\theta\|_{\infty,1} t}\Big)\ln\Big(1+\frac{u}{\|h_\theta\|_{\infty,1} t}\Big)- \frac{u}{\|h_\theta\|_{\infty,1} t} \ri) > 0,  \\ 
\const_{\mrm{conc}}(\theta, t)& \defeq (\theta_1/\theta_2)\cnst^{1/2}_2 \lf(\f{\theta_1 t}{1+\theta_2}\wedge 1\ri)^{1/2}\exp\lf((\theta_1/\theta_2) t\ri).
\end{aligned}
\end{align}
Furthermore, for $0<\ep<\|h_\theta\|_{\infty,1} T$,
\begin{align}\label{eq:conc-const}
\PP_{\theta}\Bigg(\sup_{t\le T} & \absolute{{\bar{\subst}^{(n)}}(t) - Z_{S,\theta}(t)}>\ep\Bigg) \leq 2 \lf(1+ \tilde \const_{\mrm{conc}}(\theta,T,n, \ep)\ri)\myExp{-I_{T,\theta}(\ep)n}.
\end{align}
where 
\begin{align*}
\tilde\const_{\mrm{conc}}(\theta,T,n, \ep) \defeq \f{|\ln(1-\frac{\ep}{\|h_\theta\|_{\infty,1} T})|\lf(1+ \frac{\ep}{\|h_\theta\|_{\infty,1}}\const_{\mrm{conc}}(\theta,T) n^{1/2}\ri)}{\frac{\ep}{\|h_\theta\|_{\infty,1}T}- \ln(1+\frac{\ep}{\|h_\theta\|_{\infty,1} T})}.
\end{align*}
 \end{myTheorem}

\begin{proof}[Proof of \Cref{th:conc-ineq}] Since $\theta$ is fixed, we suppress it from various notations for convenience. Note that the process $R^{(n)}(t) \defeq n\lf(\bar{\subst}^{(n)}(t) -Z_S(t)\ri)$ satisfies 
        \begin{align*}
            R^{(n)}(t)= n\int_0^t h(Z_S(s)) \differential{s} - \sum_{i=1}^{n}\int_{[0, t]\times [0, \infty)} \indicator{[0, \subst_i^{(n)}(s-) g_{\theta}(\bar{\subst}^{(n)}(s-)) ]}{v} \eta_i\left(\differential{v}\times \differential{s}\right).
        \end{align*}
For $\lambda \in \R$, by It\^o's lemma, $D^{(n)}(\lambda,t) \defeq e^{\lambda R^{(n)}(t)} \equiv e^{\lambda n(\bar{\subst}^{(n)}(t) - Z_S(t))} $ satisfies
 \begin{align}\label{eq:exp-Rn-proc}
 \begin{aligned}
e^{\l R^{(n)}(t)}& = 1+\lambda n \int_0^t e^{\l R^{(n)}(s)} h(Z_S(s)) \differential{s}
 +(e^{-\lambda}-1)\sum_{i=1}^n \int_{[0, t]\times [0, \infty)}e^{\l R^{(n)}(s-)}  \indicator{[0, \subst_i^{(n)}(s-) g(\bar{\subst}^{(n)}(s-)) ]}{u}\eta_i(\differential{u}\times\differential{s})\\
  &=1 + (\lambda+e^{-\lambda}-1)n \int_0^te^{\l R^{(n)}(s)} h(Z_S(s)) \differential{s}
  +(e^{-\lambda}-1)n \int_0^t e^{\l R^{(n)}(s)} \lf(h(\bar{\subst}^{(n)}(s))-h(Z_S(s))\ri) \differential{s}\\
  &\quad 
  +\mart^{(n)}(\lambda,t),
 \end{aligned} 
 \end{align}
where the martingale $\mart^{(n)}(\lambda,\cdot)$ is given by
$$\mart^{(n)}(\lambda,t) = (e^{-\lambda}-1)\sum_{i=1}^n \int_{[0, t]\times [0, \infty)}e^{\l R^{(n)}(s-)}  \indicator{[0, \subst_i^{(n)}(s-) g(\bar{\subst}^{(n)}(s-)) ]}{u}\bar\eta_i(\differential{u}\times\differential{s}).$$
Since $h$ is increasing, $h(\bar{\subst}^{(n)}(s))-h(Z_S(s)) \geq 0$ on $\{R^{(n)}(s) \geq 0\}$ and $h(\bar{\subst}^{(n)}(s))-h(Z_S(s)) \leq 0$ on $\{R^{(n)}(s) \leq 0\}$. 
Now splitting the integrand in the third term according as $\{R^{(n)}(s) \geq 0\}$ and $\{R^{(n)}(s) \leq 0\}$, it follows for any $\l \geq 0$ (after noting $\l \geq 0$ implies $e^{-\l}-1 \leq 0$) that
\begin{align*}
 (e^{-\lambda}-1) e^{\l R^{(n)}(s)}& \lf(h(\bar{\subst}^{(n)}(s))-h(Z_S(s))\ri) \\
 & \leq  (e^{-\lambda}-1) e^{\l R^{(n)}(s)} \lf(h(\bar{\subst}^{(n)}(s))-h(Z_S(s))\ri) \indic{\{R^{(n)}(s) \leq 0\}} \\
 & \leq |e^{-\lambda}-1| (\theta_1/\theta_2)|\bar{\subst}^{(n)}(s) - Z_S(s)|.
\end{align*}
where in the last step we used Lipschitz continuity of $h$ with $\|h\|_{\mrm{Lip}} = \theta_1/\theta_2$,  and the fact that for $\l \geq 0$, $e^{\l R^{(n)}(s)} \leq 1$ on $\{R^{(n)}(s) \leq 0\}$. Similarly, for $\l \leq 0,$
\begin{align*}
 (e^{-\lambda}-1) e^{\l R^{(n)}(s)}& \lf(h(\bar{\subst}^{(n)}(s))-h(Z_S(s))\ri) \\
 & \leq  (e^{-\lambda}-1) e^{\l R^{(n)}(s)} \lf(h(\bar{\subst}^{(n)}(s))-h(Z_S(s))\ri) \indic{\{R^{(n)}(s) \geq 0\}} \\
 & \leq (e^{-\lambda}-1) (\theta_1/\theta_2)|\bar{\subst}^{(n)}(s) - Z_S(s)|.
\end{align*}

It follows that for any $\l \in \R$  (after noting that $\lambda+e^{-\lambda}-1 \geq 0$),
\begin{align*}
    \EE[  e^{\l R^{(n)}(t)}] & \leq  1+(\lambda+e^{-\lambda}-1) \|h\|_{\infty,1} n \int_0^t  \EE(e^{\l R^{(n)}(s)})  \differential{s}
    + |e^{-\lambda}-1| (\theta_1/\theta_2)tn\ \sup_{s\leq t} \EE[|\bar{\subst}^{(n)}(s) - Z_S(s)|]\\
    & \leq 1+(\lambda+e^{-\lambda}-1) \|h\|_{\infty,1} n \int_0^t  \EE(e^{\l R^{(n)}(s)})  \differential{s}\\ 
  &{}  + |e^{-\lambda}-1| (\theta_1/\theta_2) \cnst^{1/2}_2\lf(\f{\theta_1 t}{1+\theta_2}\wedge 1\ri)^{1/2}\exp\lf((\theta_1/\theta_2) t\ri) tn^{1/2},
\end{align*}
where for the last inequality we used \Cref{prop:FLLN_Tbar}.
By Gr\"onwall's inequality it follows that for any $\lambda \in \R$
\begin{align}
    \label{eq:barY-exp-mmt-bd}
    \begin{aligned}
    \EE[ e^{\l R^{(n)}(t)}] \leq&\ \lf(1+ t|e^{-\lambda}-1|\const_{\mrm{conc}}(\theta,t) n^{1/2}\ri)  \myExp{(\lambda+e^{-\lambda}-1) \|h\|_{\infty,1} n t}.
    \end{aligned}
\end{align}
where $\const_{\mrm{conc}}(\theta, t)$ is defined in \eqref{eq:barY-conc-rate}.
 Now, notice that  
 $$\PP(\absolute{{\bar{\subst}^{(n)}}(t) - Z_{S}(t)}>\ep) \leq \mathsf{P}^{(n),+}_\ep(t)+\mathsf{P}^{(n),-}_\ep(t),$$ 
 where 
 \begin{align*}
   \mathsf{P}^{(n),+}_\ep(t) &\defeq \PP\lf(\big(\bar{\subst}^{(n)}(t) - Z_{S}(t)\big)>\ep\ri) = \PP\lf(R^{(n)}(t) > n \ep\ri) ,\\
   \mathsf{P}^{(n),-}_\ep(t) &\defeq \PP\lf(-(\bar{\subst}^{(n)}(t) - Z_{S}(t))>\ep\ri)=\PP\lf(-R^{(n)}(t) > n \ep\ri).   
 \end{align*}
By Markov's inequality and \eqref{eq:barY-exp-mmt-bd}, for any $\lambda'>0$ 
\begin{align*}
   \mathsf{P}^{(n),+}_\ep(t) & = \PP\lf(\lambda' R^{(n)}(t)>\lambda' n \ep\ri)\ =  \PP\lf(e^{\lambda' R^{(n)}(t)}>e^{\lambda' n \ep}\ri)\\
   & \leq \lf(1+ t|e^{-\lambda'}-1|\const_{\mrm{conc}}(\theta, t) n^{1/2}\ri) 
   \myExp{\lf((\lambda'+e^{-\lambda'}-1) \|h\|_{\infty,1} t - \lambda' \ep\ri)n}.
\end{align*}
Notice that the exponent in the exponential term is minimized for $\lambda'_{+,*} = - \ln\lf(1- \frac{\ep}{\|h\|_{\infty,1} t}\ri)>0$, with  $0<\ep<  \|h\|_{\infty,1} t$. Plugging $\lambda'_{+,*}$ in the above equation, we get for $0<\ep<  \|h\|_{\infty,1} t$, 
\begin{align*}
\mathsf{P}^{(n),+}_\ep(t) \leq \lf(1+\frac{ \ep }{\|h\|_{\infty,1}}\const_{\mrm{conc}}(\theta, t)n^{1/2}\ri)\myExp{-I_t(-\ep)n}\eqcomma 
\end{align*}
where $I_t(\cdot)$ is defined by \eqref{eq:barY-conc-rate}. 
 
Similarly, it follows that for any $\lambda'>0$,
$$\mathsf{P}^{(n),-}_\ep(t) \leq \lf(1+ t(e^{\lambda'}-1)\const_{\mrm{conc}}(\theta,t) n^{1/2}\ri)\myExp{\lf((-\lambda'+e^{\lambda'}-1) \|h\|_{\infty,1} t - \lambda' \ep\ri)n}.$$ 
The exponent in the exponential term is minimized for $\lambda'_{-,*} = \ln \left( 1 + \frac{\epsilon}{\|h\|_{\infty,1} t} \right)>0$, and plugging this in the above inequality shows that for any $\ep>0$, 
$$\mathsf{P}^{(n),-}_\ep(t) \leq \lf(1+ \frac{\ep}{\|h\|_{\infty,1}}\const_{\mrm{conc}}(\theta,t) n^{1/2}\ri) e^{-I_t(\ep)n}.$$
The assertion follows because $I_t(\ep)\leq I_t(-\ep)$ for any $0<\ep<\|h\|_{\infty,1} t$.

Next, note that from \eqref{eq:exp-Rn-proc}
\begin{align*}
    \EE\lf(\sup_{0\leq t\leq T} e^{\lambda R^{(n)}(t)}\ri) \leq
    \begin{cases}
        1+\lambda n \int_0^T \EE\lf(e^{\lambda R^{(n)}(s)} h(Z_S(s))\ri) \differential{s},& \quad \lambda\geq 0\eqcomma \\
        1+(e^{-\lambda}-1)n \int_0^T \EE\lf(e^{\lambda R^{(n)}(s)} h(\bar{\subst}^{(n)}(s))\ri) \differential{s},& \quad \l \leq 0.\\
    \end{cases}
 \end{align*}
It follows from \eqref{eq:barY-exp-mmt-bd} that for any $\l \in \R -\{0\}$
\begin{align*}
    \EE\lf(\sup_{0\leq t\leq T} e^{\lambda R^{(n)}(t)}\ri) \leq&\ 1 +(\lambda \vee (e^{-\lambda}-1))\|h\|_{\infty,1}n \int_0^T \EE\lf(e^{\lambda R^{(n)}(s)}\ri) \differential{s}\\
  \leq&\ 1 +\f{(\lambda \vee (e^{-\lambda}-1))\|h\|_{\infty,1}n\lf(1+ T|e^{-\lambda}-1|\const_{\mrm{conc}}(\theta,T) n^{1/2}\ri)}{(\lambda+e^{-\lambda}-1) \|h\|_{\infty,1}n}   
  \lf(e^{(\lambda+e^{-\lambda}-1) \|h\|_{\infty,1} n T} -1\ri)\\
  \leq&\ \lf(1 +\f{(\lambda \vee (e^{-\lambda}-1))\lf(1+ T|e^{-\lambda}-1|\const_{\mrm{conc}}(\theta,T) n^{1/2}\ri)}{\lambda+e^{-\lambda}-1}\ri)   
  e^{(\lambda+e^{-\lambda}-1) \|h\|_{\infty,1} n T}.
\end{align*}
Noting that
\[
\PP_{\theta}\left(\sup_{t\le T} \absolute{{\bar{\subst}^{(n)}}(t) - Z_{S,\theta}(t)} > \ep\right) \leq \PP\left(\sup_{t\le T} R^{(n)}(t) > n\ep\right) + \PP\left(\sup_{t\le T} (-R^{(n)}(t)) > n\ep\right) \eqstop 
\]
Applying Markov's inequality, we obtain for any $\lambda' > 0$
\begin{align*}
 \PP\left(\sup_{t\le T} R^{(n)}(t) > n\ep\right) &\leq e^{-\lambda' n\ep} \EE\left(\sup_{0\leq t\leq T} e^{\lambda' R^{(n)}(t)}\right), \quad
 \PP\left(\sup_{t\le T} (-R^{(n)}(t)) > n\ep\right) \leq e^{-\lambda' n\ep} \EE\left(\sup_{0\leq t\leq T} e^{-\lambda' R^{(n)}(t)}\right).
\end{align*}
The remainder of the proof follows as before.



\end{proof}

\begin{myLemma}
    \label{lem:Z-prob-meas}  The function $Z_S \equiv Z_{S,\theta}$ in  \eqref{eq:barQ_ODE} determines a probability measure $\taudist_\theta$ on $[0,\infty)$ defined by
     \begin{align*}
        \taudist_\theta([0,t]) \defeq&\ 1 - Z_{S,\theta}(t) \equiv 1 - \exp\lf(-\int_0^t g_\theta(Z_{S,\theta}(s)) \differential{s}\ri), \quad 0< t<\infty.
    \end{align*}  
\end{myLemma}

\begin{proof}[Proof of \Cref{lem:Z-prob-meas}]
    We will show that the mapping $t \in [0,\infty) \mapsto  \taudist_\theta([0,t])$ defines a valid \ac{CDF}. Since $h_{\theta}(y) \geq 0$ for $y \geq 0$, it is clear that $ Z'_{S,\theta}(t) \leq 0$ and therefore, $Z_{S,\theta}$ is decreasing with values in $[0,1]$; equivalently,  the mapping $t \in [0,\infty) \Rt \taudist_\theta([0,t]) = 1 - Z_{S,\theta}(t)$ is increasing with values in $[0,1]$, and it is clearly continuous in $t$. 
    
    We now need to show that $\lim_{t \rt \infty}\taudist_\theta([0,t]) =1$. To this end,  notice that \eqref{eq:h_defn} shows $h_\theta(y) \geq \conscon c_3 y/ (c_1+c_2)$ for $y \in [0,1]$. Therefore, 
  \begin{align*}
     \f{\differential{Z_{S,\theta}(t)}}{\differential{t}} \leq -  \f{\conscon c_3 Z_{S,\theta}(t)}{c_1+c_2}, \quad Z_{S,\theta}(0) =1,
 \end{align*}
 and hence by the comparison principle for \acp{ODE},
 \begin{align}\label{eq:Z-lbd}
 Z_{S,\theta}(t) \leq \exp\lf(-\f{\conscon c_3}{c_1+c_2}  t\ri) \equiv \exp\lf(-\f{\theta_1}{1+
 \theta_2}  t\ri).
 \end{align}
 This shows $ Z_{S,\theta}(t) \rt 0$ as $t \rt \infty$ proving the assertion.
\end{proof}

Note that $\taudist_\theta([0,t])$ represents the proportion of substrate molecules converted to product by time $t$ in the limiting system, and  $\taudist_\theta([0,\infty))=1$ reflects that the conversion process is complete by time infinity. While this holds for the \ac{MM} system, it need not be true for other systems, where the corresponding $\taudist_\theta(\{\infty\})$ may be strictly positive. For example, in the standard \ac{SIR} model in infectious disease epidemiology, the limiting proportion of susceptible individuals $S$ plays the role of $Z_S$. However, it is straightforward that $\taudist_\theta(\{\infty\}) = S(\infty)>0$, which denotes the proportion of susceptible individuals who escape infection (and are never infected). See \cite{KhudaBukhsh2020Focus,DiLauro2022Interface}.

\begin{myRemark}
    \label{rem:FLLN_Tbar}
Let $\Xi^{(n)} \defeq \f{1}{n}\sum_{i=1}^n \delta_{\tau^{(n)}_i}$ be the empirical measure of the $\tau^{(n)}_i$, where the $\tau_i^{(n)}$ are defined in \eqref{eq:tau_i_defn}. Notice that $\Xi^{(n)}([0,t]) = 1-\bar Y^{(n)}(t)$; therefore, for any $\theta \in \Theta$,  $T>0$ and $0<\ep<\|h_\theta\|_\infty T$,
\begin{align}
    \label{eq:emp-meas-conc}
    \begin{aligned}
    \EE_{\theta}\lf[\sup_{t\le T} \absolute{\Xi^{(n)}([0,t]) - \taudist_\theta([0,t])}^p\ri] \leq \cnst_p \lf(\f{\theta_1 T}{1+\theta_2} \wedge 1\ri)\exp\lf(p(\theta_1/\theta_2) T\ri) n^{-p/2}\eqcomma \\
     \PP_{\theta}\lf(\sup_{t\le T} \absolute{\Xi^{(n)}([0,t]) - \taudist_\theta([0,t])}>\ep\ri) \leq  2 \lf(1+ \tilde \const_{\mrm{conc}}(\theta,T,n,\ep)\ri)\myExp{-I_{T,\theta}(\ep)n},
\end{aligned}  
\end{align} 
where $\tilde \const_{\mrm{conc}}(\theta,T,n) = O(n^{1/2})$ is as in \eqref{eq:conc-const}.
The second inequality immediately shows that the following large deviation upper bound holds
\begin{align*}
\limsup_{n \rt \infty}\f{1}{n}\ln \PP_{\theta}\lf(\sup_{t\le T} \absolute{\Xi^{(n)}([0,t]) - \taudist_\theta([0,t])}>\ep\ri) \leq -I_{T,\theta}(\ep).
\end{align*}

Furthermore, by the Borel--Cantelli lemma,  $\bar{\subst}^{(n)}$ converges to $Z_S$, $\PP_{\theta}$-a.s.  in the following sense: for any $T>0$,
\begin{align}\label{eq:Yn-mean-conv}
\begin{aligned}
\sup_{t\le T} \absolute{\Xi^{(n)}([0,t]) - \taudist_\theta([0,t])} \equiv&\  \sup_{t\le T} \absolute{\bar{\subst}^{(n)}(t) - Z_{S,\theta}(t)} \nrt 0 \quad \PP_{\theta}\text{-a.s.}
\end{aligned}
\end{align}
and by \Cref{prop:FLLN_Tbar}, since the empirical measure of $\subst^{(n)}(t)$, $\f{1}{n} \sum_{i=1}^n \delta_{\subst^{(n)}_i(t)}$, is a probability measure on $\{0,1\}$, it is clear from \eqref{eq:Yn-mean-conv} that 
$$\sup_{t\le T}\lf\|\f{1}{n} \sum_{i=1}^n \delta_{\subst^{(n)}_i(t)} -(Z_{S,\theta}(t) \delta_1+(1-Z_{S,\theta}(t))\delta_0)\ri\|_{\mrm{TV}} \nrt 0, \quad  \PP_{\theta}\text{-a.s. and in } L^1(\PP_{\theta}).$$
\end{myRemark}

\np
The above observations lead to the following result.

\begin{myCorollary}
    \label{cor:prop-ch-tau}
 Fix $\theta \in \Theta$, $T>0$, and define the class of functions $ \mathscr{F}_{T,1}$ by
\begin{align}
    \label{eq:bdd-lip-T}
     \mathscr{F}_{T,1}=\{f\in \mathrm{Lip}([0,T],\R): \|f\|_{\mathrm{BL}, T} \leq 1\},
\end{align}
 where $\|f\|_{\mathrm{BL}, T} = \|f\|_{\infty,T}+\|f\|_{\mathrm{Lip},T}$ with 
 $$\|f\|_{\infty,T} = \sup_{t\leq T}|f(t)|, \quad \|f\|_{\mathrm{Lip},T} = \sup_{\substack{t,t' \in [0,T] \\ t\neq t'}}|f(t)-f(t')|/|t-t'|.$$
 Then,  for any $p \geq 2$ and $\theta \in \Theta$, $T>0$ and $0<\ep<(T\vee1)\|h_\theta\|_\infty T$,
\begin{align*}
  &\PP_{\theta}\lf(\sup_{f\in \mathscr{F}_{T,1}}\sup_{t\le T} \absolute{\Xi^{(n)}(f\indic{[0,t]}) - \taudist_\theta(f\indic{[0,t]})}>\ep\ri) \leq  2 \lf(1+ \tilde \const_{\mrm{conc}}(\theta,T,n, \ep/(T\vee 1))\ri)\myExp{-I_{T,\theta}(\ep/(T\vee 1))n}\eqcomma \\
 &\EE_{\theta}\lf(\sup_{f\in \mathscr{F}_{T,1}}\sup_{t\le T} \absolute{\Xi^{(n)}(f\indic{[0,t]}) - \taudist_\theta(f\indic{[0,t]})}^p\ri)  \leq \cnst_p \lf(\f{\theta_1 T}{1+\theta_2} \wedge 1\ri) (T \vee 1)^p \exp\lf(p(\theta_1/\theta_2) T\ri) n^{-p/2}\eqstop 
\end{align*}

\end{myCorollary}

\begin{proof}[Proof of \Cref{cor:prop-ch-tau}]
 Fix $f \in \mathscr{F}_{T,1}$. 
Since $f$ is Lipschitz continuous on $[0,T]$, it is differentiable almost everywhere on $[0,T]$ with $ \|f\|_{\mathrm{Lip},T} =\|f'\|_{\infty,T} \defeq \text{ess sup}_{t \leq T} |f'(t)| < \infty$.
Since for any Radon measure $\nu$,
$$\nu(f1_{[0,t]}) \equiv \int_0^t f(s)\nu(\differential{s}) = f(t)\nu([0,t])-\int_0^t f'(s)\nu([0,s])\differential{s}\eqcomma$$
we have 
\begin{align}\label{eq:emp-measure-f-ineq}
\begin{aligned}
|\Xi^{(n)}(f1_{[0,t]})-\taudist_\theta(f1_{[0,t]})| & \leq  |f(t)| \absolute{\Xi^{(n)}([0,t]) - \taudist_\theta([0,t])} 
 + t\|f'\|_{\infty,t} \sup_{s\leq t}  \absolute{\Xi^{(n)}([0,s]) - \taudist_\theta([0,s])} \eqstop 
\end{aligned}
\end{align}
Consequently, 
$$\sup_{t\leq T}|\Xi^{(n)}(f1_{[0,t]})-\taudist_\theta(f1_{[0,t]})| \leq (T \vee 1) \|f\|_{\mathrm{BL}, T} \sup_{t\leq T}\absolute{\Xi^{(n)}([0,t]) - \taudist_\theta([0,t])},$$
and the assertion follows from \eqref{eq:emp-meas-conc}.
\end{proof}

Notice, by the equivalence of convergence of empirical measures and propagation of chaos for an exchangeable system \cite[Proposition~2.2(i)]{Sznitman1991Topics}, it follows that
for any fixed $k$,
\begin{align}
    \label{eq:tau-prop-chaos}
 \bm{\tau}^{(n)}_{1:k} \defeq (\tau^{(n)}_1, \tau^{(n)}_2, \hdots, \tau^{(n)}_k ) \stackrel{n\rt \infty}\LRT \Phi^{\ot k}_{\theta} \equiv \bigotimes_{i=1}^k \Phi_\theta\eqstop 
\end{align}
We present a stronger result involving non-asymptotic bounds below.

Let $\Phi^{(n,k)}_\theta$ denote the distribution of $\bm{\tau}^{(n)}_{1:k}$, and by a slight abuse of notation, we will also denote the \ac{CDF} of  $\Phi^{(n,k)}_\theta$ and $\Phi_\theta$ by the same respective symbols, that is, $\Phi_\theta(s)  = \Phi_\theta([0,s])$ and 
\begin{align*}
\Phi^{(n,k)}_\theta(s_1,s_2,\hdots, s_k) & \equiv \Phi^{(n,k)}\lf(\prod_{i=1}^k[0,s_i]\ri) = \PP_\theta\lf(\tau^{(n)}_i \leq s_i,\  i=1,\hdots,k\ri). 
\end{align*}

\begin{myCorollary}
    \label{cor:prop-chaos-tau}
Let $T>0$. Then,  
\begin{align*}
\sup_{0\leq s_1,\hdots, s_k\leq T}\Big|\Phi^{(n,k)}_\theta(s_1,& s_2,\hdots, s_k)  - \prod_{i=1}^k \Phi_\theta(s_i) \Big| 
\leq\ \cnst_2^{1/2} \lf(\f{\theta_1 T}{1+\theta_2} \wedge 1\ri)^{1/2}\f{\theta_1 T}{\theta_2^2}\exp\lf(2\theta_1 T/\theta_2\ri)k n^{-1/2}.
\end{align*}
\end{myCorollary}

\begin{proof}[Proof of \Cref{cor:prop-chaos-tau}]
 For the proof, we construct a limiting \ac{IPS} that is coupled to the prelimit \ac{IPS} \eqref{eq:particle_i_trajectory}  through the reuse of the \acp{PRM} $\eta_i$. To be more precise, consider the collection of \ac{iid} stochastic processes $\subst_i, i=1,2,\ldots,$  defined by
    \begin{align}
        \begin{aligned}
            \subst_i(t) &{}= 1 - \int_{[0, t]\times [0, \infty)} \indicator{[0, Y_{i}(s-) g_{\theta}(Z_S(s)) ]}{v} \eta_i\left(\differential{v}\times \differential{s}\right)\eqcomma 
        \end{aligned}
        \label{eq:new_IPS}
    \end{align}
    for $t\ge 0$
    and let
    $\tau_i \defeq \inf\{t\ge 0 : \subst_i(t) = 0\}.$ Clearly,
\begin{align}
    \label{eq:tau-Y-equiv}
\subst_i^{(n)}(t) = \indic{\{\tau^{(n)}_i>t\}}, \quad  \subst_i(t) =  \indic{\{\tau_i>t\}}.
\end{align}  
Note that the representation $$\Phi^{(n,k)}_\theta(s_1,s_2,\hdots, s_k) = \EE \lf(\prod_{i=1}^k \indic{\{\tau^{(n)}_i \leq s_i\}}\ri) = \EE \lf(\prod_{i=1}^k \lf(1-\subst_i^{(n)}(s_i)\ri)\ri)$$ and the inequality $\left| \prod_{i=1}^k w_i - \prod_{i=1}^k z_i \right| \leq \sum_{i=1}^k |w_i - z_i|$ for  $|w_i|, |z_i| \leq 1$ give
\begin{align}
    \label{eq:cdf-diff}
\lf|\Phi^{(n,k)}_\theta(s_1,s_2,\hdots, s_k) - \prod_{i=1}^k \Phi_\theta(s_i) \ri| 
\leq& \sum_{i=1}^k \EE|\subst_i^{(n)}(s_i) - \subst_i(s_i)|,
\end{align}
Next from \eqref{eq:particle_i_trajectory} and \eqref{eq:new_IPS}, we have
    \begin{align*}
      \EE\lf[\sup_{t\leq T} |\subst^{(n)}_i(t) - \subst_i(t)|\ri]\leq &\   \int_0^T \EE\lf|\subst_i^{(n)}(s) g_{\theta}(\bar{\subst}^{(n)}(s)) -\subst_i(s) g_{\theta}(Z_S(s))\ri|\ \differential{s}\\
      \leq &\ A_n(T)+\int_0^T \sup_{s\leq r}g_\theta(Z_S(s))\EE\lf[\sup_{s\leq r} |\subst^{(n)}_i(s) - \subst_i(s)|\ri] \differential{r},
    \end{align*}
    where
    \begin{align*}
        A_n(T)\defeq  \int_0^T \EE \lf[\subst_i^{(n)}(s)\lf| g_{\theta}(\bar{\subst}^{(n)}(s)) - g_{\theta}(Z_S(s))\ri|\ri]\differential{s}\eqstop
    \end{align*}
Noting that $0\leq \sup_{y \geq 0}g_\theta(y) \leq \theta_1/\theta_2$ and $\|g_\theta\|_{\mrm{Lip}} = \theta_1/\theta_2^2$, we get by Gr\"onwall's inequality,
    \begin{align*}
        \EE\lf[\sup_{t\leq T} |\subst^{(n)}_i(t) - \subst_i(t)|\ri] \leq&\ \sup_{t\leq T}\EE[|\bar{\subst}^{(n)}(t) - Z_S(t)|]\f{\theta_1 T}{\theta_2^2} e^{\theta_1T/\theta_2}\eqstop 
    \end{align*}
The assertion now follows from \Cref{prop:FLLN_Tbar} and \eqref{eq:cdf-diff}.

\end{proof}

\begin{myRemark}
    \label{rem:prop-chaos-tau}
\Cref{cor:prop-chaos-tau} implies that for any $k_n =o(n^{1/2})$ and $T>0$, $\sup_{0\leq s_1,\hdots, s_{k_n}\leq T}\Big|\Phi^{(n,k_n)}_\theta(s_1,s_2,\hdots, s_{k_n})  - \prod_{i=1}^{k_n} \Phi_\theta(s_i) \Big|  \nrt 0.$ In particular, for any fixed $k$, the convergence in \eqref{eq:tau-prop-chaos} holds.
\end{myRemark}

\subsection{Approximate likelihood and $\theta$-estimator}
\label{sec:like}
Our data vector, $\bm{t}_{1:k_n} \defeq (t_1,t_2,\hdots, t_{k_n})$, collected over a fixed finite time interval $[0,T]$ for some $T>0$, is a realization of the random variable 
\begin{align}
    \label{eq:sample-rv-sym}
 \bm{\tilde \tau}^{(n)}_{1:K_n} \defeq (\tilde \tau^{(n)}_1, \tilde \tau^{(n)}_2, \hdots, \tilde \tau^{(n)}_{K_n})  
\end{align} 
consisting of a random sample of size $K_n$ of product formation times $\tau^{(n)}_i$ that lie in  $[0,T]$. Here 
\begin{align}
    \label{eq:def-Kn}
 K_n \equiv N_n(\tau^{(n)}_{1:n}) \wedge \tilde k_n,   
\end{align}
 where $N_n(\tau^{(n)}_{1:n}) \equiv \sum_{i=1}^n \indic{\{\tau^{(n)}_i \leq T\}}$ is the number of product formation times that fall in $[0,T]$ and $\tilde k_n$ is a fixed number. The motivation behind this definition is that the number of products produced in any given interval $[0, T]$ is a random variable and a priori unknown. In practice, we can only observe a small number of the  product formation times. Here, the numbers $\tilde k_n$ are our choice. We will discuss the precise conditions on $\tilde k_n$ required for the consistency of our proposed estimator later in this section.

For a fixed $k$, we first derive an expression for the distribution of $\bm{\tilde \tau}^{(n)}_{1:k}$. As before, we suppress the notation $\theta$ from the probability measure $\PP$. Given $\bm{r}_{1:n} = (r_1,r_2,\hdots, r_n) \in [0,\infty)^n$, let 
$$\SC{I}^{(n)}(\bm{r}_{1:n}) \defeq \{ i\in \{1,2,\hdots,n\}: r_i \in [0,T]\}$$ 
denote the set of indices for which the given times $r_i$ lie in the observation window $[0,T]$, and let 
$$N_n(\bm{r}_{1:n}) = |\SC{I}^{(n)}(\bm{r}_{1:n})| =\sum_{i=1}^n \indic{\{r_i\leq T\}}$$ 
be the size of the set $\SC{I}^{(n)}(\bm{r}_{1:n})$.
 For integers $k\leq M$, let
$$\text{Perm}(M,k) =\lf\{\bm{j}_{1:k} =(j_1,\hdots, j_k): j_i \in \{1,2,\hdots,M\}, j_i \neq j_l, \text{ for } i\neq l \ri\}$$ 
denote the set of $k$-permutations from the set $\{1,2,\hdots,M\}$
and 
$$P(M,k) \defeq |\text{Perm}(M,k)| = \frac{M!}{(M-k)!}.$$

Let $\bm{J}^{(n)}_{1:k} = (J^{(n)}_1,J^{(n)}_2,\hdots, J^{(n)}_{k})$ be a random vector  taking values in $\text{Perm}(n,k)$ with distribution given by
\begin{align*}
    \PP\lf(J^{(n)}_i = j_i,\ i=1,2,\hdots,k \mid \bm{\tau}^{(n)}_{1:n}\ri)(\om) = \f{\prod_{i=1}^{k} \indic{\{\tau^{(n)}_{j_i}(\om) \leq T\}}}{ P(N_n(\bm{\tau}^{(n)}_{1:n}(\om)), k)}.
\end{align*}
This simply means that conditional on $\bm{\tau}^{(n)}_{1:n}$,   $\bm{J}^{(n)}_{1:k}$ is a random vector uniformly distributed over $k$ permutations of $\SC{I}^{(n)}(\bm{\tau}^{(n)}_{1:n})$. Define $\bm{\tilde \tau}^{(n)}_{1:k}$ as
\begin{align} \label{eq:sample-rv-def}
   \tilde \tau^{(n)}_i(\om) = \tau^{(n)}_{J^{(n)}_i(\om)}(\om). 
\end{align}

\np
{\em Distribution of }  $\bm{\tilde \tau}^{(n)}_{1:k}$: 
Let $\rho^{(n)}_{1:k}$ denote the distribution of $\bm{\tilde \tau}^{(n)}_{1:k}$.  Note that for $\bm{s}_{1:k} =(s_1,\hdots, s_{k}) \in [0,T]^{k}$,
\begin{align*}
 \rho^{(n)}_{1:k}\lf(\prod_{i=1}^{k}[0,s_i]\ri)& \equiv   \PP\lf( \bm{\tilde \tau}^{(n)}_{1:k }\leq \ \bm{s}_{1:k_n}\ri) =  \EE\lf[\sum_{\bm{j}_{1:k} \in \text{Perm}(n,k)} \f{\prod_{i=1}^{k} \indic{\{\tau^{(n)}_{j_i} \leq s_i\}}}{ P(N_n(\bm{\tau}^{(n)}_{1:n}), k)} \ri]\\
     =& \ P(n,k_n) \EE\lf[\f{\prod_{i=1}^{k} \indic{\{\tau^{(n)}_{i} \leq s_i\}}}{ P(N_n(\bm{\tau}^{(n)}_{1:n}), k)}\ri]\\
     =& \ P(n,k) \int_{\left(\prod_{i=1}^{k}[0,s_i]\right)\times[0,\infty)^{n-k}} \frac{\tauden^{(n)}_{n}(\bm{r}_{1:k}, \bm{u}_{k+1:n})}{P(N_n(\bm{r}_{1:k}, \bm{u}_{k+1:n}), k_n)} \differential{\bm{r}_{1:k}}\differential{\bm{u}_{k+1:n}}\eqcomma 
\end{align*}
where the middle equality is by the exchangeability of $\{\tau^{(n)}_i\}$, and $\tauden^{(n)}_{k}$ denotes the density of $\bm{\tau}^{(n)}_{1:k}.$

It follows that the density $q^{(n)}_{1:k}$ of $\bm{\tilde\tau}^{(n)}_{1:k}$ is given by
\begin{align*}
 q^{(n)}_{1:k}(\bm{s}_{1:k}) =&\ P(n,k) \int_{[0,\infty]^{n-k}} \frac{\tauden^{(n)}_{n}(\bm{s}_{1:k}, \bm{u}_{k+1:n})}{P(N_n(\bm{s}_{1:k}, \bm{u}_{k+1:n}), k)} \differential{\bm{u}_{k+1:n}}\\
 =&\ P(n,k) \tauden^{(n)}_{k}(\bm{s}_{1:k}) \int_{[0,\infty)^{n-k}} \frac{1}{P(N_n(\bm{s}_{1:k}, \bm{u}_{k+1:n}), k_n)} \f{\tauden^{(n)}_{n}(\bm{s}_{1:k}, \bm{u}_{k_n+1:n})}{\tauden^{(n)}_{k}(\bm{s}_{1:k}) }\differential{\bm{u}_{k+1:n}}\\
 =&\  P(n,k) \tauden^{(n)}_{k}(\bm{s}_{1:k}) \EE\lf[\f{1}{P(N_n(\bm{\tau}^{(n)}_{1:n}, k)} \ \Big|\  \bm{\tau}^{(n)}_{1:k} = \bm{s}_{1:k}\ri] \eqstop 
\end{align*}
Notice that the random variables   $\tilde \tau^{(n)}_i$ are exchangeable. Next, \eqref{eq:Yn-mean-conv} and \Cref{cor:prop-chaos-tau} suggest that for large $n$, under $\PP_\theta$
\begin{align*}
  \tauden^{(n)}_{k}(\bm{s}_{1:k}) \approx \prod_{i=1}^{k} \tauden_\theta(s_i), \quad \f{N_n(\bm{\tau}^{(n)}_{1:n})}{n}  \approx \taudist_\theta([0,T]),
\end{align*}
where $\tauden_\theta$ is the density of the probability measure $\taudist_\theta$; thus
$$\tauden_\theta(t) = \f{\differential{\taudist_{\theta}([0,t])}}{\differential{t}} = -Z'_{S,\theta}(t) = h_\theta(Z_{S,\theta}(t)).$$
Since $\f{P(n,k)}{P(N_n, k)} \approx (n/N_n)^{k}$, we see that  
\begin{align*}
 q^{(n)}_{1:k}(\bm{s}_{1:k}) \approx \f{\prod_{i=1}^{k} \tauden_\theta(s_i)}{\taudist^{k}_\theta([0,T])} = \f{\prod_{i=1}^{k} h_\theta(Z_{S,\theta}(s_i))}{(1-Z_{S,\theta}(T))^{k}}.
 \end{align*}
The approximations above are formal, but they motivate estimating $\theta$ via the approximate likelihood function $\like(\theta \mid \bm{\tilde \tau}^{(n)}_{1:K_n})$, defined as
\begin{align}\label{eq:approx-like}
\like(\theta\mid \bm{\tilde \tau}^{(n)}_{1:K_n}) 
\defeq \f{\prod_{i=1}^{K_n} h_\theta(Z_{S,\theta}(\tilde{\tau}^{(n)}_i))}{(1-Z_{S,\theta}(T))^{K_n}}.
\end{align}
The main objective in the remainder of the paper is to show that, given data $\bm{t}_{1:k_n}$ arising as a realization of $\bm{\tilde \tau}^{(n)}_{1:K_n}$, the estimator $\hat\theta_n$ defined by
\begin{align}
\label{eq:est-def}
\hat\theta_n \equiv \hat\theta_n(\bm{\tilde \tau}^{(n)}_{1:K_n}) 
= \argmax_{\theta \in \Theta} \ln \like(\theta\mid \bm{\tilde \tau}^{(n)}_{1:K_n})
\end{align}
is consistent. 

\subsection{Consistency} \label{sec:const}


\begin{myDefinition}\label{def:consistency}
Consider the \ac{IPS} defined by \eqref{eq:particle_i_trajectory}. An adapted sequence of random variables $\{\hat \theta^{(n)}: n \ge 1 \}$ is  a consistent estimator for $\theta^* \in \Theta$ if
$
\hat \theta^{(n)} \stackrel{\PP_{\theta^*}}\Rt \theta^*
$ as $n \rt \infty$.
\end{myDefinition}

\np
We now present the main result of this section.

\begin{myTheorem}
\label{th:const}
Let $\Theta \subset (0,\infty)^2$ be compact, and fix $\theta^* \in \Theta$. Suppose $\tilde k_n \leq \alpha\ \taudist_{\theta^*}([0,T]) n$ for some $0<\alpha <1$, and $\tilde k_n \to \infty$, where $K_n, \tilde k_n$ are as in \eqref{eq:def-Kn}. Then, the sequence of estimators $\{\hat \theta^{(n)} \equiv \hat\theta_n(\bm{\tilde \tau}^{(n)}_{1:K_n}) : n\ge 1\}$ defined in \eqref{eq:est-def}, with $\bm{\tilde \tau}^{(n)}_{1:K_n}$  given by \eqref{eq:sample-rv-def}, is consistent for any $\theta^* \in \Theta$.
\end{myTheorem}
Since the true $\theta^*$ is unknown, we can take $\tilde k_n = o(n)$  to ensure that the condition on $\tilde k_n$ is satisfied. The proof of the above result hinges on the following important theorem establishing the propagation of chaos for the exchangeable sequence $\bm{\tilde \tau}^{(n)}_{1:k_n}$, which is interesting in its own right.

For a fixed $T>0$, define the probability measure $\tilde{\taudist}_{\theta,T} \in \mathcal{P}([0,T])$ by 
\begin{align}
    \label{eq:cond-dist}
    \tilde{\taudist}_{\theta,T}([0,t]) = \taudist_\theta([0,t])/\taudist_\theta([0,T]) = \frac{1-Z_{S,\theta}(t)}{1-Z_{S,\theta}(T)}, \quad t \in [0,T]\eqstop 
\end{align}
The density $\tilde{\tauden}_{\theta,T}$ of $\tilde{\taudist}_{\theta,T}$ is given by
\begin{align}
    \label{eq:cond-den}
\tilde{\tauden}_{\theta,T}(t) \defeq \f{\differential{\tilde{\taudist}_{\theta,T}([0,t])}}{\differential{t}} = \f{h_\theta(Z_{S,\theta}(t))}{1-Z_{S,\theta}(T)}. 
\end{align}

\begin{myTheorem}
    \label{th:prop-chaos}
    Let $m \in \N$ be fixed and $\psi: \N_0 \rt [0,\infty)$ be a bounded function. Then, for any $\theta \in \Theta$ and $T>0$, there exists a constant $\bar C_m(T,\theta)$ such that 
\begin{align*}
  \sup_{f_1,\hdots, f_m \in \mathscr{F}_{T,1}}\lf| \EE_{\theta}\lf[\lf(\prod_{i=1}^m \bar{f}_i(\tilde \tau^{(n)}_i)\ri)\psi(N_n(\bm{\tau}^{(n)}_{1:n}))\ri]\ri| \leq \|\psi\|_\infty \bar C_m(T,\theta)n^{-m/2}\eqcomma 
\end{align*}
where $\bar{f}_{i,\theta} = f_i - \tilde \taudist_{\theta,T}(f_i) = f_i -(\int_0^Tf_i\differential{}\taudist_\theta)/\taudist_\theta([0,T])$.
\end{myTheorem}
\begin{myCorollary}
     Let $m \in \N$ be fixed. Then, for any $\theta \in \Theta$ and $T>0$, there exists a constant $C_m(T,\theta)$ such that
    $$\sup_{f_1,\hdots, f_m \in \mathscr{F}_{T,1}}\lf|\EE_{\theta} \lf(\prod_{i=1}^m f_i(\tilde \tau^{(n)}_i)\ri) - \prod_{i=1}^m \tilde{\taudist}_{\theta,T}(f_i)\ri| \leq C_m(T,\theta)n^{-1/2}.$$
\end{myCorollary}

\begin{proof}[Proof of \Cref{th:prop-chaos}]
For convenience, we establish the result for $m=2$; the proof for general finite $m$ is essentially the same, except for slightly more cumbersome notation. As before, we will drop $\theta$  from the notation.

Notice that 
$$ \EE\lf(\left(\prod_{i=1}^2 \bar{f}_i(\tilde \tau^{(n)}_i) \right)\psi(N_n(\bm{\tau}^{(n)}_{1:n}))\ri) = \EE\lf(\EE\lf(I^{(n)}_T(f_1,f_2) \mid \bm{\tau}^{(n)}_{1:n}\ri)\psi(N_n(\bm{\tau}^{(n)}_{1:n}))\ri) = \EE\lf(I^{(n)}_T(f_1,f_2)\psi(N_n(\bm{\tau}^{(n)}_{1:n}))\ri),$$ 
where 
\begin{align*}
   I^{(n)}_T(f_1,f_2) & \defeq \f{1}{N_n(\bm{\tau}^{(n)}_{1:n})(N_n(\bm{\tau}^{(n)}_{1:n})-1)}\sum_{\substack{i,i'=1 \\ i\neq i'}}^n \bar{f}_1(\tau^{(n)}_{i'})\bar{f}_2(\tau^{(n)}_i)\indic{\{\tau^{(n)}_{i} \leq T\}}\indic{\{\tau^{(n)}_{i'} \leq T\}}\\
   & = \f{1}{N_n(\bm{\tau}^{(n)}_{1:n})(N_n(\bm{\tau}^{(n)}_{1:n})-1)} \lf(n^2 \Xi^{(n)}(\bar{f}_1\indic{[0,T]})\Xi^{(n)}(\bar{f}_2\indic{[0,T]})- n \Xi^{(n)}(\bar{f}_1\bar{f}_2\indic{[0,T]})\ri)\\
   & = \f{1}{N_n(\bm{\tau}^{(n)}_{1:n})/n(N_n(\bm{\tau}^{(n)}_{1:n})/n-1/n)} \tilde I^{(n)}_T(f_1,f_2).
\end{align*}
The middle equality above used the identity $\sum_{i\neq i'}a_ib_{i'} = (\sum_i a_i)(\sum_i b_i)-\sum_i a_ib_i$, and $\tilde I^{(n)}_T(f_1,f_2)$ is given by 
$$\tilde  I^{(n)}_T(f_1,f_2) \defeq \lf(\Xi^{(n)}(\bar{f}_1\indic{[0,T]})\Xi^{(n)}(\bar{f}_2\indic{[0,T]})- n^{-1} \Xi^{(n)}(\bar{f}_1\bar{f}_2\indic{[0,T]})\ri).$$
Now, \eqref{eq:emp-measure-f-ineq} shows that
\begin{align*}
|\Xi^{(n)}(\bar f_1 \indic{[0,T]})| =&\ \lf|\Xi^{(n)}(f_1\indic{[0,T]})-\f{\taudist(f_1\indic{[0,T]})}{\taudist([0,T])} \Xi^{(n)}([0,T])\ri|\\
=&\ \lf|\lf(\Xi^{(n)}(f_1\indic{[0,T]}) - \taudist(f_1\indic{[0,T]})\ri)+ \f{\taudist(f_1\indic{[0,T]})}{\taudist([0,T])} \lf(\taudist([0,T])-\Xi^{(n)}([0,T])\ri)\ri|\\
\leq&\ 2\lf((T\vee 1)\|f_1\|_{\text{BL},T} \sup_{t\leq T}  \absolute{\Xi^{(n)}([0,t]) - \taudist([0,t])} \wedge \|f_1\|_{\infty,T}\ri). 
\end{align*}
Hence, for $f_1, f_2 \in \mathscr{F}_{T,1}$, 
\begin{align*}
 |\tilde  I^{(n)}_T(f_1,f_2)| \leq 4 \lf((T\vee 1)^2\sup_{t\leq T}  \absolute{\Xi^{(n)}([0,t]) - \taudist([0,t])}^2\wedge 1+  n^{-1}\ri)\eqstop 
\end{align*}
Now, choose a $\delta$ such that $2\delta<\taudist([0,T])$. Then, for every $n>\delta^{-1}$, 
\begin{align}\label{eq:I-part1}
\begin{aligned}
| I^{(n)}_T(f_1,f_2)|&\indic{\{|N_n(\bm{\tau}^{(n)}_{1:n})/n-\taudist([0,T]| \leq  \delta\}} \leq \f{1}{(\taudist([0,T])-\delta)(\taudist([0,T])-\delta-n^{-1})}  |\tilde I^{(n)}_T(f_1,f_2)|\\
& \leq \f{4}{(\taudist([0,T])-\delta)(\taudist([0,T])-2\delta)} \Big((T\vee 1)^2\sup_{t\leq T}  \absolute{\Xi^{(n)}([0,t]) - \taudist([0,t])}^2+n^{-1}\Big)\eqstop 
\end{aligned}
\end{align}
Next, since $\om \in \{\tau^{(n)}_{i} \leq T\} \cap \{\tau^{(n)}_{i'} \leq T\}$ implies that $N_n(\bm{\tau}^{(n)}_{1:n}(\om)) \geq 2$, we have 
\begin{align}\label{eq:I-part2}
\begin{aligned}
| I^{(n)}_T(f_1,f_2)|\indic{\{|N_n(\bm{\tau}^{(n)}_{1:n})/n-\taudist([0,T]|> \delta\}} & \leq \f{n^2}{2} |\tilde  I^{(n)}_T(f_1,f_2)|\indic{\{|N_n(\bm{\tau}^{(n)}_{1:n})/n-\taudist([0,T]|> \delta\}}\\
&\leq 4n^2\indic{\{|N_n(\bm{\tau}^{(n)}_{1:n})/n-\taudist([0,T]|> \delta\}}\eqstop 
\end{aligned}
\end{align}
Since $N_n(\bm{\tau}^{(n)}_{1:n})/n  = \Xi^{(n)}([0,T])$, it follows by \eqref{eq:I-part1}, \eqref{eq:I-part2},  the first inequality in \eqref{eq:emp-meas-conc} along with Markov inequality that
\begin{align*}
 \Big|\EE  \Big(\Big(\prod_{i=1}^2 & \bar{f}_i(\tilde \tau^{(n)}_i)\Big) \psi(N_n(\bm{\tau}^{(n)}_{1:n}))\Big)\Big|   = \lf|\EE\lf(I^{(n)}_T(f_1,f_2)\psi(N_n(\bm{\tau}^{(n)}_{1:n}))\ri)\ri| \leq \|\psi\|_\infty\EE|I^{(n)}_T(f_1,f_2)| \\
 \leq&  4\|\psi\|_\infty\lf(\f{(T\vee 1)^2\cnst_2 \lf(\f{\theta_1 T}{1+\theta_2} \wedge 1\ri)\exp\lf(2(\theta_1/\theta_2) T\ri)n^{-1}+n^{-1}}{(\taudist([0,T])-\delta))(\taudist([0,T])-2\delta)} +  \f{n^2 \cnst_6 \lf(\f{\theta_1 T}{1+\theta_2} \wedge 1\ri)\exp\lf(6(\theta_1/\theta_2) T\ri)n^{-3} }{\delta^{6}}\ri) \\
\equiv &\ \|\psi\|_\infty \bar C_2(T,\theta)n^{-1}, 
\end{align*}
with 
$$\bar C_2(T,\theta)   \defeq  4 \lf(\f{(T\vee 1)^2\cnst_2 \lf(\f{\theta_1 T}{1+\theta_2} \wedge 1\ri) \exp\lf(2(\theta_1/\theta_2) T\ri)+1}{(\taudist([0,T])-\delta))(\taudist([0,T])-2\delta)} +  \f{ \cnst_6 \lf(\f{\theta_1 T}{1+\theta_2} \wedge 1\ri) \exp\lf(6(\theta_1/\theta_2) T\ri)}{\delta^{6}}\ri).$$
\end{proof}  
Define the empirical measure of $\bm{\tilde \tau}^{(n)}_{1:K_n} \defeq (\tilde \tau^{(n)}_1, \tilde \tau^{(n)}_2, \hdots, \tilde \tau^{(n)}_{K_n})$ by
 \begin{align}
     \label{eq:emp-meas-data}
   \tilde \Xi^{(n)} \defeq \f{1}{K_n} \sum_{i=1}^{K_n} \delta_{\tilde\tau^{(n)}_i} \eqstop 
 \end{align}
Thus, $\tilde \Xi^{(n)}(f) = \f{1}{K_n}\sum_{i=1}^{K_n} f(\tilde \tau^{(n)}_i)$. \Cref{th:prop-chaos} immediately leads to the following result on the convergence of $\tilde \Xi^{(n)}$.

\begin{myCorollary}\label{cor:emp-meas-conv}
Fix  $\theta \in \Theta$. Suppose $\tilde k_n \leq \alpha\ \taudist_\theta([0,T]) n$, for some $0<\alpha <1$. Then for $f\in\mathrm{Lip}([0,T],\mathbb R)$,  
\begin{align}\label{eq:emp-meas-L2-bd}
    \EE_{\theta}\lf(\lf|\tilde \Xi^{(n)}(f) - \tilde\taudist_{\theta,T}(f)\ri|^2\ri) \leq \bar C_3(T, \theta)\| f\|^2_{\text{BL},T} \Big/ \tilde k_n,
\end{align}
where
$$\bar C_3(T, \theta) = 4 +\bar C_2(T,\theta) + \f{4 \cnst_2 \lf(\f{\theta_1 T}{1+\theta_2} \wedge 1\ri)\exp\lf( 2\theta_1 T/\theta_2\ri)}{(1-\alpha)^2\taudist^2_\theta([0,T])}$$ 
with $\bar C_2(T,\theta)$ as in \Cref{th:prop-chaos}.
In particular, when $\tilde k_n \to \infty$, $\tilde \Xi^{(n)}  \stackrel{\PP_{\theta}}\Rt \tilde \taudist_{\theta, T}$, i.e., as $n \rt \infty$,
\begin{align}\label{eq:emp-meas-conv-prob}
d_{\mrm{BL}}(\tilde \Xi^{(n)},\ \tilde\taudist_{\theta,T}) \equiv \sup_{f \in \mathscr{F}_{T,1}} |\tilde \Xi^{(n)} (f) - \tilde \taudist_{\theta, T}(f)| \stackrel{\PP_{\theta}}\Rt 0.
\end{align}
\end{myCorollary}

\begin{proof}[Proof of \Cref{cor:emp-meas-conv}]
   Letting $\bar f = f - \tilde \taudist_T(f),$
 \begin{align*}
 \EE\lf(\lf| \tilde \Xi^{(n)}(\bar f)\ri|^2\ri) =  \EE\lf(\lf| \tilde \Xi^{(n)}(\bar f)\ri|^2 \indic{\{\tilde k_n \leq N_n(\bm{\tau}^{(n)}_{1:n})\}}\ri) + \EE\lf(\lf| \tilde \Xi^{(n)}(\bar f)\ri|^2 \indic{\{\tilde k_n > N_n(\bm{\tau}^{(n)}_{1:n})\}}\ri).
 \end{align*}
 Now by the exchangeability of $\tilde \tau^{(n)}_1, \hdots, \tilde \tau^{(n)}_{\tilde k_n}$ and \Cref{th:prop-chaos}
 \begin{align*}
 \EE\lf(\lf| \tilde \Xi^{(n)}(\bar f)\ri|^2 \indic{\{\tilde k_n \leq N_n\}}\ri)=
 & \f{1}{\tilde k_n^2}\Big( \sum_{i=1}^{\tilde k_n} \EE\lf(\bar f(\tilde \tau^{(n)}_i)^2 \indic{\{\tilde k_n \leq N_n\}}\ri) +  \sum_{i\neq i'}^{\tilde k_n} \EE\lf(\bar f(\tilde \tau^{(n)}_i)\bar f(\tilde \tau^{(n)}_{i'})\indic{\{\tilde k_n \leq N_n\}}\ri) \Big) \\
 \leq &\ \f{1}{\tilde k_n} \EE\lf(\bar f(\tilde \tau^{(n)}_1)^2\ri) + \f{\tilde k_n(\tilde k_n-1)}{\tilde k^2_n} \f{\bar C_2(T,\theta) \| f\|^2_{\text{BL},T}}{n} \\
 \leq &\  \f{ (4 +\bar C_2(T,\theta))\| f\|^2_{\text{BL},T}}{ \tilde k_n} \eqcomma,
 \end{align*}
as $\tilde k_n \leq \alpha\ \taudist_\theta([0,T]) n < n.$ 
Next since $0<\alpha<1$ and  $N_n(\bm{\tau}^{(n)}_{1:n})/n \equiv \Xi^{(n)}([0,T])$, we have
 \begin{align*}
     \EE\lf(\lf| \tilde \Xi^{(n)}(\bar f)\ri|^2 \indic{\{\tilde k_n > N_n(\bm{\tau}^{(n)}_{1:n})\}}\ri) \leq &\  4 \|f\|^2_{\infty, T} \PP\lf(N_n(\bm{\tau}^{(n)}_{1:n}) < \tilde k_n\ri)\\
     \leq&\ 4 \|f\|^2_{\infty, T} \PP\lf(N_n(\bm{\tau}^{(n)}_{1:n})/n - \taudist([0,T])  < (\alpha-1)\taudist([0,T])\ri)\\
     \leq&\ 4 \|f\|^2_{\infty, T}\PP\lf(\absolute{\Xi^{(n)}([0,T]) - \taudist([0,T])} > (1-\alpha)\taudist([0,T])\ri)\\
     \leq &\ 4 \|f\|^2_{\infty, T}\cnst_2 \lf(\f{\theta_1 T}{1+\theta_2} \wedge 1\ri)\exp\lf( 2\theta_1 T/\theta_2\ri) n^{-1} \Big/ (1-\alpha)^2\taudist^2([0,T]),
 \end{align*}
 where the last step used \eqref{eq:emp-meas-conc} and Markov inequality. This proves \eqref{eq:emp-meas-L2-bd}. The convergence in \eqref{eq:emp-meas-conv-prob} now follows by \cite[Lemma 4.8]{Kallenberg2017RandomMeasures}. 
\end{proof}

\begin{myRemark}\label{rem:emp-meas-conv}
 Notice that by the \ac{DCT}, $d_{\mrm{BL}}(\tilde \Xi^{(n)}, \tilde\taudist_{\theta,T})$ also converges to $0$ in $L^1(\PP_{\theta})$, that is,
$$\EE_{\theta}\lf(d_{\mrm{BL}}(\tilde \Xi^{(n)}, \tilde\taudist_{\theta,T})\ri) \equiv \EE_{\theta}\lf[\sup_{f \in \mathscr{F}_{T,1}} |\tilde \Xi^{(n)} (f) - \tilde \taudist_T(f)|\ri]  \nrt 0.$$
\end{myRemark}

From \eqref{eq:approx-like},  the normalized log-likelihood function $\ell_n(\theta) \defeq K_n^{-1}\ln \like(\cdot|\tilde{\bm{\tau}}^{(n)}_{1:K_n})$ is given by 
\begin{align} \label{eq:log-like}
\ell_n(\theta) \defeq \f{1}{K_n} \ln \like(\cdot|\tilde{\bm{\tau}}^{(n)}_{1:K_n}) =  \f{1}{K_n} \sum_{i=1}^{K_n} \ln\lf(\tilde \tauden_{\theta,T}(\tilde \tau^{(n)}_i)\ri) \equiv \f{1}{K_n} \sum_{i=1}^{K_n} \ln\lf(\f{h_\theta(Z_{S,\theta}(\tilde \tau^{(n)}_i))}{1-Z_{S,\theta}(T)}\ri),  
\end{align}
%
and that the parameter of interest $\theta$ is given by \eqref{eq:para-red}. Define
\begin{align}
    \theta_{i,\mrm{max}} = \max\{\theta_i: \theta = (\theta_1, \theta_2) \in \Theta\}, \quad \theta_{i,\mrm{min}} = \min\{\theta_i: \theta = (\theta_1, \theta_2) \in \Theta\},
\end{align}
which are finite when $\Theta $ is compact. 

\begin{myTheorem}\label{th:conv-log-like}
    Let $\Theta \subset (0,\infty)^2$ be compact, and fix $\theta^* \in \Theta$. Then under the assumption of \Cref{th:const},
\begin{align*}
    \EE_{\theta^*}\lf(\sup_{\theta \in \Theta}\lf|\ell_n(\theta) + \mrm{CE}(\tilde\taudist_{\theta^*,T}, \tilde\taudist_{\theta,T})\ri|\ri) \nrt 0\eqcomma 
\end{align*}
where $\mrm{CE}(\tilde\taudist_{\theta^*,T}, \tilde\taudist_{\theta,T})$ is the cross-entropy of the probability measure $\tilde\taudist_{\theta^*,T}$ relative to $\tilde\taudist_{\theta,T}$, and is given by
\begin{align*}
  \mrm{CE}(\tilde\taudist_{\theta^*,T}, \tilde\taudist_{\theta,T}) &\defeq -\int_0^T \ln (\tilde \tauden_{\theta,T}(t)) \tilde \tauden_{\theta^*,T}(t) \differential{t} 
  = \mrm{Ent}(\tilde\taudist_{\theta^*,T}) + \mrm{RE}(\tilde \taudist_{\theta^*,T}\|\tilde \taudist_{\theta,T}),
\end{align*}
with  $\tilde{\taudist}_{\theta,T}$ and its density $\tilde{\tauden}_{\theta,T}$ given by \eqref{eq:cond-dist} and \eqref{eq:cond-den}, respectively.

\end{myTheorem}

\np
Recall that the relative entropy $\mrm{RE}(\cdot\|\cdot)$, entropy $\mrm{Ent}(\cdot)$, and the cross entropy $\mrm{CE}(\cdot, \cdot)$ are defined in \Cref{sec:out-not}.

\begin{proof}[Proof of \Cref{th:conv-log-like}]
Notice that the hypothesis on $\Theta$ implies $0< \theta_{i,\mrm{min}} \leq  \theta_{i,\mrm{max}} < \infty$, $i=1,2$.
As noted before, from the expression of $h_\theta$ (see \eqref{eq:h_defn}), it is clear that $h_\theta$ is increasing, and 
 \begin{align} \label{eq:h-bdd-lip}
 \begin{aligned}
 0\leq    h_\theta(y) & \leq \min\lf\{(\theta_1/\theta_2)  y,\   \theta_1/(\theta_2+1)\ri\}, \quad 0\leq y \leq 1,\quad  
 \text{and } \|h_\theta\|_{\mrm{Lip}} 
 \equiv \sup_{y\geq 0}|h'_\theta(y)| =  \theta_1/\theta_2.
 \end{aligned}
\end{align}
 Therefore, $h_\theta$ is a bounded Lipschitz function with uniform in $\theta$-bound and  Lipschitz constant  
 \begin{align}
     \label{eq:h-bdd-lip-unif}
 \begin{aligned}    
\|h\|_{\infty,1,\mrm{max}} \defeq&\  \sup_{\theta \in \Theta} \sup_{y \in [0,1]} h_\theta(y) \leq \theta_{1,\mrm{max}}/(\theta_{2,\mrm{min}}+1), \\
\|h\|_{\mrm{Lip,max}} \defeq&\  \sup_{\theta \in \Theta} \|h_\theta\|_{\mrm{Lip}}  \leq \theta_{1,\mrm{max}}/\theta_{2,\mrm{min}} \eqstop 
\end{aligned}
 \end{align}
The first inequality immediately implies that $Z_{S,\theta}$ is Lipschitz continuous with Lipschitz constant uniformly bounded in $\theta$:
\begin{align}\label{eq:Z-lip-unif}
    \|Z\|_{\mrm{Lip, max}} \defeq \sup_{\theta \in \Theta} \|Z_{S,\theta}\|_{\mrm{Lip}} \equiv  \sup_{\theta \in \Theta} \sup_{t \in [0,\infty)} \|Z'_{S,\theta}(t)\| \leq \|h\|_{\infty,1,\mrm{max}}. 
\end{align}
Next, it follows from \eqref{eq:Z-lbd} that 
$$Z_{T,\mrm{max}} \defeq \max_{\theta \in \Theta} Z_{S,\theta}(T) \leq \exp\lf(-\f{\theta_{1,\mrm{min}}}{1+
 \theta_{2,\mrm{max}}}  T\ri) <1,$$
and thus by \eqref{eq:h-bdd-lip-unif}, the density $\tilde \tauden_{\theta,T}$ is uniformly bounded in $t \in [0,T]$ and $\theta \in \Theta$, that is,
\begin{align}\label{eq:cond-den-max}
\|\tilde \tauden\|_{\infty,T,\mrm{max}} \defeq \sup_{\theta \in \Theta} \sup_{t\leq T} \tilde \tauden_{\theta,T}(t) \leq \|h\|_{\infty,1,\mrm{max}}/(1-Z_{T,\mrm{max}}).
\end{align}
Next we show that $\ln\lf(\tilde \tauden_{\theta,T}(\cdot)\ri)$ is Lipschitz continuous with Lipschitz constant uniformly bounded in $\theta \in \Theta$. First observe that $\tilde \tauden_{\theta,T}$ is Lipschitz continuous with Lipschitz constant satisfying
\begin{align}
    \label{eq:cond-den-lip}
    \|\tilde\tauden_{\theta,T}\|_{\mrm{Lip}} \leq  \|h_\theta\|_{\mrm{Lip}}\|Z_{S,\theta}\|_{\mrm{Lip}}/(1-Z_{S,\theta}(T)). 
\end{align}
Next, by \eqref{eq:h-bdd-lip}
 \begin{align*}
     \f{\differential{Z_{S,\theta}(t)}}{\differential{t}} \geq - (\theta_1/\theta_2) Z_{S,\theta}(t), \quad Z_{S,\theta}(0) =1. 
 \end{align*}
By the comparison principle for \acp{ODE}, we have $Z_{S,\theta}(t) \geq \myExp{-(\theta_1/\theta_2)  t}$. It follows that 
\begin{align*}
 Z_{\mrm{min}}\defeq \min_{\theta \in \Theta}\min_{t\leq T}Z_{S,\theta}(t) \geq  e^{-(\theta_{1,\mrm{max}}/\theta_{2,\mrm{min}}) T} >0\eqstop 
\end{align*}
Consequently, because  $h_\theta$ is increasing, 
\begin{align*}
    h_{\mrm{min}} \defeq \min_{\theta \in \Theta}\min_{t \leq T} h_\theta(Z_{S,\theta}(t)) \geq \min_{\theta \in \Theta} h_\theta(Z_{\mrm{min}}) \geq \f{\theta_{1, \mrm{min}} Z_{\mrm{min}}}{Z_{\mrm{min}}+\theta_{2,\mrm{max}}}  >0,  
\end{align*}
and
\begin{align}\label{eq:tauden-cond-min}
    \tilde \tauden_{\mrm{min}} \defeq \min_{\theta \in \Theta}\min_{t \leq T} \tilde \tauden_{\theta,T}(t) = \min_{\theta \in \Theta}\min_{t \leq T} \f{h_{\theta}(Z_{S,\theta}(t))}{1-Z_{S,\theta}(T)} \geq h_{\mrm{min}}/(1-Z_{\mrm{min}}) > 0.
\end{align}
Therefore by \eqref{eq:cond-den-lip}, \eqref{eq:h-bdd-lip-unif}, \eqref{eq:Z-lip-unif}, and the Mean Value Theorem, we have for some $\alpha \in [0,1]$
\begin{align*}
    |\ln\lf(\tilde \tauden_{\theta,T}(t)\ri) - \ln\lf(\tilde \tauden_{\theta,T}(t')\ri)| & = \f{1}{\alpha \tilde \tauden_{\theta,T}(t)+(1-\alpha)\tilde \tauden_{\theta,T}(t')}|\tilde \tauden_{\theta,T}(t) -\tilde \tauden_{\theta,T}(t') | 
    \leq \f{\|h\|_{\mrm{Lip,max}}    \|Z\|_{\mrm{Lip, max}}}{\tilde \tauden_{\mrm{min}} (1-Z_{T,\mrm{max}})}|t - t'| \eqstop 
\end{align*}
It follows from \eqref{eq:log-like},   \eqref{eq:cond-den-max}, and the above inequality that
\begin{align*}
\sup_{\theta \in \Theta} \Big|\ell_n(\theta)  &+ \mrm{CE}(\tilde\taudist_{\theta^*,T}, \tilde\taudist_{\theta,T})\Big| = \sup_{\theta \in \Theta}\lf| \tilde \Xi^{(n)}(\ln \tilde \tauden_{\theta,T}) - \tilde \taudist_{\theta^*,T}(\ln \tilde \tauden_{\theta,T})\ri|\\ 
\leq&\ \sup_{\theta \in \Theta} \|\ln \tilde \tauden_{\theta,T}\|_{\mrm{BL, T}}\ d_{\mrm{BL}}(\tilde \Xi^{(n)}, \tilde\taudist_{\theta^*,T}) \\
\leq& \  \lf(\lf|\ln \tilde \tauden_{\mrm{min}}\ri| \vee\lf|\ln \|\tilde \tauden\|_{\infty,T,\mrm{max}}\ri|+ \f{\|h\|_{\mrm{Lip,max}} \|h\|_{\infty,1,\mrm{max}}}{\tilde \tauden_{\mrm{min}} (1-Z_{T,\mrm{max}})}\ri) d_{\mrm{BL}}(\tilde \Xi^{(n)}, \tilde\taudist_{\theta^*,T}),
\end{align*}
and the result follows from \Cref{rem:emp-meas-conv}.
\end{proof}

Before proceeding to the proof of the consistency of $\hat \theta_n$,  we need the following result showing that the probability measure $\tilde\taudist_{\theta,T}$ is identifiable in $\theta$.

\begin{myLemma}
    \label{lem:cond-dist-para-ident}
 The mapping $\theta \in \Theta \mapsto \tilde{\taudist}_{\theta,T} \in \mathcal{P}([0,T])$ is injective.
\end{myLemma}

\begin{proof}[Proof of \Cref{lem:cond-dist-para-ident}]
    Suppose that $\tilde{\taudist}_{\theta,T} = \tilde{\taudist}_{\theta',T}$ for $\theta =(\theta_1,\theta_2),\ \theta' = (\theta'_1,\theta'_2) \in \Theta$. Because of the continuity of the density function $ \tilde{\tauden}_{\theta,T}(\cdot)$ it follows that for all $t \in [0,T]$
   $$\tilde{\tauden}_{\theta,T}(t) \equiv \f{h_\theta(Z_{S,\theta}(t))}{1-Z_{S,\theta}(T)} = \f{h_{\theta'}(Z_{S,\theta'}(t))}{1-Z_{S,\theta'}(T)} \equiv \tilde{\tauden}_{\theta',T}(t).$$
  Differentiating the above equation with respect to $t$ and using the fact that $\differential{Z_{S,\theta}(t)}/\differential{t} = - h_\theta(Z_{S,\theta}(t))$, we see that
 \begin{align}
     \label{eq:h-deriv-eql}
      h'_\theta(Z_{S,\theta}(t)) = h'_{\theta'}(Z_{S,\theta'}(t)).
 \end{align}
  Using the initial condition that $Z_{S,\theta}(0) = Z_{S,\theta'}(0)=1$, and observing $h'(y) = \theta_1 \theta_2 / (\theta_2 + y)^2$, we get
  \begin{align}
      \label{eq:theta-1}
      \frac{\theta_1 \theta_2}{(\theta_2 + 1)^2} =  \frac{\theta'_1 \theta'_2}{(\theta'_2 + 1)^2} \eqstop 
  \end{align}
Further differentiating \eqref{eq:h-deriv-eql} gives
 \begin{align*}
      h''_\theta(Z_{S,\theta}(t))h_\theta(Z_{S,\theta}(t)) = h''_{\theta'}(Z_{S,\theta'}(t))h_{\theta'}(Z_{S,\theta'}(t)).
 \end{align*}
 Again,  using  $Z_{S,\theta}(0) = Z_{S,\theta'}(0)=1$ and noting that $h''_\theta(y) = -2 \theta_1 \theta_2 / (\theta_2 + y)^3$, we get
 \begin{align}\label{eq:theta-2}
 \frac{\theta^2_1 \theta_2}{(\theta_2 + 1)^4} = \frac{(\theta'_1)^2 \theta'_2}{(\theta'_2 + 1)^4} \eqstop 
 \end{align}
 Solving \eqref{eq:theta-1} and \eqref{eq:theta-2} shows that $\theta \equiv(\theta_1,\theta_2) = (\theta'_1,\theta'_2) \equiv \theta'.$
\end{proof}

\np
We are now ready to prove consistency of $\hat \theta_n$.
\begin{proof}[Proof of \Cref{th:const}]
We first note that the mapping $\theta \in \Theta \mapsto \mrm{CE}(\tilde\taudist_{\theta^*,T}, \tilde\taudist_{\theta,T})$ is Lipschitz continuous. To verify this, note that for $\theta,\theta'\in\Theta \subset (0, \infty)^2$ ,
\begin{align*}
\begin{aligned}
\sup_{y \in [0,1]}|h_\theta(y)-h_{\theta'}(y)|
&\leq \sup_{y \in [0,1]} \lf(
\frac{y}{\theta_2+y}|\theta_1-\theta_1'|
+
\frac{\theta_1'y}{(\theta_2+y)(\theta_2'+y)}
|\theta_2-\theta_2'|\ri) \\
& \leq\
\frac{1}{1+\theta_{2,\mathrm{min}}}|\theta_1-\theta_1'|
+
\frac{\theta_{1,\mathrm{max}}}{4\theta_{2,\mathrm{min}}}
|\theta_2-\theta_2'| \ \leq\ L_{h,\mathrm{Lip}}\|\theta-\theta'\|_1,
\end{aligned}
\end{align*}
where
$$
L_{h,\mathrm{Lip}}
=
\max\left\{
\frac{1}{1+\theta_{2,\mathrm{min}}},
\frac{\theta_{1,\mathrm{max}}}{4\theta_{2,\mathrm{min}}}
\right\}.    
$$
and the second inequality uses
$$
\frac{y}{(\theta_2+y)(\theta_2'+y)}
\leq
\frac{y}{(\theta_{2,\mathrm{min}}+y)^2}
\leq \frac{1}{4\theta_{2,\mathrm{min}}}.
$$
Thus,
\begin{align*}
    |Z_{S,\theta}(t) -  Z_{S,\theta'}(t)| =&\ \lf|\int_0^t \big(h_{\theta}(Z_{S,\theta}(s))-h_{\theta'}(Z_{S,\theta'}(s))\big)\differential{s}\ri| \\
    \leq &\ 
    L_{h, \mrm{Lip}}T\|\theta-\theta'\|_1 + \|h\|_{\mrm{Lip,max}} \int_0^t |Z_{S,\theta}(s) -  Z_{S,\theta'}(s)| \differential{s},
\end{align*}
where $\|h\|_{\mrm{Lip,max}}$ is as in \eqref{eq:h-bdd-lip-unif}. Gr\"onwall's inequality now shows that
\begin{align*}
    \sup_{t\leq T}|Z_{S,\theta}(t) -  Z_{S,\theta'}(t)| \leq e^{\|h\|_{\mrm{Lip,max}}T}L_{h, \mrm{Lip}}T\|\theta-\theta'\|_1\eqstop 
\end{align*}    
A straightforward calculation now shows that the map $\theta \in \Theta \mapsto \tilde \tauden_\theta(t)$ is Lipschitz continuous with the  Lipschitz constant 
$$L_{\tilde\tauden, \mrm{Lip}} \leq \f{(\|h\|_{\infty,1,\mrm{max}}+\|h\|_{\mrm{Lip,max}})e^{\|h\|_{\mrm{Lip,max}}T}L_{h, \mrm{Lip}}T+ L_{h, \mrm{Lip}}}{(1-Z_{T,\mrm{max}})^2}\eqcomma $$
bounded in $t \in [0,T]$. Hence, it follows that
\begin{align*}
|\mrm{CE}(\tilde\taudist_{\theta^*,T}, \tilde\taudist_{\theta,T}) - \mrm{CE}(\tilde\taudist_{\theta^*,T}, \tilde\taudist_{\theta',T})|
\leq & \int_0^T\lf|\ln (\tilde \tauden_{\theta,T}(t)) -\ln (\tilde \tauden_{\theta',T}(t))\ri| \tilde \tauden_{\theta^*,T}(t) \differential{t} 
\leq
\lf(L_{\tilde\tauden, \mrm{Lip}} \big/ \tilde \tauden_{\mrm{min}}\ri) \|\theta-\theta'\|_1,
\end{align*}
with $\tilde \tauden_{\mrm{min}}$ as in \eqref{eq:tauden-cond-min}.
Next, by the definition of $\hat \theta_n$, for any $\theta \in \Theta$
\begin{align}
    \label{eq:like-max}
 \ell_n(\hat \theta_n) \equiv \f{1}{K_n}\ln \like(\hat \theta_n|\tilde{\bm{\tau}}^{(n)}_{1:K_n}) \geq  \ell_n(\theta) \equiv \f{1}{K_n} \ln \like( \theta|\tilde{\bm{\tau}}^{(n)}_{1:K_n}).   
\end{align}
Since the sequence of random variables, $\{\hat \theta_n\}$, takes values in a compact set $\Theta$, it is tight. Let $\tilde\theta$ be a limit point of this sequence, with the subsequence $\hat \theta_{n_j} \stackrel{j\rt \infty} \LRT  \tilde\theta$. Note that, at this stage, we do not \emph{a priori} know whether $\tilde\theta$ is deterministic.  We will actually show that $\tilde\theta \equiv \theta^*$, $\PP_{\theta^*}$-a.s. Thus, the limit point  is independent of the subsequence (and obviously deterministic); hence, $\hat \theta^{(n)}\stackrel{\PP_{\theta^*}} \Rt \theta^*$ along the full sequence.  We now work toward establishing this.

By the continuity of the mapping $\theta \in \Theta \mapsto  \mrm{CE}(\tilde\taudist_{\theta^*,T}, \tilde\taudist_{\theta,T})$, under $\PP_{\theta^*}$,
$$\mrm{CE}(\tilde\taudist_{\theta^*,T}, \tilde\taudist_{\hat \theta_{n_j},T}) \stackrel{j \rt \infty}\LRT \mrm{CE}(\tilde\taudist_{\theta^*,T}, \tilde\taudist_{\tilde\theta,T}).$$
Next write
$$\ell_n(\hat \theta_n) = \ell_n(\hat \theta_n)+ \mrm{CE}(\tilde\taudist_{\theta^*,T}, \tilde\taudist_{\hat \theta_{n},T}) - \mrm{CE}(\tilde\taudist_{\theta^*,T}, \tilde\taudist_{\hat \theta_{n},T}),$$
and observe by \Cref{th:conv-log-like}, 
\begin{align*}
    |\ell_n(\hat \theta_n)+ \mrm{CE}(\tilde\taudist_{\theta^*,T}, \tilde\taudist_{\hat \theta_{n},T})| &\ \leq \sup_{\theta \in \Theta}\lf|\ell_n( \theta) + \mrm{CE}(\tilde\taudist_{\theta^*,T}, \tilde\taudist_{\theta,T})\ri| 
    \nrt 0, \quad \text{ in } L^1(\PP_{\theta^*})\eqstop 
\end{align*}
Hence, as $j\rt \infty$,  $\ell_n(\hat \theta_n) \RT - \mrm{CE}(\tilde\taudist_{\theta^*,T}, \tilde\taudist_{\tilde\theta,T})$, and taking limit in \eqref{eq:like-max} along the subsequence $\{n_j\}$ gives 
\begin{align*}
 \mrm{CE}(\tilde\taudist_{\theta^*,T}, \tilde\taudist_{\tilde\theta,T}) \leq  \mrm{CE}(\tilde\taudist_{\theta^*,T}, \tilde\taudist_{\theta,T}), \quad \PP_{\theta^*}\text{-a.s.} \quad \text{ for any } \theta \in \Theta.  
\end{align*}
Choosing $\theta = \theta^*$ in the above equation shows that 
$$\mrm{RE}(\tilde\taudist_{\theta^*,T}\| \tilde\taudist_{\tilde\theta,T}) =0, \quad \PP_{\theta^*}\text{-a.s.}\eqcomma $$
that is $\tilde\taudist_{\theta^*,T} = \tilde\taudist_{\tilde\theta,T}$, and  by \Cref{lem:cond-dist-para-ident}, $\tilde\theta = \theta^*$,  $\PP_{\theta^*}$-a.s.

\end{proof}


     \begin{figure}[t]
        \centering
        \includegraphics[width=0.49\textwidth]{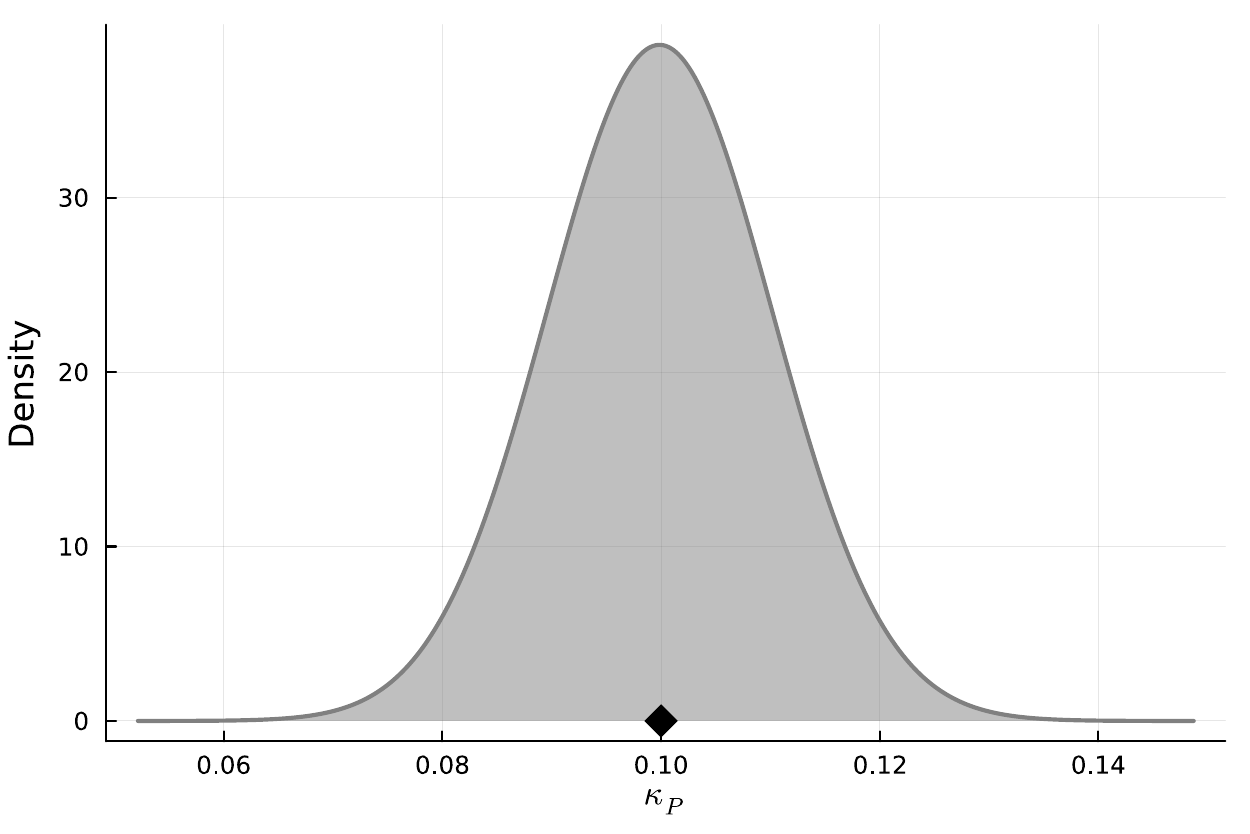}
        \includegraphics[width=0.49\textwidth]{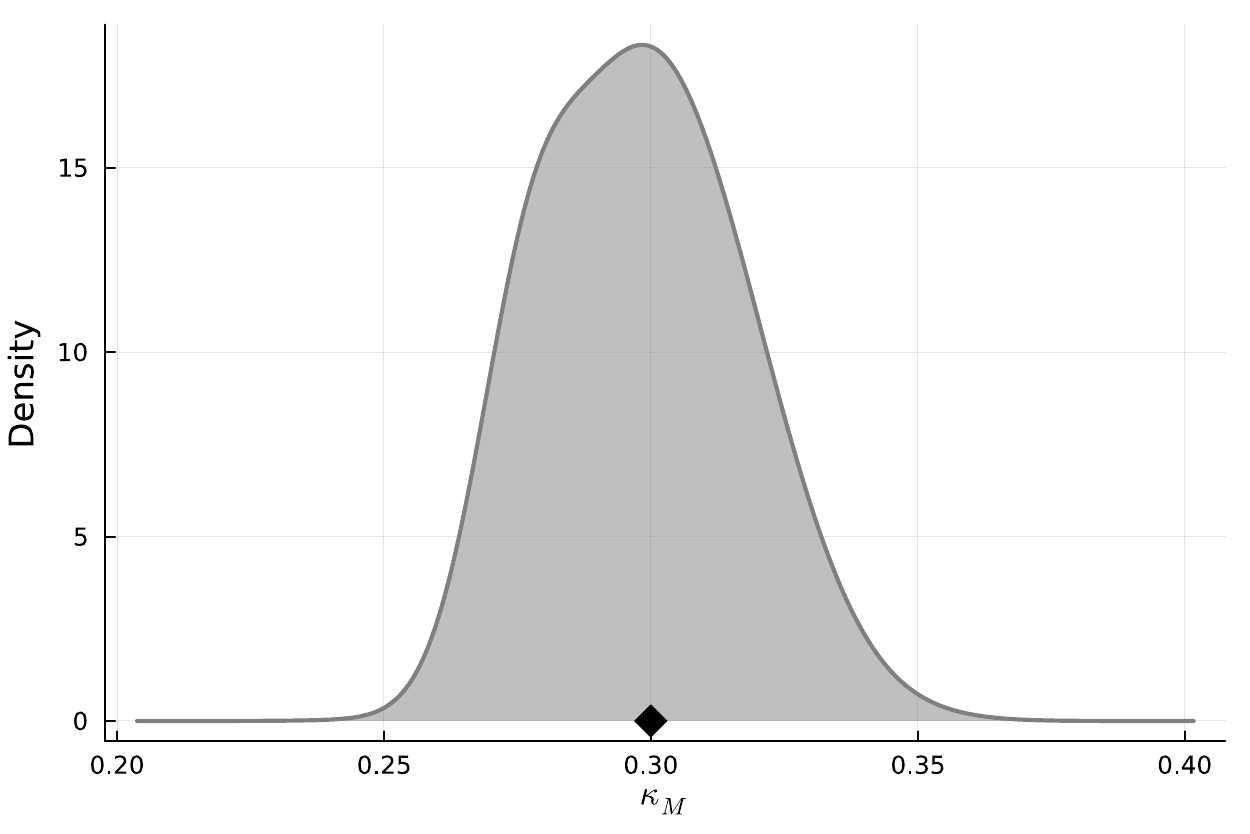}
        \caption{Densities of the \acp{MLE} $\kappa_P$ and $\kappa_M$ in the standard \ac{MM} kinetic reaction network considered in \Cref{example:standard_MM_inference} obtained by using \ac{KDE} on $5000$ \acp{MLE} each obtained from a fresh random sample of size $10^3$ of product formation times. (Left) Density of $\kappa_P$. (Right) Density of $\kappa_M$. True parameter values used in the simulation are $n = 10^6, J=10, \kappa_1 = 1, \kappa_{-1} = 0.2, \kappa_P = 0.1$, and $T=3.0$ with the true \ac{MM} constant being $\kappa_M = 0.3$. 
}
        \label{fig:MM_IPS_mle}
    \end{figure}

\begin{figure}[tbh]
        \centering
        \includegraphics[width=0.49\textwidth]{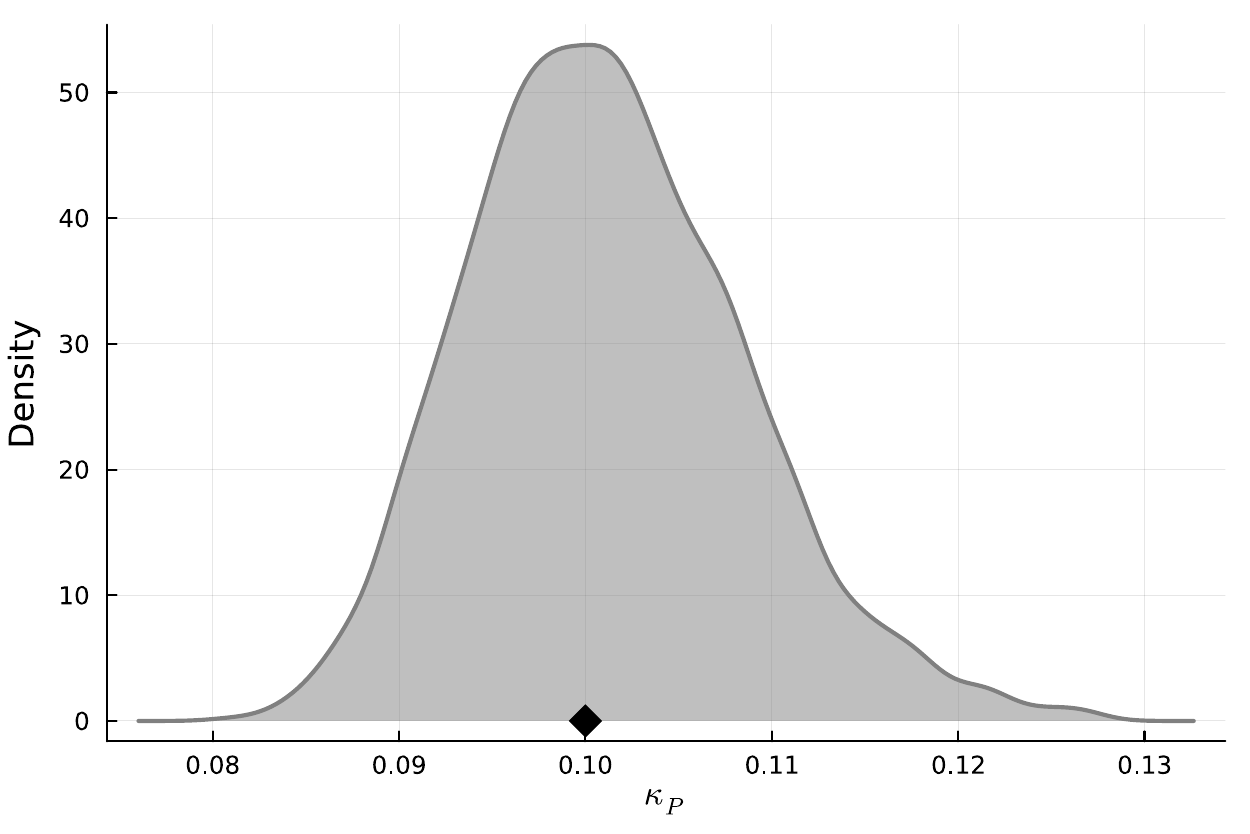}
        \includegraphics[width=0.49\textwidth]{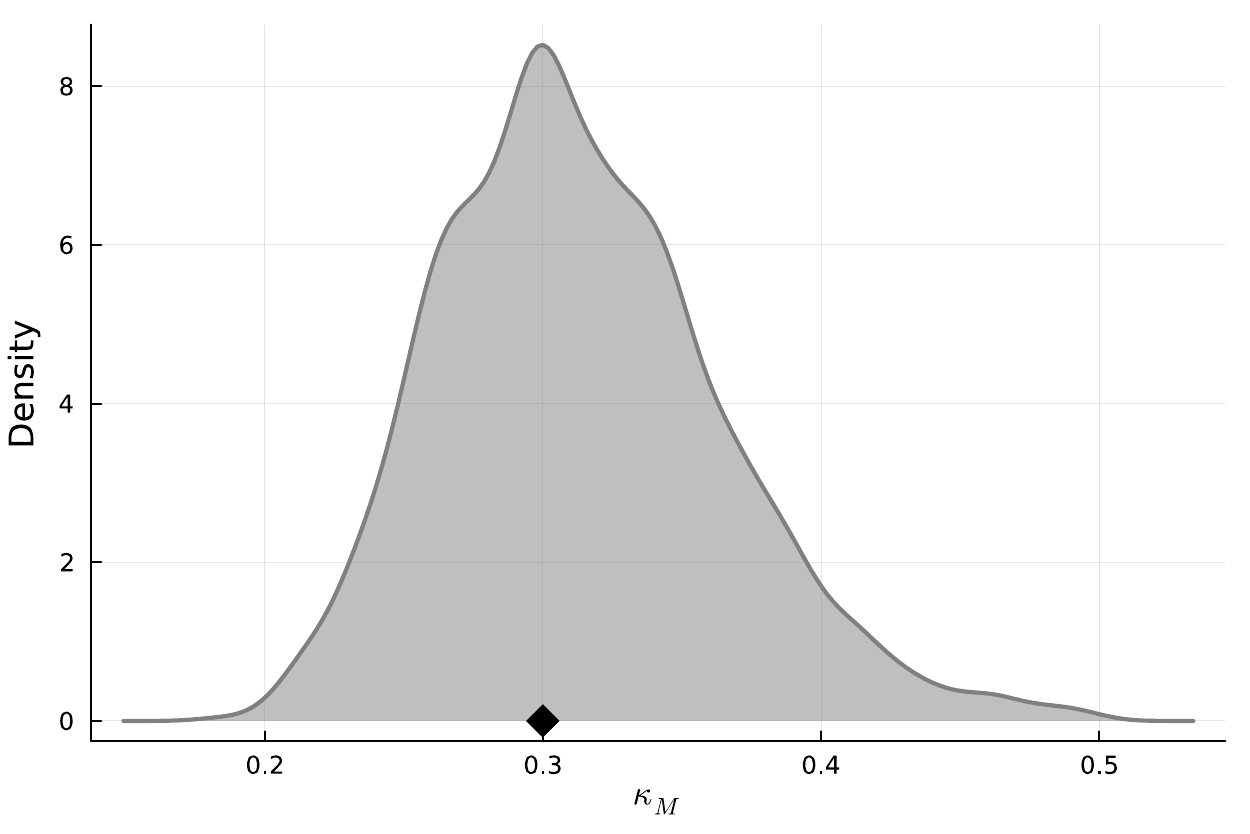}
        \caption{Posterior densities of the parameters $\kappa_P$ and $\kappa_M$ in the standard \ac{MM} kinetic reaction network considered in \Cref{example:standard_MM_inference}. (Left) Posterior density of $\kappa_P$. (Right) Posterior density of $\kappa_M$. True parameter values used in the simulation are $n = 10^6, J=10, \kappa_1 = 1, \kappa_{-1} = 0.2, \kappa_P = 0.1$, and $T=4.0$. Therefore, the true \ac{MM} constant is $\kappa_M = 0.3$. 
        }
        \label{fig:MM_IPS_posterior_density}
    \end{figure}
\begin{figure}[tbh]
        \centering
        \includegraphics[width=0.49\textwidth]{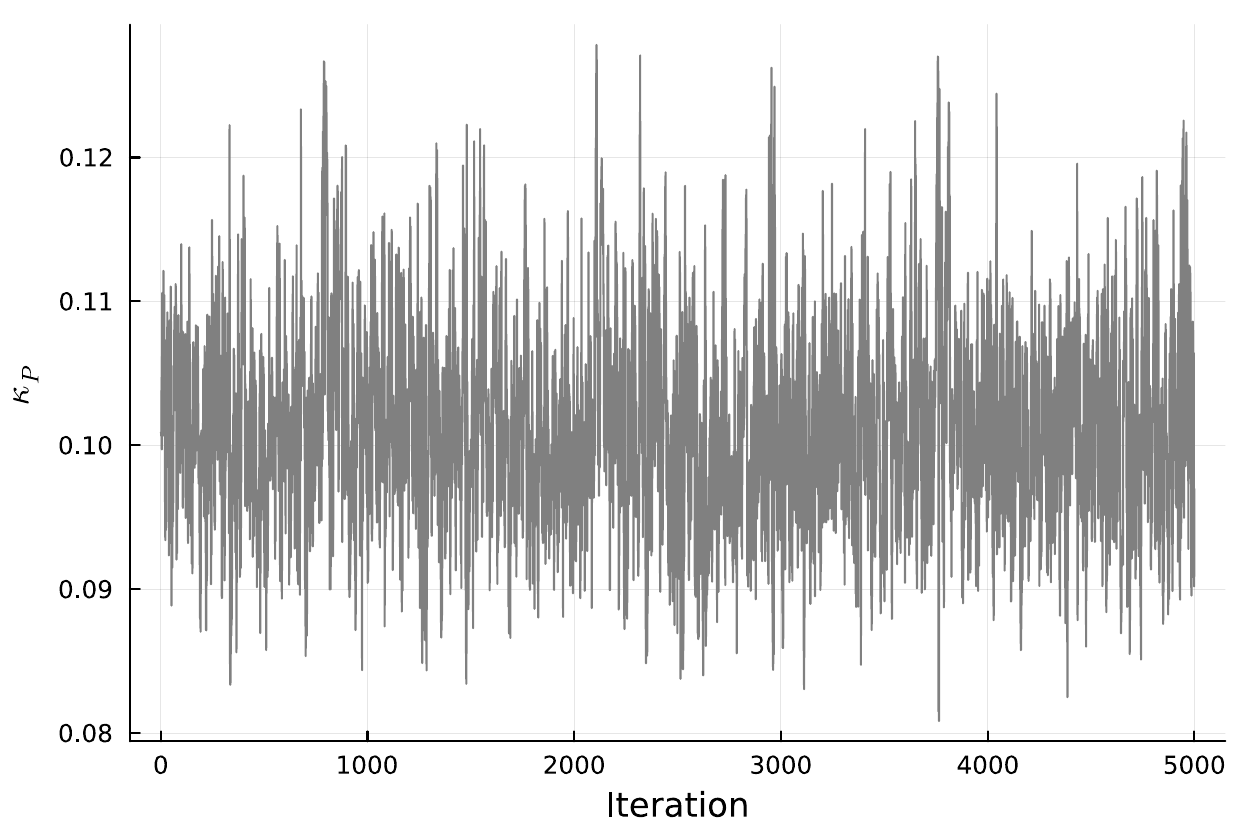}
        \includegraphics[width=0.49\textwidth]{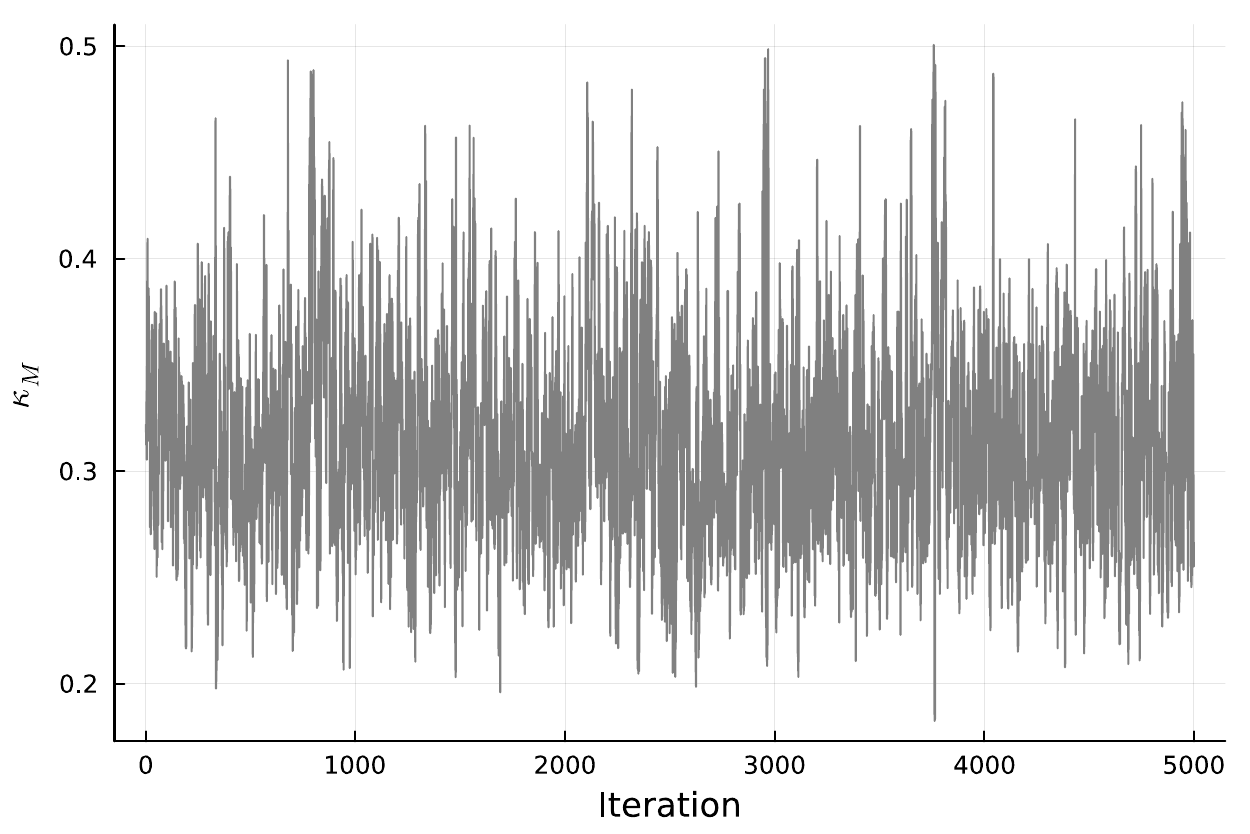}
        \caption{Trace plots of the parameters $\kappa_P$ and $\kappa_M$ in the standard \ac{MM} kinetic reaction network considered in \Cref{example:standard_MM_inference}. (Left) Trace plot of $\kappa_P$. (Right) Trace plot of $\kappa_M$. True parameter values used in the simulation are $n = 10^6, J=10, \kappa_1 = 1, \kappa_{-1} = 0.2, \kappa_P = 0.1$, and $T=4.0$. Therefore, the true \ac{MM} constant is $\kappa_M = 0.3$. 
        }
        \label{fig:MM_IPS_trace_plots}
    \end{figure}

    \begin{myExample}[Numerical results for the standard \ac{MM} model]\label{example:standard_MM_inference}
        Consider the simplest possible form of the \ac{MM} kinetic reaction in \eqref{eq:mm_standard_det}, which is \eqref{eq:mm_det} with $r=1$. We are presenting this case because this is the most commonly used form of \ac{MM} enzyme kinetic reaction, and we hope that the inference method will be readily used by applied mathematicians and statisticians. In this case, the function $h_\theta$ is given by $h_\theta(y) = \conscon\kappa_Py/(\kappa_M+y)$, where $\kappa_M = (\kappa_{-1} + \kappa_{P})/\kappa_{1}$. 
         In practice, the most important parameter is the \ac{MM} constant $\kappa_M$. However, for the purpose of illustration, we will estimate both $\kappa_M$ and $\kappa_P$, that is, we set $\theta = (\kappa_M, \kappa_P)$. We assume that we know the conservation constant $\conscon$. 

        In \Cref{fig:MM_IPS_mle}, we show a numerical example of the \ac{MLE} of $\theta$. We take  the initial amount of $S$ to be large so that the \ac{FLLN} is an accurate approximation. In the simulation, we take $n=10^6$. The standard frequentist approach is to report a point estimate along with a confidence interval. Here, we take a slightly different approach to uncertainty quantification. We first take a random sample of product formation times of size $10^3$ (approximately, $10^{-3}$-th of all product formation times if we were to observe till all $S$ molecules were converted into $P$ molecules). Then, we calculate the \ac{MLE} by maximizing the likelihood function. We then repeat this process $5000$ times, each time with a fresh random sample of product formation times of size $10^3$, to generate $5000$ different point estimates. We then use the \ac{KDE} methodology to construct a density of the obtained point estimates, which we show in \Cref{fig:MM_IPS_mle}. As we can see, the \acp{MLE} are very accurate.

        In addition to the \acp{MLE}, we also implement a \ac{HMC} method with  $\kappa_P \sim \sys{Uniform}(0, 0.5)$ and $ \kappa_M \mid \kappa_P \sim \sys{Uniform}(\kappa_P, 1)$ as prior distributions. We take a single random sample of product formation times of size $10^3$ (again, roughly $10^{-3}$-th of all product formation times) to evaluate the likelihood function, and run a single \ac{HMC} chain. We choose the standard values for the tuning parameters, \ie, 1000 burn-in, and then 5000 sample, 0.65 target acceptance probability. \Cref{fig:MM_IPS_posterior_density}  shows that the method is able to identify the true parameter values with very high accuracy. The trace plots in \Cref{fig:MM_IPS_trace_plots} show that the \ac{HMC} chain mixes well. Additional summary statistics and diagnostics are provided in \Cref{table:HMC_summary_statistics}. 

        The estimation tasks  are performed in the Julia programming language \cite{Bezanson2017Julia}. The programming scripts used in this example are available at the GitHub repository \cite{MM_inference_github_repo}.

    \end{myExample}


   
\begin{table}[ht]
\centering
\begin{tabular}{lrrrrrrr}
\hline
Parameter & Mean & SD & 2.5\% & 25\% & 50\% & 75\% & 97.5\% \\
\hline
$\kappa_P$ & 0.1011 & 0.0074 & 0.0883 & 0.0960 & 0.1006 & 0.1059 & 0.1172 \\
$\kappa_M$ & 0.3135 & 0.0506 & 0.2279 & 0.2769 & 0.3082 & 0.3447 & 0.4263 \\
\hline
\end{tabular}
\vspace{10pt}
\centering
\vspace{10pt}
\begin{tabular}{lrrrrrrr}
\hline
Parameter &  MCSE & ESS bulk & ESS tail & $\hat{R}$ & ESS/s \\
\hline
$\kappa_P$ &  0.0002 & 985.6747 & 696.2071 & 1.0027 & 36.6804 \\
$\kappa_M$ & 0.0018 & 953.5874 & 666.5672 & 1.0027 & 35.4863 \\
\hline
\end{tabular}
\caption{Summary statistics for the posterior samples and diagnostics for the \ac{HMC} chain.}
\label{table:HMC_summary_statistics}
\end{table}


\section{Concluding remarks}
\label{sec:conclusion}
In this paper, we present a stochastic model for a multistage Michaelis--Menten enzyme kinetic reaction network using a system of \acp{SDE} with jumps. Under a suitable scaling regime, we first establish a stochastic averaging principle yielding a reduced model for the product-substrate dynamics. For the purpose of statistical inference, we then construct a novel \ac{IPS} that approximates the product-substrate process at the particle level. We then construct an estimator based on a  product-form approximate likelihood requiring only a random sample of product formation times. This approach does not need access to the system states, and we rigorously prove the consistency of the estimator.

As mentioned at the beginning of \Cref{sec:IPS}, our \ac{IPS}-based inference approach is inspired by the \ac{DSA} methodology from infectious disease epidemiology 
\cite{DiLauro2022Interface,KhudaBukhsh2020Focus,KhudaBukhsh2023JTB,KhudaBukhsh2024Howto,Rempala2026DSA,Roy2026MarginalDSA}. Although the \ac{DSA} methodology has been applied either directly or after necessary modification to various data types  (\eg, a single snapshot of mass-testing data in a prison in Ohio \cite{KhudaBukhsh2023Prison}, individual incidence and hospitalization data for Ebola victims in the  Democratic Republic of Congo \cite{Vossler2022Ebola}, extensions to a marginal likelihood approach in \cite{Roy2026MarginalDSA}),  the most common use case of the method  is when the  likelihood function is based on a random sample of infection times (and recovery times, if available)  of individuals in a large population. The ability to work with a small amount of data has been the method's biggest selling point during and since the COVID-19 pandemic.  However, its applications and derivations are still largely heuristic and nonrigorous. In particular, a rigorous mathematical analysis of the  statistical and the asymptotic properties (\eg,  consistency) of the \ac{DSA}-based estimators  remains unavailable even in the context of epidemiology. For the \ac{MM} enzyme kinetic reaction network, our contribution in this paper is to provide a rigorous justification for the approximate product-form likelihood function using propagation of chaos, and then prove the consistency of the estimators in \Cref{th:const} along with several other asymptotic results and non-asymptotic concentration inequalities. 

The stochastic averaging principle in \Cref{thm:FLLN} corresponds to the \ac{sQSSA}. Several other types of \acp{QSSA} are available in the literature for \ac{MM} enzyme kinetic reaction networks, \eg, the \ac{tQSSA} and the \ac{rQSSA} (see \cite{Kang2019QSSA}). The validity of and the fluctuations around the \ac{tQSSA} have recently been rigorously analysed in \cite{Ganguly2025tQSSA}. The (scaling) parameter regimes for the validity of these other types of \acp{QSSA} are, in general, different from those assumed for the \ac{sQSSA} and could also require additional timescale separation.  As such, which particular \ac{QSSA} is applicable in a given situation is entirely decided by the validity of the underlying scaling assumption. Therefore, the asymptotic properties of our estimators are tied to the \ac{sQSSA}, and may not extend immediately to other types of \acp{QSSA} such as the \ac{tQSSA}. Nevertheless, one can see that a similar \ac{IPS} could be constructed under the \ac{tQSSA} and a similar inference method could be developed. A thorough analysis of this new inference method in other scaling regimes is beyond the scope of the present paper and will be pursued in the future. 










\appendix


\section{Additional derivations}
\label{sec:additional_derivations}

\begin{myLemma} \label{lem:conv-int}
Consider the sequence $\{(\phi^{(n)}, \nu^{(n)}) : n\ge 1\} \subset D([0,T],\R^d)\times \mathcal{M}(\mathbb{F}\times [0,T])$, where $\mathbb{F}$ is a finite set. Assume that $\phi^{(n)} \to \phi^{(\infty)} \in C([0,T], \R^d)$ in uniform metric and  $\nu^{(n)} \to  \nu^{(\infty)}$ in $\mathcal{M}(\mathbb{F}\times [0,T])$. Let $h: \R^d \times \mathbb{F}\times [0,T] \rt \R$  be a continuous function that is Lipschitz in the first argument, that is, 
\begin{align} \label{eq:h-lip}
|h(u,y,s)-h(u',y,s)| \leq L_h\|u-u'\|, \quad u, u' \in \R^d, y \in \mathbb F, s \in [0,T],    
\end{align}
for some constant $L_h \equiv L_h(T) >0$ independent of $y$ and $s$. Then, 
\begin{align*}
\int_{\mathbb{F}\times [0,T]} h(\phi^{(n)}(s), y, s)\nu^{(n)}(\differential{y}\times \differential{s})&\ \equiv \sum_{y  \in \mathbb{F} }\int_0^T h(\phi^{(n)}(s),y,s) \nu^{(n)}(y \times  \differential{s}) \\
&\ 
\stackrel{n\rt \infty}\Rt \int_{\mathbb{F}\times [0,T]} h(\phi^{(\infty)}(s), y, s)\nu^{(\infty)}(\differential{y}\times \differential{s}) \eqstop 
\end{align*}
    
\end{myLemma}

\begin{proof}
    The proof follows by first writing the integrand as
$$h(\phi^{(n)}(s),y,s) = h(\phi^{(\infty)}(s),y,s) + \lf(h(\phi^{(n)}(s),y,s) - h(\phi^{(\infty)}(s),y,s)\ri)\eqcomma $$  
and then using Lipschitz continuity of $h$ in \eqref{eq:h-lip},  uniform convergence of $\phi^{(n)} \stackrel{n\rt \infty}\Rt \phi^{(\infty)} $ (i.e.,  $\sup_{t\leq T} \|\phi^{(n)}(s) - \phi^{(\infty)}(s)\| \stackrel{n\rt \infty}\Rt 0$) and the  (weak) convergence of $\nu^{(n)}$ to $\nu^{(\infty)}$. Notice that since $\mathbb F$ is finite and $[0,T]$ is compact, $L_h$ in \eqref{eq:h-lip} can be chosen independently of $y$ and $s$  whenever for each \(y\in\mathbb F\) the Lipschitz constants in the first argument are locally bounded as functions of \(s\in[0,T]\).  
\end{proof}

\begin{myLemma}\label{lem:int-eq-tight}
For each $n\geq 0$, let $U^{(n)}, A^{(n)}$ and $B^{(n)}$ be stochastic processes  satisfying
     $$U^{(n)}(t)=A^{(n)}(t)+\int_0^t B^{(n)}(s)U^{(n)}(s) \differential{s}.$$
Assume that the  $\{\sup_{t\leq T}|A^{(n)}(t)|\}$ and $\{\sup_{t\leq T}|B^{(n)}(t)|\}$ are tight in $\setOfPositiveReals$ and $\{A^{(n)}\}$ is tight in $D([0,T],\R)$. Then $\{U^{(n)}\}$ is tight in $D([0,T],\R)$. If $\{A^{(n)}\}$ is $C$-tight in $D([0,T],\R)$, then so is $\{U^{(n)}\}$. 
\end{myLemma}
The idea of the proof is as follows: Gr\"onwall's inequality shows that $\{\sup_{t\leq T}|U^{(n)}(t)|\}$ is tight. It is now easy to show that for each $\delta$, $\modu(U^{(n)},T,\delta) \defeq \sup\{|U^{(n)}(t_1)-U^{(n)}(t_2)|: t_1,t_2\in [0,T], |t_2-t_1|\leq \delta\}$, the modulus of continuity of $U^{(n)}$, is tight. More details can be found in \cite[Corollary A.1]{Ganguly2025tQSSA}.

    \section{Stationary distribution}
    \label{sec:stationary_distribution}

    Consider a \ac{CTMC} $(S(t))_{t\ge 0}$ on the state space 
    \begin{align*}
        \mathcal{S}\defeq \{ (s_1, s_2, \ldots, s_{r+1}) \in \{0, 1, 2, \ldots, M\}^{r+1} \mid  \sum_{i=1}^{r+1} s_i  =  M\}\eqcomma 
    \end{align*}  
    with the generator $Q$ given by 
    \begin{align}
        \begin{aligned}
            \label{eq:CTMC_generator}
            Qf(s) &{}\defeq  l_1s_{r+1} \left(f(s+e_1^{(r+1)}-e_{r+1}^{(r+1)}) - f(s)\right) + l_{-1}s_1 \left(f(s-e_1^{(r+1)}+e_{r+1}^{(r+1)}) -f(s) \right) \\
            &{} +  \sum_{i=2}^{r+1} l_i s_{i-1} \left(f(s - e_{i-1}^{(r+1)} + e_i^{(r+1)}) - f(s)\right) 
            + \sum_{i=2}^{r} l_{-i} s_{i} \left(f(s + e_{i-1}^{(r+1)} - e_i^{(r+1)}) - f(s)\right) \eqcomma 
        \end{aligned}
    \end{align}
    for functions  $f : \{0, 1, 2, \ldots, M\}^{r+1} \mapsto \setOfReals$, and positive constants $l_1, l_{-1}, \ldots, l_r, l_{-r}, l_{r+1}$, where $e_i^{(r+1)}$ is the $i$-th unit basis vector in $\setOfReals^{r+1}$, \ie, the vector all components of which are zeroes except the $i$-th component, which is one. In practice, it will be easier to treat $Q$ as a matrix whose elements $q_{s, s'}$, for $s\neq s'$, are the jump intensities given in \eqref{eq:CTMC_generator}, and $q_{ss} = - \sum_{s' \neq s} q_{s, s'}$. The generator in \eqref{eq:CTMC_generator} can be rewritten as 
    \begin{align*}
        Qf(s) = \sum_{s' \in \mathcal{S}} q_{s, s'} \left(f(s') - f(s)\right) \eqstop 
    \end{align*}

    \begin{myLemma}
        The \ac{CTMC} $(S(t))_{t\ge 0}$ with generator $Q$ given in \eqref{eq:CTMC_generator} has a unique stationary distribution $\pi$ given by 
        \begin{align}
            \pi(s) & = \indicator{\mathcal{S}}{s} \frac{M!}{  \prod_{i=1}^{r+1} s_{i}!} \prod_{i=1}^{r+1} p_i^{s_{i}}  \eqcomma \label{eq:CTMC_stationary_distribution}
        \end{align}
        for $s = (s_1, s_2, \ldots, s_{r+1})$ where 
        \begin{align}
            \begin{aligned}
                p_1 &= \left(1+ a_1 + \sum_{i=2}^{r}  \frac{1}{\prod_{j=2}^{i} a_j} \right)^{-1}\eqcomma \ \
                p_{i} = \frac{p_1}{\prod_{j=2}^{i} a_j } \text{ for } i=2, 3, \ldots, r\eqcomma \text{ and } \ 
                 p_{r+1} 
                 = 1- \sum_{i=1}^{r}p_i \eqcomma \\
            \end{aligned}
            \label{eq:ctmc_stationary}
        \end{align}
        where the numbers $a_1, a_2, \ldots, a_{r+1}$ satisfy the following recursive relations
        \begin{align}
            \begin{aligned}
                a_{r} & = \frac{(l_{-r}+l_{r+1})}{l_r } \eqcomma \quad a_1  = \frac{l_{-1}}{l_1} + \frac{1}{a_2 a_3 \cdots a_r} \frac{l_{r+1}}{l_1}\eqcomma \quad  
                a_{i} 
                = \frac{(l_{-i} + l_{i+1})}{l_i } - \frac{l_{-(i+1)}}{a_{i+1} l_i}  \text{ for } i=2, 3,  \ldots, r-1\eqcomma \\
                  a_{r+1} & = \left(\prod_{i=1}^{r}a_i\right)^{-1}\eqstop 
            \end{aligned}
        \end{align}
        
        \label{lemma:stationary_distribution}
    \end{myLemma}
    \begin{proof}[Proof of \Cref{lemma:stationary_distribution}]
        It is straightforward to verify that the \ac{CTMC} $(S(t))_{t\ge 0}$ is irreducible, and positive recurrent \cite[Chapter~3]{Norris1997MarkovChains}. Therefore, it has a unique stationary distribution $\pi$ satisfying $\sum_{s'} \pi(s') q_{s', s} = 0$ for all $s \in \mathcal{S}$. Based on the literature on urn models, let us assume $\pi$ is of the form given by \eqref{eq:CTMC_stationary_distribution} for some $p_1, p_2, \ldots, p_{r+1}$. 
        Then, putting $s =  (M, 0, 0, \ldots, 0)$, $(0, M, 0, \ldots, 0), \ldots (0, 0, \ldots, M)$ in $\sum_{s'} \pi(s') q_{s', s} = 0$ gives us 
        \begin{align}
            \begin{aligned}
                p_r l_{r+1} + p_1 l_{-1} & = p_{r+1}  l_1 \eqcomma \\
                p_{r+1}l_1 + p_2 l_{-2} &= p_1 (l_{-1}+l_2) \eqcomma \\
                p_{i-1} l_i + p_{i+1}l_{-(i+1)} & = p_i (l_{-i} + l_{i+1}) \text{ for } i=2, 3, \ldots, r-1\eqcomma \\
                p_{r-1}l_r & = p_r (l_{-r} + l_{r+1}) \eqstop 
            \end{aligned}
            \label{eq:_pr_recursions}
        \end{align}
        Now, define $$a_1 \defeq p_{r+1}/p_1,\; a_2 \defeq p_{1}/p_{2},\; a_3 \defeq p_{2}/p_{3},\; \ldots, a_{r+1} \defeq p_{r}/p_{r+1}$$ so that $\prod_{i=1}^{r+1} a_i=1$. Then,  \eqref{eq:_pr_recursions} yields immediately $a_r = (l_{-r}+l_{r+1})/l_r$. Now, from \eqref{eq:_pr_recursions}, we get 
        \begin{align*}
            a_{i} = \frac{(l_{-i} + l_{i+1})}{l_i} - \frac{l_{-(i+1)}}{a_{i+1} l_i}  \text{ for } i=2, 3,  \ldots, r-1\eqcomma
        \end{align*}
        from which we can find $a_{r-1}$ by plugging in $a_r = (l_{-r}+l_{r+1})/l_r$. Continuing this process, we can find $a_{r-2}, a_{r-3}, \ldots, a_2$. We find $a_1$ from the first equation in \eqref{eq:_pr_recursions}
        \begin{align*}
            a_1 = \frac{l_{-1}}{l_1} + \frac{1}{a_2 a_3 \cdots a_r} \frac{l_{r+1}}{l_1} \eqstop 
        \end{align*}
        Finally, $a_{r+1} = (\prod_{i=1}^{r}a_i)^{-1}$. 
        
        Let us now fix $p_1$. Then, notice that $p_2 = p_1/a_2, p_3 = p_2/(a_3) = p_1/(a_2a_3), \ldots, p_{r} = p_{r-1}/a_r = p_1/(a_2\cdots a_r)$, and $p_{r+1} = a_1 p_1$. Since $\sum_{i=1}^{r+1} p_i = 1$, we must have 
        \begin{align*}
            p_1 = \left(1+ a_1 + \sum_{i=2}^{r}  \frac{1}{\prod_{j=2}^{i} a_j} \right)^{-1} \eqcomma
        \end{align*}
        which completes the proof.

    \end{proof}

\begin{acronym}[OWL-QN]
	\acro{ABM}{Agent-based Model}
	\acro{BA}{Barab\'asi--Albert}
    \acro{CDC}{Centers for Disease Control and Prevention}
	\acro{CDF}{Cumulative Distribution Function}
	\acro{CLT}{Central Limit Theorem}
	\acro{CM}{Configuration Model}
	\acro{CME}{Chemical Master Equation}
	\acro{CRM}{Conditional Random Measure}
	\acro{CRN}{Chemical Reaction Network}
	\acro{CTBN}{Continuous Time Bayesian Network}
	\acro{CTMC}{Continuous Time Markov Chain}
	\acro{DCT}{Dominated Convergence Theorem}
	\acro{DSA}{Dynamic Survival Analysis}
	\acro{DTMC}{Discrete Time Markov Chain}
	\acro{DRC}{Democratic Republic of Congo}
	\acro{ER}{Erd\"{o}s--R\'{e}nyi}
	\acro{ESI}{Enzyme-Substrate-Inhibitor}
	\acro{FCLT}{Functional Central Limit Theorem}
	\acro{FLLN}{Functional Law of Large Numbers}
	\acrodefplural{FLLN}[FLLNs]{Functional Laws of Large Numbers}
	\acro{FPT}{First Passage Time}
	\acro{GP}{Gaussian Process}
	\acrodefplural{GP}[GPs]{Gaussian Processes}
	\acro{HJB}{Hamilton–Jacobi–Bellman}
	\acro{HMC}{Hamiltonian Monte Carlo}
	\acro{iid}{independent and identically distributed}
	\acro{IPS}{Interacting Particle System}
        \acro{KDE}{Kernel Density Estimation}
	\acro{KL}{Kullback-Leibler}
	\acro{LDP}{Large Deviations Principle}
	\acro{LLN}{Law of Large Numbers}
	\acrodefplural{LLN}[LLNs]{Laws of Large Numbers}
	\acro{LNA}{Linear Noise Approximation}
	\acro{MAPK}{Mitogen-activated Protein Kinase}
	\acro{MCLT}{Martingale Central Limit Theorem}
	\acro{MCMC}{Markov Chain Monte Carlo}
	\acro{MFPT}{Mean First Passage Time}
	\acro{MGF}{Moment Generating Function}
	\acro{MLE}{Maximum Likelihood Estimate}
	\acro{MM}{Michaelis--Menten}
	\acro{MPI}{Message Passing Interface}
	\acro{MSE}{Mean Squared Error}
	\acro{ODE}{Ordinary Differential Equation}
	\acro{PDE}{Partial Differential Equation}
	\acro{PDF}{Probability Density Function}
	\acro{PGF}{Probability Generating Function}
	\acro{PMF}{Probability Mass Function}
	\acro{PPM}{Poisson Point Measure}
        \acro{PRM}{Poisson Random Measure}
	\acro{psd}{positive semi-definite}
	\acro{PT}{Poisson-type}
	\acro{QSSA}{Quasi-Steady State Approximation}
	\acro{rQSSA}{reverse QSSA}
	\acro{SD}{Standard Deviation}
        \acro{SDE}{Stochastic Differential Equation}
	\acro{SEIR}{Susceptible-Exposed-Infected-Recovered}
	\acro{SI}{Susceptible-Infected}
	\acro{SIR}{Susceptible-Infected-Recovered}
	\acro{SIS}{Susceptible-Infected-Susceptible}
	\acro{sQSSA}{standard QSSA}
	\acro{tQSSA}{total QSSA}
	\acro{TK}{Togashi--Kaneko}
	\acro{WS}{Watts-Strogatz}
	\acro{whp}{with high probability}
        \acro{BDG}{Burkholder--Davis--Gundy}
 
\end{acronym}





\section*{Acknowledgment}
The authors thank the Editorial Board and the anonymous reviewers for their constructive feedback, which improved the paper substantially.




\printcredits

\bibliographystyle{cas-model2-names}

\bibliography{references}



\end{document}